\documentclass[11pt]{article}
\usepackage{amsmath,amsthm,amsfonts,amssymb,latexsym}
\usepackage{amsmath}

\usepackage{a4wide}
\newtheorem{theorem}{Theorem}[section]

\newtheorem{lemma}[theorem]{Lemma}
\newtheorem{remark}[theorem]{Remark}
\newtheorem{proposition}[theorem]{Proposition}

\newtheorem{example}[theorem]{Example}

\def\l{\lambda}

\def\j{\mathfrak{j}}

\def\L{{\cal   L}}

\def\hb#1{\hbox{#1}}

\def\hb #1{\hbox{#1}}

\def\hb#1{\hbox{#1}}

\def\dim#1{\hb{dim}(#1)}

\def\L1#1{L^1(#1)}

\def\lef({\left(}
\def\rig){\right)}

\begin{document}
\title{ Cohomology of the Lie Superalgebra of Contact Vector Fields on
$\mathbb{R}^{1|1} $ and Deformations of the Superspace of Symbols }

\author{ Imed Basdouri $^1$, Mabrouk Ben Ammar $^1$,\\ Nizar Ben Fraj
$^2$, Maha Boujelben $^1$ and Kaouthar Kammoun $^1$\\
$^1$D\'epartement de Math\'ematiques, Facult\'e des
Sciences de Sfax,\\ Route de Soukra 3018 Sfax BP 802, Tunisie\\
E-mails: basdourimed@yahoo.fr, mabrouk.benammar@fss.rnu.tn,\\
Maha.Boujelben@fss.rnu.tn, lkkaouthar@yahoo.com\\[10pt]
$^2$Institut Sup\'{e}rieur de Sciences Appliqu\'{e}es et
Technologie, Sousse, Tunisie \\
E-mail:~benfraj\_nizar@yahoo.fr\\[10pt] }

\maketitle

\begin{abstract}
Following Feigin and Fuchs, we compute the first cohomology of the
Lie superalgebra $\mathcal{K}(1)$ of contact vector fields on the
(1,1)-dimensional real superspace with coefficients in the
superspace of linear differential operators acting on the
superspaces of weighted densities. We also compute the same, but
$\mathfrak{osp}(1|2)$-relative, cohomology. We explicitly give
1-cocycles spanning these cohomology. We classify generic formal
$\mathfrak{osp}(1|2)$-trivial deformations of the
$\mathcal{K}(1)$-module structure on the superspaces of symbols of
differential operators. We prove that any generic formal
$\mathfrak{osp}(1|2)$-trivial deformation of this
$\mathcal{K}(1)$-module is equivalent to a polynomial one of
degree $\leq4$. This work is the simplest superization of a result
by Bouarroudj [On $\mathfrak{sl}$(2)-relative cohomology of the
Lie algebra of vector fields and differential operators, J.
Nonlinear Math. Phys., no.1, (2007), 112--127]. Further
superizations correspond to $\mathfrak{osp}(N|2)$-relative
cohomology of the Lie superalgebras of contact vector fields on
$1|N$-dimensional superspace.
\end{abstract}

\maketitle

\section{Introduction}
For motivations, see Bouarroudj's paper \cite{b} of which this
work is the most natural superization, other possibilities being
cohomology of polynomial versions of various infinite dimensional
\lq\lq stringy" Lie superalgebras (for their list, see
\cite{gls}). This list contains several infinite series and
several exceptional superalgebras, but to consider cohomology
relative a \lq\lq middle" subsuperalgebra similar, in a sense, to
$\mathfrak{sl}(2)$ is only possible when such a subsuperalgebra
exists which only happens in a few cases. Here we consider the
simplest of such cases.

Let $\mathfrak{vect}(1)$ be the Lie algebra of polynomial vector
fields on $\mathbb{K}:=\mathbb{R}$ or $\mathbb{C}$. Consider the
1-parameter deformation of the $\mathfrak{vect}(1)$-action  on
$\mathbb{K}[x]$:
\begin{equation*}
L_{X\frac{d}{dx}}^\lambda(f)= Xf'+\lambda X'f,
\end{equation*}
where $X, f\in\mathbb{K}[x]$ and $X':=\frac{dX}{dx}$. This
deformation shows that on the level of Lie algebras (and similarly
below, for Lie superalgebras) it is natural to choose $\mathbb{C}$
as the ground field.

Denote by $\mathcal{F}_\lambda$ the $\mathfrak{vect}(1)$-module
structure on $\mathbb{K}[x]$ defined by $L^\lambda$ for a fixed
$\lambda$. Geometrically, ${\cal F}_\lambda=\left\{fdx^{\lambda}\mid
f\in \mathbb{K}[x]\right\}$ is the space of polynomial weighted
densities of weight $\lambda\in\mathbb{C}$. The space ${\cal
F}_\lambda$ coincides with the space of vector fields, functions and
differential 1-forms for $\lambda = -1,\, 0$ and $1$, respectively.

Denote by $\mathrm{D}_{\nu,\mu}:=\mathrm{Hom}_{\rm{diff}}({\cal
F}_\nu, {\cal F}_\mu)$ the $\mathfrak{vect}(1)$-module of linear
differential operators with the natural
$\mathfrak{vect}(1)$-action denoted $L_X^{\nu,\mu}(A)$. Each
module $\mathrm{D}_{\nu,\mu}$ has a natural filtration by the
order of differential operators; the graded module ${\cal
S}_{\nu,\mu}:=\mathrm{gr}\mathrm{D}_{\nu,\mu}$ is called the {\it
space of symbols}. The quotient-module
$\mathrm{D}^k_{\nu,\mu}/\mathrm{D}^{k-1}_{\nu,\mu}$ is isomorphic
to the module of weighted densities $\mathcal{F}_{\mu-\nu-k}$; the
isomorphism is provided by the principal symbol map $\sigma_{\rm
pr}$ defined by:
\begin{equation*}
A=\sum_{i=0}^ka_i(x)\left(\frac{\partial}{\partial
x}\right)^i\mapsto\sigma_{\rm pr}(A)=a_k(x)(dx)^{\mu-\nu-k},
\end{equation*}
(see, e.g.,\cite{gmo}). Therefore, as a
$\mathfrak{vect}(1)$-module, the space ${\cal S}_{\nu,\mu}$
depends only on the
difference $\beta=\mu-\nu$,
so that ${\cal S}_{\nu,\mu}$ can be written as ${\cal S}_{\beta}$,
and we have
\begin{equation*}
{\cal S}_{\beta} = \bigoplus_{k=0}^\infty \mathcal{F}_{\beta-k}
\end{equation*}
as $\mathfrak{vect}(1)$-modules. The space of symbols of order
$\leq n$ is
\begin{equation*}
{\cal S}_\beta^n:=\bigoplus_{k=0}^n{\cal F}_{\beta-k}.
\end{equation*}

In the last two decades, deformations of various types of
structures have assumed an ever increasing role in mathematics and
physics. For each such deformation problem a goal is to determine
if all related deformation obstructions vanish and many beautiful
techniques were developed to determine when this is so.
Deformations of Lie algebras with base and versal deformations
were already considered by Fialowski in 1986 \cite{f1}. In 1988,
Fialowski \cite{f2} further introduced deformations whose base is
a complete local algebra (the algebra is said to be {\it local} if
it has a unique maximal ideal). Also, in \cite{f2}, the notion of
miniversal (or formal versal) deformation was introduced in
general, and it was proved that under some cohomology
restrictions, a versal deformation exists. Later Fialowski and
Fuchs, using this framework, gave a construction for a versal
deformation.  Formal deformations of the
$\mathfrak{vect}(1)$-module ${\cal S}_\beta^n$ were studied in
\cite{aalo,bbdo}. Moreover, the formal deformations that become
trivial once the action is restricted to $\mathrm{\frak {sl}}(2)$
were completely described in \cite{bb}.

According to Nijenhuis-Richardson the space
$\mathrm{H}^1\left(\mathfrak{g};\mathrm{End}(V)\right)$ classifies
the infinitesimal deformations of a $\mathfrak{g}$-module $V$ and
the obstructions to integrability of a given infinitesimal
deformation of $V$ are elements of
$\mathrm{H}^2\left(\mathfrak{g};\mathrm{End}(V)\right)$. More
generally, if $\frak h$ is a subalgebra of $\frak g$, then the
$\frak h$-relative cohomology
$\mathrm{H}^1\left(\mathfrak{g},\frak h;\mathrm{End}(V)\right)$
measures the infinitesimal deformations that become trivial once
the action is restricted to $\frak h$ ($\frak h$-{\it trivial
deformations}), while the obstructions to extension of any $\frak
h$-trivial infinitesimal deformation to a formal one are related
to $\mathrm{H}^2\left(\mathfrak{g},\frak
h;\mathrm{End}(V)\right)$. Similarly, in the infinite dimensional
setting, the infinitesimal deformations of the
$\mathfrak{vect}(1)$-module ${\mathcal S}^n_\beta$ are classified,
from a certain point of view, by the space

\begin{equation}\label{1}
\mathrm{H}^1_{\rm diff}\left(\mathfrak{vect}(1);
\mathrm{D}\right)=\bigoplus_{0\leq i, j\leq n}\mathrm{H}^1_{\rm
diff}\left(\mathfrak{vect}(1);
\mathrm{D}_{\beta-j,\beta-i}\right),
\end{equation}
where $\mathrm{D}:=\mathrm{D}(n,\beta)$ is the
$\mathfrak{vect}(1)$-module of differential operators in ${\cal
S}_\beta^n$ and where $\mathrm{H}^i_\mathrm{diff}$ denotes the
differential cohomology; that is, only cochains given by
differential operators are considered. The $\mathrm{\frak{
sl}}(2)$-trivial infinitesimal deformations are classified by the
space
\begin{equation}\label{2}
{\rm H}^1_{\rm diff}\left(\mathfrak{vect}(1),\mathrm{\frak{
sl}}(2); \mathrm{D}\right)=\bigoplus_{0\leq i, j\leq n}{\mathrm
H}^1_{\rm diff}\left(\mathfrak{vect}(1),\mathrm{\frak{ sl}}(2);
\mathrm{D}_{\beta-j,\beta-i}\right).
\end{equation}

Feigin and Fuchs computed  $\mathrm{H}^1_{\rm
diff}\left(\mathfrak{vect}(1);
\mathrm{D}_{\lambda,\lambda'}\right)$, see \cite{ff}. They showed
that non-zero cohomology $\mathrm{H}^1_{\rm
diff}\left(\mathfrak{vect}(1);\mathrm{D}_{\lambda,\lambda'}\right)$
only appear for particular values of weights that we call {\it
resonant} which satisfy  $\lambda'-\lambda\in\mathbb{N}$.
Therefore, in formulas (\ref{1}) and (\ref{2}), the summations are
only over $i$ and $j$ such that $i\leq j$. Bouarroudj and Ovsienko
\cite{bo} computed ${\mathrm H}^1_{\rm
diff}\left(\mathfrak{vect}(1),\mathrm{{\frak sl}}(2);
\mathrm{D}_{\lambda,\lambda'}\right)$, and Bouarroudj \cite{b1}
solved a multi-dimensional version of the same problem on
manifolds.

In this paper we study the simplest super analog of the problem
solved in \cite{ff,bo, b1}, namely, we consider the superspace
$\mathbb{K}^{1|1}$ equipped with the contact structure determined by
a 1-form $\alpha$, and the Lie superalgebra $\mathcal{K}(1)$ of
contact polynomial vector fields on $\mathbb{K}^{1|1}$. We introduce
the $\mathcal{K}(1)$-module $\mathfrak{F}_\lambda$ of
$\lambda$-densities on $\mathbb{K}^{1|1}$ and the
$\mathcal{K}(1)$-module of linear differential operators,
$\frak{D}_{\nu,\mu}
:=\mathrm{Hom}_{\rm{diff}}(\mathfrak{F}_{\nu},\mathfrak{F}_{\mu})$,
which are super analogues of the spaces $\mathcal{F}_\lambda$ and
$\mathrm{D}_{\nu,\mu}$, respectively. The Lie superalgebra
$\mathfrak{osp}(1|2)$, a super analogue of $\mathrm{\frak {sl}}(2)$,
can be realized as a subalgebra of $\mathcal{K}(1)$. We compute
$\mathrm{H}^1_{\rm
diff}\left(\mathcal{K}(1);\mathfrak{D}_{\lambda,\lambda'}\right)$
and $\mathrm{H}^1_{\rm diff}\left(\mathcal{K}(1), \mathfrak{osp}(1|
2);\mathfrak{D}_{\lambda,\lambda'}\right)$ and we show that, as in
the classical setting, non-zero cohomology $\mathrm{H}^1_{\rm
diff}\left(\mathcal{K}(1);\mathfrak{D}_{\lambda,\lambda'}\right)$
only appear for resonant values of weights which satisfy
$\lambda'-\lambda\in{1\over2}\mathbb{N}$. So, the super analogue of
the space ${\cal S}_\beta^n$ is naturally  the superspace (see
\cite{gmo}):
\begin{equation*}
{\frak
S}^n_{\beta}=\bigoplus_{k=0}^{2n}\mathfrak{F}_{\beta-\frac{k}{2}},\quad\text{where}\quad
n\in{1\over2}\mathbb{N}.
\end{equation*}
We use the result to study formal deformations of the
$\mathcal{K}(1)$-module structure on ${\frak S}^n_{\beta}$. Denote
by $\mathfrak{D}:=\mathfrak{D}(n,\beta)$  the
$\mathcal{K}(1)$-module of linear differential operators in
${\frak S_\beta^n}$. The infinitesimal deformations of the
$\mathcal{K}(1)$-module structure on ${\frak S}^n_{\beta}$ are
classified by the space
\begin{equation*}
\mathrm{H}^1_{\rm
diff}\left(\mathcal{K}(1);\mathfrak{D}\right)=\bigoplus_{0\leq
i\leq j \leq 2n}\mathrm{H}^1_{\rm
diff}\left(\mathcal{K}(1);\frak{D}_{\beta-{j\over2},\beta-{i\over2}}\right).
\end{equation*}
The $\mathrm{\frak{ osp}}(1|2)$-trivial infinitesimal deformations
are classified by the space
\begin{equation*}
\mathrm{H}^1_{\rm diff}\left(\mathcal{K}(1), \mathfrak{osp}(1|
2);\mathfrak{D}\right)=\bigoplus_{0\leq i\leq j \leq
2n}\mathrm{H}^1_{\rm diff}\left(\mathcal{K}(1), \mathfrak{osp}(1|
2);\frak{D}_{\beta-{j\over2},\beta-{i\over2}}\right).
\end{equation*}

In this work, we study only the generic formal
$\mathfrak{osp}(1|2)$-trivial deformations of the action of
$\mathcal{K}(1)$ on the space ${\frak S}^n_{\beta}$. In order to
study the integrability of a given $\mathrm{\frak{
osp}}(1|2)$-trivial infinitesimal deformation, we need the
description of $\mathfrak{osp}(1|2)$-invariant bilinear
differential operators
$\mathfrak{F}_{\tau}\otimes\mathfrak{F}_\lambda\longrightarrow\mathfrak{F}_{\mu}$.

\section{Definitions and Notations}

\subsection{The Lie superalgebra of contact vector fields on
$\mathbb{K}^{1|n}$}

Let $\mathbb{K}^{1\mid n}$ be the superspace with coordinates
$(x,~\theta_1,\ldots,\theta_n),$ where the $\theta_i$ are odd
indeterminates equipped with the standard contact structure given by
the following $1$-form:
\begin{equation*}
\alpha_n=dx+\sum_{i=1}^n\theta_id\theta_i.
\end{equation*} On
$\mathbb{K}[x,\theta]:=\mathbb{K}[x,\theta_1,\dots,\theta_n]$, we
consider the contact bracket
\begin{equation}
\{F,G\}=FG'-F'G-\frac{1}{2}(-1)^{p(F)}\sum_{i=1}^n\overline{\eta}_i(F)\cdot
\overline{\eta}_i(G),
\end{equation}where
$\overline{\eta}_i=\frac{\partial}{\partial
{\theta_i}}-\theta_i\frac{\partial}{\partial x}$ and $p(F)$ is the
parity of $F$.

Let $\mathrm{Vect_{Pol}}(\mathbb{K}^{1|n})$ be the superspace of
polynomial vector fields on ${\mathbb{K}}^{1|n}$:
\begin{equation*}\mathrm{Vect_{Pol}}(\mathbb{K}^{1|n})=\left\{F_0\partial_x
+ \sum_{i=1}^n F_i\partial_i \mid ~F_i\in\mathbb{K}[x,\theta]~
\text{ for all } i  \right\},\end{equation*} where
$\partial_i=\frac{\partial}{\partial\theta_i}$ and
$\partial_x=\frac{\partial}{\partial x} $, and consider the
superspace $\mathcal{K}(n)$ of contact polynomial vector fields on
${\mathbb{K}}^{1|n}$. That is, $\mathcal{K}(n)$ is the superspace of
vector fields on $\mathbb{K}^{1|n}$ preserving the distribution
singled out by the $1$-form $\alpha_n$: $$
\mathcal{K}(n)=\big\{X\in\mathrm{Vect_{Pol}}(\mathbb{K}^{1|n})~|~\hbox{there
exists}~F\in {\mathbb{K}}[x,\,\theta]~ \hbox{such
that}~{L}_X(\alpha_n)=F\alpha_n\big\}. $$ The Lie superalgebra
$\mathcal{K}(n)$ is spanned by the fields of the form:
\begin{equation*}
X_F=F\partial_x
-\frac{1}{2}\sum_{i=1}^n(-1)^{p(F)}\overline{\eta}_i(F)\overline{\eta}_i,\;\text{where
$F\in \mathbb{K}[x,\theta]$.}
\end{equation*}
In particular, we have $\mathcal{K}(0)=\mathfrak{vect}(1)$.
Observe that ${L}_{X_F}(\alpha_n)=X_1(F)\alpha_n$. The bracket in
$\mathcal{K}(n)$ can be written as: $ [X_F,\,X_G]=X_{\{F,\,G\}}$.

\subsection{The subalgebra $\mathfrak{osp}(1|2)$}
In $\mathcal{K}(1)$, there is a subalgebra $\mathfrak{osp}(1|2)$
of projective transformations
\begin{equation*}\mathfrak{osp}(1|2)=\text{Span}\left(X_1,\,X_{\theta},\, X_{x},\,X_{x\theta},\,
X_{x^2}\right);\quad(\mathfrak{osp}(1|2))_{\bar{0}}=\text{Span}(X_1,\,X_{x},\,X_{x^2})
\cong\mathfrak{sl}(2) .\end{equation*}

\subsection{The space of polynomial weighted densities on
$\mathbb{K}^{1|1}$}

From now on, $n=1$ and we will denote $\alpha_1$ and
$\overline{\eta}_1$ respectively by $\alpha$ and
$\overline{\eta}$. We have analogous definition of weighted
densities in super setting (see \cite{ab}) with $dx$ replaced by
$\alpha$. The elements of these spaces are indeed (weighted)
densities since all spaces of generalized tensor fields have just
one parameter relative $\mathcal{K}(1)$
--- the value of $X_x$ on the lowest weight vector (the one
annihilated by $X_\theta$). From this point of view the volume
element (roughly speaking, $\lq\lq
dx\frac{\partial}{\partial\theta}"$) is indistinguishable from
$\alpha^{\frac{1}{2}}.$

Consider the $1$-parameter action of $\mathcal{K}(1)$ on
$\mathbb{K}[x,\theta]$ given by the rule:
\begin{equation}
\label{superaction} \mathfrak{L}^{\lambda}_{X_F}= X_F + \lambda
F',
\end{equation}
where $F'=\partial_{x}F$, or, in components:
\begin{equation}
\label{deriv}
\frak{L}^{\lambda}_{X_F}(G)=L^{\lambda}_{a\partial_x}(g_0)+\frac{1}{2}~bg_1
+\left(L^{\lambda+\frac{1}{2}}_{a\partial_x}(g_1)+\lambda
g_0b'+{1\over2} g'_0 b\right)\theta,
\end{equation}
where $F=a+b\theta,\,G=g_0+g_1\theta \in\mathbb{K}[x,\theta]$. We
denote this $\mathcal{K}(1)$-module by $\mathfrak{F}_{\lambda}$, the
space of all polynomial weighted densities on $\mathbb{K}^{1|1}$ of
weight $\lambda$:
\begin{equation}
\label{densities}
\mathfrak{F}_\lambda=\left\{f(x,\theta)\alpha^{\lambda} \mid
f(x,\theta) \in\mathbb{K}[x,\theta]\right\}.
\end{equation}
Obviously:
\begin{itemize}
\item[1)] The adjoint $\mathcal{K}(1)$-module, is isomorphic to $\mathfrak{F}_{-1}.$
\item[2)] As a $\mathfrak{vect}(1)$-module,
$\mathfrak{F}_{\lambda}\simeq\mathcal{F}_\lambda \oplus
\Pi(\mathcal{F}_{\lambda+{1\over2}})$.
\end{itemize}
Any differential operator $A$ on $\mathbb{K}^{1|1}$ can be viewed as
a linear mapping $F\alpha^\lambda\mapsto(AF)\alpha^\mu$ from
$\mathfrak{F}_{\lambda}$ to $\mathfrak{F}_\mu$, thus the space of
differential operators becomes a $\mathcal{K}(1)$-module denoted
$\mathfrak{D}_{\lambda,\mu}$ for the natural action:
\begin{equation}\label{d-action}
\mathfrak{L}^{\lambda,\mu}_{X_F}(A)=\mathfrak{L}^{\mu}_{X_F}\circ
A-(-1)^{p(A)p(F)} A\circ \mathfrak{L}^{\lambda}_{X_F}.
\end{equation}

\begin{proposition}
\label{decom} As a $\mathfrak{vect}(1)$-module, we have
\begin{equation*}(\frak{D}_{\lambda,\mu})_{\bar 0}\simeq \mathrm{D}_{\lambda,\mu}
\oplus \mathrm{D}_{\lambda+\frac{1}{2},\mu+\frac{1}{2}}\\ \hbox{and}\\
(\frak{D}_{\lambda,\mu})_{\bar1}\simeq\Pi(\mathrm{D}_{\lambda+\frac{1}{2},\mu}
\oplus \mathrm{D}_{\lambda,\mu+\frac{1}{2}}).
\end{equation*}
\end{proposition}
\begin{proof}
It is clear that the map
\begin{equation}\label{ph}\begin{array}{ll}
\varphi_\lambda:\mathfrak{F}_\lambda&\longrightarrow\mathcal{F}_{\lambda}\oplus
\Pi(\mathcal{F}_{\lambda+{1\over2}})\\
F\alpha^{\lambda}&\mapsto\left((1-\theta\partial_{\theta})(F)
(dx)^{\lambda},~
\Pi(\partial_{\theta}(F)(dx)^{\lambda+{1\over2}})\right)
\end{array}\end{equation} is
$\mathfrak{vect}(1)$-isomorphism, see formulae (\ref{deriv}). So,
we deduce a $\mathfrak{vect}(1)$-isomorphism:
\begin{equation}\label{Phi}
\begin{array}{lcll}\Phi_{\lambda,\mu}:&\frak{D}_{\lambda,\mu}&\longrightarrow&
\mathrm{D}_{\lambda,\mu}\oplus
\mathrm{D}_{\lambda+{1\over2},\mu+{1\over2}}\oplus\Pi(
\mathrm{D}_{\lambda,\mu+{1\over2}})\oplus\Pi(
\mathrm{D}_{\lambda+{1\over2},\mu})\\&A&\mapsto&\varphi_\mu\circ
A\circ\varphi_\lambda^{-1}.
\end{array}
\end{equation}Here, we identify
the $\mathfrak{vect}(1)$-modules via the following isomorphisms:
\begin{gather*}\begin{array}{llllllll}
\rm{Hom}_{\rm
diff}\left(\mathcal{F}_\lambda,\Pi(\mathcal{F}_{\mu+\frac{1}{2}})\right)
&\longrightarrow&\Pi\left(\mathrm{D}_{\lambda,\mu+\frac{1}{2}}\right),
\quad &A&\mapsto&\Pi(\Pi\circ A),\\[10pt] \rm{Hom}_{\rm
diff}\left(\Pi(\mathcal{F}_{\lambda+\frac{1}{2}}),\mathcal{F}_{\mu}\right)
&\longrightarrow&\Pi\left(\mathrm{D}_{\lambda+\frac{1}{2},\mu}\right),
\quad &A&\mapsto&\Pi(A\circ\Pi),\\[10pt] \rm{Hom}_{\rm
diff}\left(\Pi(\mathcal{F}_{\lambda+\frac{1}{2}}),\Pi(\mathcal{F}_{\mu+\frac{1}{2}})\right)
&\longrightarrow&\mathrm{D}_{\lambda+\frac{1}{2},\mu+\frac{1}{2}},
\quad &A&\mapsto&\Pi\circ A\circ\Pi.\\[10pt]
\end{array}
\end{gather*}
Note that the change of parity map $\Pi$ commutes with the
$\mathfrak{vect}(1)$-action. \end{proof}

Consider a family of $\mathfrak{vect}(1)$-modules on the space
$\mathrm{D}_{(\lambda_1,\dots,\lambda_m);\mu}$ of linear
differential operators: $~A: {\cal
F}_{\lambda_1}\otimes\cdots\otimes\mathcal{F}_{\lambda_m}\longrightarrow{\cal
F}_\mu.$ The Lie algebra $\mathfrak{vect}(1)$ naturally
acts on $\mathrm{D}_{(\lambda_1,\dots,\lambda_m);\mu}$ (by the Leibniz rule). 
We similarly  consider a family of ${\rm \mathcal{K}}(1)$-modules
on the space $\mathfrak{ D}_{(\lambda_1,\dots,\lambda_m);\mu}$ of
linear differential operators: $~A: {\frak
F}_{\lambda_1}\otimes\cdots\otimes\frak{F}_{\lambda_m}\longrightarrow{\frak
F}_\mu$. 
\section{$\mathfrak{sl}(2)$- and $\mathfrak{osp}(1|2)$-invariant bilinear differential
operators }
\begin{proposition}\label{trans2}{\rm(Gordon, \cite{pg})}  There exist
$\mathfrak{sl}(2)$-invariant bilinear differential operators,
called {\rm transvectants},
\begin{equation*}
J_k^{\tau,\lambda}:
\mathcal{F}_\tau\otimes\mathcal{F}_\lambda\longrightarrow\mathcal{F}_{\tau+\lambda+k},\quad
(\varphi dx^\tau,\phi dx^\lambda)\mapsto
J_k^{\tau,\lambda}(\varphi,\phi)dx^{\tau+\lambda+k}
\end{equation*} given by
\begin{equation*}
J_k^{\tau,\lambda}(\varphi,\phi)=\sum_{0\leq i\leq k,
i+j=k}c_{i,j}\varphi^{(i)}\phi^{(j)},
\end{equation*} where $k\in\mathbb{N}$ and the coefficients $c_{i,j}$ are characterized
as follows:
\begin{itemize}
  \item [i)] If $\tau, \lambda\not\in\{0,\,-{1\over2},\,-1,\,\dots,\,-{k-1\over2}\}$, then
$ c_{i,j}=\Gamma_{i,j,k-1}^{\tau,\lambda}$,  see  (\ref{coe}).
  \item [ii)] If  $\tau$ or $\lambda\in\{0,\,-{1\over2},\,-1,\,\dots,\,-{k-1\over2}\}$,
the coefficients $c_{i,j}$  satisfy the recurrence relation
\begin{equation}\label{cij}(i+1)(i + 2\tau)c_{i+1,j} + (j+1)(j+2\lambda)c_{i,j+1} = 0.\end{equation}
Moreover, the space of solutions of the system (\ref{cij}) is
two-dimensional if $2\lambda=-s$ and $2\tau=-t$ with $t > k-s-2$,
and one-dimensional otherwise.
\end{itemize}
\end{proposition}

Gieres and Theisen \cite{gt} listed the
$\mathfrak{osp}(1|2)$-invariant bilinear differential operators,
from $\mathfrak{F}_{\tau}\otimes\mathfrak{F}_\lambda$ to
$\mathfrak{F}_{\mu}$, called {\it supertransvectants}. Gargoubi
and Ovsienko \cite{go} gave an interpretation of these operators.
In \cite{gt}, the supertransvectants are expressed in terms of
supercovariant derivative. Here, the supertransvectants appear in
the context of the $\frak{ osp}(1|2)$-relative cohomology. More
precisely, we need to describe the $\mathfrak{osp}(1|2)$-invariant
linear differential operators from $\mathcal{K}{(1)}$ to
$\frak{D}_{\lambda,\lambda+k-1}$ vanishing on
$\mathfrak{osp}(1|2)$. Thus, using the Gordan's transvectants and
the isomorphism (\ref{ph}), we give, in the following theorem,
another description and other explicit formulas.

\begin{theorem}\label{main}  i) There are only the
following $\mathfrak{osp}(1|2)$-invariant bilinear differential
operators acting in the spaces $\frak{F}_{\lambda}$:
\begin{equation*}\begin{array}{ll}
\frak{J}_{k}^{\tau,\lambda}:
\frak{F}_\tau\otimes\frak{F}_\lambda&\longrightarrow\frak{F}_{\tau+\lambda+k}\\
(F \alpha^\tau, G \alpha^\lambda)&\mapsto
\frak{J}_k^{\tau,\lambda}(F,G)\alpha^{\tau+\lambda+k},\end{array}
\end{equation*}  where $k\in{1\over2}\mathbb{N}$.
The operators $\frak{J}_{k}^{\tau,\lambda}$ labeled by
semi-integer $k$ are odd; they are given by
\begin{equation*}\begin{array}{lllll}
\frak{J}_{k}^{\tau,\lambda}(F,G)&=\displaystyle\sum_{i+j=[k]}
\Gamma_{i,j,k}^{\tau,\lambda}\left((-1)^{p(F)}(2\tau+[k]-j)F^{(i)}\overline{\eta}(G^{(j)})-
(2\lambda+[k]-i)\overline{\eta}(F^{(i)})G^{(j)}\right).
\end{array}\end{equation*}
The operators $\frak{J}_{k}^{\tau,\lambda}$, where $k\in
\mathbb{N}$, are even; set $\frak{J}_0^{\tau,\lambda}(F,G)=FG$ and
\begin{equation*}\begin{array}{lllll}
\frak{J}_k^{\tau,\lambda}(F,G)&=\displaystyle\sum_{i+j=k-1}(-1)^{p(F)}
\Gamma_{i,j,k-1}^{\tau,\lambda}\overline{\eta}(F^{(i)})\overline{\eta}(G^{(j)})
-\displaystyle\sum_{i+j=k}\Gamma_{i,j,k-1}^{\tau,\lambda}F^{(i)}G^{(j)},\end{array}
\end{equation*}where $\big(^x_i\big)=\frac {x(x-1)\cdots (x-i+1)}{i!}$ and $[k]$ denotes
the integer part of $k$, $k>0$, and
\begin{equation}
\label{coe} \Gamma_{i,j,k}^{\tau,\lambda}=(-1)^{j}
\begin{pmatrix}2\tau+[k]\\j\end{pmatrix}
\begin{pmatrix}2\lambda+[k]\\i\end{pmatrix}.
\end{equation}

ii) If $\tau,
\lambda\not\in\{0,\,-{1\over2},\,-1,\,\dots,\,-{[k]\over2}\}$,
then $\frak{J}_{k}^{\tau,\lambda}$ is the unique (up to a scalar
factor) bilinear $\mathfrak{osp}(1|2)$-invariant bilinear
differential operator  $\frak{F}_\tau\otimes\frak{F}_\lambda
\longrightarrow\frak{F}_{\tau+\lambda+k}$.

iii) For $k\in{1\over2}(\mathbb{N}+5)$, the space of
$\mathfrak{osp}(1|2)$-invariant linear differential operators from
$\mathcal{K}{(1)}$ to $\frak{D}_{\lambda,\lambda+k-1}$ vanishing
on $\mathfrak{osp}(1|2)$ is one dimensional.
\end{theorem}

\begin{proof} i) Let $\mathcal{T}:
\frak{F}_\tau\otimes\frak{F}_\lambda\longrightarrow\frak{F}_\mu$
be an $\mathfrak{osp}(1|2)$-invariant differential operator. Using
the fact that, as $\mathfrak{vect}(1)$-modules,
\begin{equation}\label{decomposition}\frak{F}_\tau\otimes\frak{F}_\lambda\simeq
\mathcal{F}_\tau\otimes\mathcal{F}_\lambda\oplus\Pi(\mathcal{F}_{\tau+\frac{1}{2}}\otimes
\mathcal{F}_{\lambda+\frac{1}{2}})\oplus
\mathcal{F}_\tau\otimes\Pi(\mathcal{F}_{\lambda+\frac{1}{2}})\oplus\Pi(\mathcal{F}_{\tau+\frac{1}{2}})\otimes
\mathcal{F}_\lambda
\end{equation}
and $$
\frak{F}_\mu\simeq\mathcal{F}_\mu\oplus\Pi(\mathcal{F}_{\mu+\frac{1}{2}}),
$$ we can deduce that the restriction of $\mathcal{T}$ to each
component of the right hand side of (\ref{decomposition}) is a
transvectant. So, the parameters $\tau,$ $\lambda$ and $\mu$ must
satisfy $$\mu=\lambda+\tau+k,\quad\hbox{ where }\quad
k\in{1\over2}\mathbb{N}.$$ The corresponding operators will be
denoted $\frak{J}_{k}^{\tau,\lambda}$. Obviously, if $k$ is
integer, then the operator $\frak{J}_{k}^{\tau,\lambda}$ is even
and its restriction to each component of the right hand side of
(\ref{decomposition}) coincides (up to a scalar factor) with the
respective transvectants:

\begin{equation}\label{restric1}
  \begin{cases}
    & \text{${J}_k^{\tau,\lambda}~~~~~~~:
\mathcal{F}_\tau\otimes\mathcal{F}_\lambda\longrightarrow\mathcal{F}_{\mu},$}\\
     & \text{${J}_{k-1}^{\tau+\frac{1}{2},\lambda+\frac{1}{2}}:
\Pi(\mathcal{F}_{\tau+\frac{1}{2}})\otimes\Pi(\mathcal{F}_{\lambda+\frac{1}{2}})\longrightarrow\mathcal{F}_{\mu}$},
\\ & \text{${J}_{k}^{\tau,\lambda+\frac{1}{2}}~~~:
\mathcal{F}_\tau\otimes\Pi(\mathcal{F}_{\lambda+\frac{1}{2}})\longrightarrow\Pi(\mathcal{F}_{\mu+\frac{1}{2}}$}),
\\ & \text{${J}_{k}^{\tau+\frac{1}{2},\lambda}~~~:
\Pi(\mathcal{F}_{\tau+\frac{1}{2}})\otimes\mathcal{F}_\lambda\longrightarrow\Pi(\mathcal{F}_{\mu+\frac{1}{2}}$}).
\\
  \end{cases}
  \end{equation}
If $k$ is semi-integer, then the operator
$\frak{J}_{k}^{\tau,\lambda}$ is odd and its  restriction to each
component of the right hand side of (\ref{decomposition})
coincides (up to a scalar factor ) with the respective
transvectants:
\begin{equation}\label{restric2}
  \begin{cases}
    & \text{$J_{[k]+1}^{\tau,\lambda}~~~~~~~:
\mathcal{F}_\tau\otimes\mathcal{F}_\lambda\longrightarrow\Pi(\mathcal{F}_{\mu+\frac{1}{2}}),$}\\
     & \text{$J_{[k]}^{\tau+\frac{1}{2},\lambda+\frac{1}{2}}:
\Pi(\mathcal{F}_{\tau+\frac{1}{2}})\otimes\Pi(\mathcal{F}_{\lambda+\frac{1}{2}})
\longrightarrow\Pi(\mathcal{F}_{\mu+\frac{1}{2}})$}, \\ &
\text{$J_{[k]}^{\tau,\lambda+\frac{1}{2}}~~~:
\mathcal{F}_\tau\otimes\Pi(\mathcal{F}_{\lambda+\frac{1}{2}})\longrightarrow\mathcal{F}_{\mu}$},
\\ & \text{$J_{[k]}^{\tau+\frac{1}{2},\lambda}~~~:
\Pi(\mathcal{F}_{\tau+\frac{1}{2}})\otimes\mathcal{F}_\lambda\longrightarrow\mathcal{F}_{\mu}$}.
\\
  \end{cases}
  \end{equation}

More precisely, let $F\alpha^\tau\otimes G
\alpha^\lambda\in\frak{F}_\tau\otimes\frak{F}_\lambda$, where
$F=f_0+\theta f_1$ and $G=g_0+\theta g_1$, with
$f_0,\,f_1,\,g_0,\,g_1\in\mathbb{K}[x]$. Then if $k$ is integer, we
have
\begin{equation}\label{integer}
\begin{array}{llll}
\frak{J}_{k}^{\tau,\lambda}(\varphi,\psi)=&\Big[a_1J_k^{\tau,\lambda}(f_0,g_0)+
a_2J_{k-1}^{\tau+\frac{1}{2},\lambda+\frac{1}{2}}(f_1,g_1)
+\\[6pt]&\theta\left(a_3J_{k}^{\tau,\lambda+\frac{1}{2}}(f_0,g_1)+a_4J_{k}^{\tau+\frac{1}{2},\lambda}(f_1,g_0)
\right)\Big]\alpha^\mu
\end{array}
\end{equation}
and if $k$ is semi-integer, we have
\begin{equation}\label{semi}
\begin{array}{llll}
\frak{J}_{k}^{\tau,\lambda}(\varphi,\psi)=&\Big[b_1J_{[k]}^{\tau,\lambda+\frac{1}{2}}(f_0,g_1)+
b_2J_{[k]}^{\tau+\frac{1}{2},\lambda}(f_1,g_0)
+\\[6pt]&\theta\left(b_3J_{[k]+1}^{\tau,\lambda}(f_0,g_0)+b_4J_{[k]}^{\tau+\frac{1}{2},\lambda+\frac{1}{2}}(f_1,g_1)
\right)\Big]\alpha^\mu,
\end{array}
\end{equation}where the $a_i$ and $b_i$ are constants.
The invariance of $\frak{J}_{k}^{\tau,\lambda}$ with respect to
$X_{\theta}$ and $X_{x\theta}$ reads:
\begin{equation}\label{T1}
{\frak L}_{X_\theta}^\mu\circ
\frak{J}_{k}^{\tau,\lambda}-(-1)^{2k}\frak{J}_{k}^{\tau,\lambda}\circ
{\frak L}_{X_\theta}^{(\tau,\lambda)}={\frak
L}_{X_{x\theta}}^\mu\circ\frak{J}_{k}^{\tau,\lambda}
-(-1)^{2k}\frak{J}_{k}^{\tau,\lambda}\circ {\frak
L}_{X_{x\theta}}^{(\tau,\lambda)}=0.\end{equation}

The formula (\ref{T1}) allows us to determine the coefficients
$a_i$ and $b_i$. More precisely, the invariance property with
respect to $X_{\theta}$ and $X_{x\theta}$ yields
{\small\begin{equation*}
a_1=a_2=a_3=a_4,~b_2=-\frac{2\lambda+k-1}{2\tau+k-1}b_1,~ b_3
=\frac{k}{2\tau+k-1}b_1\text{ and }
b_4=-(1+\frac{2\lambda}{2\tau+k-1})b_1.
\end{equation*}}

ii) The uniqueness of supertansvectants follows from the
uniqueness of transvectants.

iii) In the non-super case, according to formulae (\ref{cij}), if
$2\tau=-1$ and $k\geq2$, the space of $\mathfrak{sl}(2)$-invariant
bilinear differential operators ${\cal F}_\tau\otimes{\cal
F}_\lambda\longrightarrow{\mathcal F}_{\tau+\lambda+k}$ is
2-dimensional if and only if $2\lambda=-s$, where
$s\in\{k-1,\,k-2\}$. This space is spanned by
$J_k^{-{1\over2},-{s\over2}}$ and $I_k^{-{1\over2},-{s\over2}}$,
where
\begin{equation*}I_k^{-{1\over2},-{s\over2}}(\varphi,\phi)=
\left\{\begin{array}{ll}\varphi\phi^{(k)}\quad&\text{if}\quad
s=k-1\\
\varphi\phi^{(k)}+k\varphi'\phi^{(k-1)}\quad&\text{if}\quad s=k-2
\end{array}\right.\end{equation*} and
\begin{equation*}
J_k^{-{1\over2},-{s\over2}}(\varphi,\phi)=\sum_{i+j=k,\;i\geq
k-s+1 }c_{i,j}\varphi^{(i)}\phi^{(j)},
\end{equation*} where the coefficients $c_{i,j}$
satisfy (\ref{cij}). We see that only the operators
$J_k^{-{1\over2},-{s\over2}}$  vanish on the space of affine
functions, i.e., of the form  $\varphi(x)=ax+b$.

If $k\geq3$, the space of $\mathfrak{sl}(2)$-invariant bilinear
differential operators ${\cal F}_{-1}\otimes{\cal
F}_\lambda\longrightarrow\mathcal{F}_{\lambda+k-1}$ is
2-dimensional if and only if $2\lambda=-s$, where
$s\in\{k-1,\,k-2,\,k-3\}$. This space is spanned by
$J_k^{-1,-{s\over2}}$ and $I_k^{-1,-{s\over2}}$, where
\begin{equation*}I_k^{-1,-{s\over2}}(\varphi,\phi)=\begin{cases}
\varphi\phi^{(k)}&\text{if $s=k-1$}\\
\varphi\phi^{(k)}+{k\over2}\varphi'\phi^{(k-1)}&\text{if $
s=k-2$}\\
\varphi\phi^{(k)}+k\varphi'\phi^{(k-1)}+{k(k-1)\over2}\varphi''\phi^{(k-2)}&\text{if
$ s=k-3$}\end{cases}\end{equation*} and where
\begin{equation*}
J_k^{-1,-{s\over2}}(\varphi,\phi)=\sum_{i+j=k,\;i\geq 3
}c_{i,j}\varphi^{(i)}\phi^{(j)}.
\end{equation*}
We see that the operator $I_k^{-1,-{s\over2}}$ does not vanish on
$\frak{sl}(2)$, but the operator $J_k^{-1,-{s\over2}}$ vanishes on
$\frak{sl}(2)$.

Now, if $\tau=-1, -{1\over2}$ and
$2\lambda\notin\{1-k,\,2-k,\,3-k\}$ with $k\geq3$, the space of
$\mathfrak{sl}(2)$-invariant bilinear differential operators
${\cal F}_\tau\otimes{\cal F}_\lambda\longrightarrow{\cal
F}_{\tau+\lambda+k}$ is 1-dimensional. But, in this case, we see
that the coefficients $c_{i,j}$ satisfying (\ref{cij}) are such
that $c_{i,j}=0$ if $i\leq2$ for $\tau=-1$ and $c_{i,j}=0$ if
$i\leq1$ for $\tau=-{1\over2}$.

Thus, in the super setting, if  $2k\geq5$, according to equations
(\ref{integer}) and (\ref{semi}), we see that the space of
$\mathfrak{osp}(1|2)$-invariant linear differential operator from
$\mathcal{K}{(1)}$ to $\frak{D}_{\lambda,\lambda+k-1}$ vanishing
on $\mathfrak{osp}(1|2)$ is one-dimensional.
\end{proof}
\section{Cohomology}
Let us first recall some fundamental concepts from cohomology
theory~(see, e.g., \cite{Fu}). Let $\frak{g}=\frak{g}_{\bar
0}\oplus \frak{g}_{\bar 1}$ be a Lie superalgebra acting on a
superspace $V=V_{\bar 0}\oplus V_{\bar 1}$ and let $\mathfrak{h}$
be a subalgebra of $\mathfrak{g}$. (If $\frak{h}$ is omitted it
assumed to be $\{0\}$.) The space of $\frak h$-relative
$n$-cochains of $\frak{g}$ with values in $V$ is the
$\frak{g}$-module
\begin{equation*}
C^n(\frak{g},\frak{h}; V ) := \mathrm{Hom}_{\frak
h}(\Lambda^n(\frak{g}/\frak{h});V).
\end{equation*}
The {\it coboundary operator} $ \delta_n: C^n(\frak{g},\frak{h}; V
)\longrightarrow C^{n+1}(\frak{g},\frak{h}; V )$ is a
$\frak{g}$-map satisfying $\delta_n\circ\delta_{n-1}=0$. The
kernel of $\delta_n$, denoted $Z^n(\mathfrak{g},\frak{h};V)$, is
the space of $\frak h$-relative $n$-{\it cocycles}, among them,
the elements in the range of $\delta_{n-1}$ are called $\frak
h$-relative $n$-{\it coboundaries}. We denote
$B^n(\mathfrak{g},\frak{h};V)$ the space of $n$-coboundaries.

By definition, the $n^{th}$ $\frak h$-relative  cohomolgy space is
the quotient space
\begin{equation*}
\mathrm{H}^n
(\mathfrak{g},\frak{h};V)=Z^n(\mathfrak{g},\frak{h};V)/B^n(\mathfrak{g},\frak{h};V).
\end{equation*}
We will only need the formula of $\delta_n$ (which will be simply
denoted $\delta$) in degrees 0 and 1: for $v \in
C^0(\frak{g},\,\frak{h}; V) =V^{\frak h}$,~ $\delta v(g) : =
(-1)^{p(g)p(v)}g\cdot v$, where
\begin{equation*}
V^{\frak h}=\{v\in V~\mid~h\cdot v=0\quad\text{ for all }
h\in\frak h\},
\end{equation*}
and  for  $ \Upsilon\in C^1(\frak{g}, \frak{h};V )$,
\begin{equation*}\delta(\Upsilon)(g,\,h):=
(-1)^{p(g)p(\Upsilon)}g\cdot
\Upsilon(h)-(-1)^{p(h)(p(g)+p(\Upsilon))}h\cdot
\Upsilon(g)-\Upsilon([g,~h])\quad\text{for any}\quad g,h\in
\frak{g}.
\end{equation*}
According to the $\mathbb{Z}_2$-grading (parity) of $\frak g$, for
any $\Upsilon\in Z^1(\frak{g}, V)$, we have
\begin{equation*}
\Upsilon=\Upsilon'+\Upsilon''\in Z^1(\frak{g}_{\bar 0};\, V)\oplus
\mathrm{Hom}(\frak{g}_{\bar 1},\, V)
\end{equation*}
subject to the following three equations:
\begin{gather}
\label{coc1} \Upsilon'([g_1,g_2]) -g_1\cdot \Upsilon'(g_2) +
g_2\cdot \Upsilon'(g_1)= 0 \quad\text{for any}\quad
g_1,\,g_2\in\frak{g}_{\bar 0}, \\[10pt] \label{coc2}
\Upsilon''([g,\,h]) - g\cdot \Upsilon''(h) +(-1)^{p(\Upsilon)}
h\cdot \Upsilon'(g)= 0 \quad\text{for any}\quad g\in\frak{g}_{\bar
0},\,h\in\frak{g}_{\bar 1},\\[10pt] \label{coc3}
\Upsilon'([h_1,h_2]) - (-1)^{p(\Upsilon)}\left(h_1\cdot
\Upsilon''(h_2)+ h_2\cdot \Upsilon''(h_1)\right)=0 \quad\text{for
any}\quad h_1,\,h_2\in\frak{g}_{\bar 1}.
\end{gather}

Formulas (\ref{coc1})--(\ref{coc3}) show that $\mathrm{H}^1_{\rm
diff}(\mathcal{K}(1);\frak{D}_{\lambda,\mu})$ and
$\mathrm{H}^1_{\rm diff}(\mathfrak{vect}(1);
\mathrm{D}_{\lambda,\mu})$ are closely related. Similarly,
$\mathrm{H}^1_{\rm
diff}(\mathcal{K}(1),\frak{osp}(1|2);\frak{D}_{\lambda,\mu})$ is
related to $\mathrm{H}^1_{\rm
diff}(\mathfrak{vect}(1),\frak{sl}(2); \mathrm{D}_{\lambda,\mu})$.
The\-refore, for comparison and to build upon, we recall the
description of $\mathrm{H}^1_{\rm diff}(\mathfrak{vect}(1);
\mathrm{D}_{\lambda,\mu})$. Note that $\mathrm{H}^1_{\rm
diff}(\mathcal{K}(1),\frak{osp}(1|2);\frak{D}_{\lambda,\mu})$ is
also computed by Conley, see \cite{c}.

\subsection{Relationship between $\mathrm{H}^1_{\rm
diff}(\mathfrak{vect}(1); \mathrm{D}_{\lambda,\mu})$ and
$\mathrm{H}^1_{\rm diff}(\mathcal{K}(1);\frak{D}_{\lambda,\mu})$}
\label{FirstSect} Feigin and Fuchs \cite{ff} calculated
$\mathrm{H}^1_{\rm diff}(\mathfrak{vect}(1);
\mathrm{D}_{\lambda,\mu})$. The result is as follows
\begin{equation}
\label{CohSpace2} \mathrm{H}^1_{\rm diff}(\mathfrak{vect}(1);
\mathrm{D}_{\lambda,\mu})\simeq\left\{
\begin{array}{ll}
\mathbb{K}&\hbox{ if }~~ \mu-\lambda=0,2,3,4 \hbox{ for all
}\lambda,\\[2pt] \mathbb{K}^2& \hbox{ if }~~\lambda=0\hbox{ ~and~
}\mu=1 ,\\[2pt] \mathbb{K}&\hbox{ if }~~ \lambda=0  \hbox{ or }
\lambda=-4\hbox{ ~and~ }\mu-\lambda=5, \\[2pt] \mathbb{K}& \hbox{
if }~~ \lambda=-\frac{5\pm \sqrt{19}}{2}\hbox{ ~and~ }\mu-\lambda
=6,\\[2pt] 0 &\hbox{ otherwise. }
\end{array}
\right.
\end{equation}
For $X\frac{d}{dx}\in\mathfrak{vect}(1)$ and
$f{dx}^{\lambda}\in{\cal F}_\lambda$, we write
\begin{align*}\begin{array}{llll}
C_{\lambda,\lambda+k
}(X\frac{d}{dx})(f{dx}^{\lambda})=C_{\lambda,\lambda+k
}(X,f){dx}^{\lambda+k}. \end{array}\end{align*} The spaces
${\mathrm H}^1_{\rm dif\/f}(\mathfrak{vect}(1),
\mathrm{D}_{\lambda,\lambda+k})$ are generated by the cohomology
classes of the following 1-cocycles:
\begin{equation}\label{cocycles}\begin{array}{llllllllll}
  C_{\lambda,\lambda}(X,f)&=&X'f \\
C_{0,1}(X,f)&=&X''f \\
   {\widetilde C}_{0,1}(X,f)&=&(X'f)' \\
   C_{\lambda,\lambda+2}(X,f)&=&X^{(3)}f+2X''f' \\
  C_{\lambda,\lambda+3}(X,f)&=&X^{(3)}f'+X''f'' \\
  C_{\lambda,\lambda+4}(X,f)&=&-\lambda
X^{(5)}f+X^{(4)}f'-6X^{(3)}f''-4X''f^{(3)}\\
C_{0,5}(X,f)&=&2X^{(5)}f'-5X^{(4)}f''+10X^{(3)}f^{(3)}+5X''f^{(4)}\\
C_{-4,1}(X,f)&=&12X^{(6)}f+22X^{(5)}f'+5X^{(4)}f''-10X^{(3)}f^{(3)}-5X''f^{(4)}\\
  C_{a_i,a_i+6}(X,f)&=&\alpha_i X^{(7)}f-\beta_i X^{(6)}f'-\gamma_i
X^{(5)}f''-
  5X^{(4)}f^{(3)}+5X^{(3)}f^{(4)}+~&2X''f^{(5)},
\end{array}
\end{equation}
where\begin{equation*}
\begin{array}{llllllll}a_1=-\frac{5+ \sqrt{19}}{2},
&\alpha_1=-\frac{22+ 5\sqrt{19}}{4}, &\beta_1=\frac{31+
7\sqrt{19}}{2}, &\gamma_1=\frac{25+ 7\sqrt{19}}{2}\\[4pt]
a_2=-\frac{5- \sqrt{19}}{2}, &\alpha_2=-\frac{22- 5\sqrt{19}}{4},
& \beta_2=\frac{31- 7\sqrt{19}}{2},&\gamma_2=\frac{25-
7\sqrt{19}}{2}.\end{array}\end{equation*}

Now, let us study the relationship between any 1-cocycle of
${\mathcal K}(1)$ and its restriction to the subalgebra
$\mathfrak{vect}(1)$. More precisely, we study the relationship
between $\mathrm{H}_{\rm diff}^1({\mathcal
K}(1);\mathfrak{D}_{\lambda,\mu})$ and $\mathrm{H}_{\rm
diff}^1(\mathfrak{vect}(1); \mathrm{D}_{\lambda,\mu})$. According
to Proposition \ref{decom}, we see that $\mathrm{H}_{\rm
diff}^1(\mathfrak{vect}(1);\mathfrak{D}_{\lambda,\mu})$ can be
deduced from the spaces $\mathrm{H}_{\rm
diff}^1(\mathfrak{vect}(1); \mathrm{D}_{\lambda,\mu})$:
\begin{equation}\label{coho}\begin{array}{ll}\mathrm{H}_{\rm
diff}^1\left(\mathfrak{vect}(1);\mathfrak{D}_{\lambda,\mu}\right)&\simeq\mathrm{H}_{\rm
diff}^1\left(\mathfrak{vect}(1);
\mathrm{D}_{\lambda,\mu}\right)\oplus\mathrm{H}_{\rm
diff}^1\left(\mathfrak{vect}(1);
\mathrm{D}_{\lambda+\frac{1}{2},\mu+\frac{1}{2}}\right)\oplus\\[8pt]
&\mathrm{H}_{\rm diff}^1\left(\mathfrak{vect}(1);\Pi(
\mathrm{D}_{\lambda,\mu+\frac{1}{2}})\right) \oplus\mathrm{H}_{\rm
diff}^1\left(\mathfrak{vect}(1);\Pi(
\mathrm{D}_{\lambda+\frac{1}{2},\mu})\right).
\end{array}\end{equation}
Moreover, the following lemma shows the close relationship between
the cohomolgy spaces
$\mathrm{H}^1(\mathcal{K}(1);\mathfrak{D}_{\lambda,\mu})$ and
$\mathrm{H}^1(\mathfrak{vect}(1);\mathrm{D}_{\lambda,\mu})$.
\begin{lemma}\label{sa}
The 1-cocycle $\Upsilon$ of $\mathcal{K}(1)$ is a coboundary if
and only if its restriction $\Upsilon'$ to $\mathfrak{vect}(1)$ is
a coboundary.
\end{lemma}
\begin{proof} It is easy to see that if $\Upsilon$ is a coboundary
of $\mathcal{K}(1)$, then  $\Upsilon'$ is a coboundary of
$\mathfrak{vect}(1)$. Now, assume that $\Upsilon'$ is a coboundary
of $\mathfrak{vect}(1)$, that is, there exist
${A}\in\frak{D}_{\lambda,\mu}$ such that $\Upsilon'$ is defined by
\begin{equation*}
\Upsilon'(X_f)=\mathfrak{L}_{X_f}^{\lambda,\mu}{A}\quad\text{ for
all } f\in\mathbb{K}[x].
\end{equation*}
By replacing $\Upsilon$ by $\Upsilon-\delta{A}$, we can suppose
that $\Upsilon'=0$. But, in this case, the map $\Upsilon$ must
satisfy the following equations
\begin{align}
\label{sltr1}
&\mathfrak{L}^{\lambda,\mu}_{X_g}\Upsilon(X_{h\theta})-\Upsilon([X_g,X_{h\theta}])=0
\quad\text{ for all } g,\,h\in\mathbb{K}[x].\\[6pt] \label{sltr2}
&\mathfrak{L}^{\lambda,\mu}_{X_{h_1\theta}}\Upsilon(X_{h_2\theta})+
\mathfrak{L}^{\lambda,\mu}_{X_{h_2\theta}}\Upsilon(X_{h_1\theta})=0
\quad\text{ for all } h_1,\,h_2\in\mathbb{K}[x].
\end{align}
The equation (\ref{sltr1}) expresses the
$\mathfrak{vect}(1)$-invariance of the map
$\Upsilon:\Pi(\mathcal{F}_{-{1\over2}})\times
\mathcal{F}_\lambda\longrightarrow\mathcal{F}_\mu$. Therefore, if
$\Upsilon$ is an even 1-cocycle, then, according to Proposition
\ref{decom}, we can easily deduce the expression of $\Upsilon$
from the work of P.~Grozman \cite{G4}. More precisely, $\Upsilon$
has, {\it a priori}, the following form:
\begin{equation*}
\Upsilon(X_{h\theta })(F\alpha^\lambda)=\left\{\begin{array}{llll}
(a_1hf\theta)\alpha^{\lambda-1}&\text{if
}\quad\mu=\lambda-1\\[6pt] (a_2hg+a_3({1\over2}hf'+ \lambda
h'f)\theta)\alpha^\lambda &\text{if }\quad\mu=\lambda\\[6pt]
a_4({1\over2}hg'+ \lambda h'g)\alpha^{\lambda+1} &\text{if
}\quad\mu=\lambda+1, \,\lambda\neq0,-{1\over2} \\[6pt]
a_5({1\over2}hg''+ h'g')\alpha^{{3\over2}}&\text{if
}\quad(\lambda,\mu)=(-{1\over2},{3\over2})\\[6pt]
\left(a_6(hg'+h'g)+a_7({1\over2}hf''+ h'f')\theta\right)\alpha
&\text{if }\quad(\lambda,\mu)=(0,1)\\[6pt] a_8(hg''-
h''g)\alpha\quad&\text{if }(\lambda,\mu)=(-1,1)\\[6pt] (a_9hg'+
a_{10}(hf''- h''f)\theta)\alpha^{1\over2}&\text{if
}\quad(\lambda,\mu)=(-{1\over2},{1\over2})\\[6pt] 0
&\text{otherwise},\end{array}\right.
\end{equation*}
where $a_i\in\mathbb{K}$, $f,\,g\in\mathbb{K}[x]$ and $F=f+g\theta$.
But, the map $\Upsilon$ must satisfy the equation (\ref{sltr2}), so
we obtain $a_1=a_4=a_5=a_8=0,~a_3=-2a_2,~a_7=-2a_6$ and
$a_{10}=-a_9$. More precisely, up to a scalar factor, $\Upsilon$ is
given by:
\begin{equation*}\Upsilon=\left\{\begin{array}{lllll}
\delta  ((1-\theta\partial_\theta)\partial_x)&\text{ if
}\quad(\lambda,\mu)=(0,1),\\
\delta(\theta\partial_\theta\partial_x)&\text{ if
}\quad(\lambda,\mu)=(-{1\over2},{1\over2}),\\ \delta
(\theta\partial_\theta)&\text{ if }\quad\lambda=\mu,\\ 0 &\text{
otherwise. }
\end{array}\right.
\end{equation*}
Similarly, if $\Upsilon$ is an odd 1-cocycle, then, $\Upsilon$
has, a priori, the following form (see \cite{G4}):
\begin{equation*}
\Upsilon(X_{h\theta })(F\alpha^\lambda)=\left\{\begin{array}{llll}
(b_1hf+b_2hg)\alpha^{\lambda-{1\over2}}&\text{if
}\quad\mu=\lambda-{1\over2}\\[6pt] (b_3({1\over2} hf'+\lambda
h'f)+b_4({1\over2}
hg'+(\lambda+{1\over2})h'g)\theta)\alpha^{\lambda+{1\over2}}&\text{if
}\quad\mu=\lambda+{1\over2}\\[6pt] b_5({1\over2}
hf''+h''f)\alpha^{\frac{3}{2}}&\text{if
}\quad(\lambda,\mu)=(0,{3\over2})\\[6pt] (b_6(hf''-h''f)+b_7({1\over2}
hg''+h'g')\theta)\alpha &\text{if
}\quad(\lambda,\mu)=(-{1\over2},1)\\[6pt]
(b_8(hg''-h''g)\theta\alpha^{{1\over2}}&\text{if
}\quad(\lambda,\mu)=(-1,{1\over2})
\\[6pt] 0 &\text{otherwise,}\end{array}\right.
\end{equation*}
where $b_i\in\mathbb{K}$. But, the map $\Upsilon$ must satisfy the
equation (\ref{sltr2}), so we obtain $b_5=b_8=0,~b_1=b_2,~b_3=b_4$
and $b_{7}=2b_6$. More precisely, up to a scalar factor, $\Upsilon$
is given by:
\begin{equation*}\Upsilon=\left\{\begin{array}{lllll}
\delta  (\partial_\theta)&\text{ if}\quad\mu=\lambda+{1\over2},\\
\delta  (\theta)&\text{ if }\quad\mu=\lambda-{1\over2},\\ \delta
(\partial_\theta\partial_x)&\text{ if
}\quad(\lambda,\mu)=(-{1\over2},1),\\ 0 &\text{ otherwise. }
\end{array}\right.
\end{equation*}
This completes the proof. \end{proof}

The following lemma gives the general form of any 1-cocycle of
$\mathcal{K}(1)$.
\begin{lemma}
\label{sd} Let $\Upsilon\in\mathrm{Z}^1_{\rm
diff}(\mathcal{K}(1);\frak{D}_{\lambda,\mu})$. Up to a coboundary,
the map $\Upsilon$ has the following general form
\begin{equation}\label{coef}
\Upsilon(X_{F})=\sum_{m,k}(a_{m,k}+\theta
b_{m,k})\overline{\eta}^{k}(F)\overline{\eta}^m,
\end{equation}where the coefficients  $a_{m,k}$ and $b_{m,k}$
are constants.
\end{lemma}

\begin{proof} Since $-\overline{\eta}^2={\partial _x}$,
the operator $\Upsilon$ has the form (\ref{coef}), where, {\it a
priori}, the coefficients $a_{m,k}$ and ${b}_{m,k}$ are functions
(see \cite{gmo}), but we will prove that, up to a coboundary,
$\Upsilon$ is invariant with respect the vector field
$X_1={\partial_x}$. The 1-cocycle condition reads:
\begin{equation}
\begin{array}{lll}\label{partial1}
\mathfrak{L}^{\lambda,\mu}_{X_1}(\Upsilon(X_{F}))-
(-1)^{p(F)p(\Upsilon)}\mathfrak{L}^{\lambda,\mu}_{X_{F}}(\Upsilon(X_1))-
\Upsilon([X_1,X_{F}])=0.
\end{array}
\end{equation}
But, from  (\ref{cocycles}), up to a coboundary, we have
$\Upsilon(X_1)=0$, and therefore the equation (\ref{partial1})
becomes
\begin{equation*}
\begin{array}{lll}\label{}
\mathfrak{L}^{\lambda,\mu}_{X_1}(\Upsilon(X_{F}))-
\Upsilon([X_1,X_{F}])=0
\end{array}
\end{equation*}
which is nothing but the invariance property of $\Upsilon$ with
respect the vector field $X_1$.
\end{proof}

\begin{lemma}\label{osp} Any 1-cocycle $\Upsilon\in
Z^1_{\mathrm{diff}}(\mathcal{K}(1);\frak{D}_{\lambda,\mu})$
vanishing on $\frak{osp}(1|2)$ is $\mathfrak{osp}(1|2)$-invariant.
\end{lemma}
\begin{proof} The 1-cocycle relation of $\Upsilon$ reads:
\begin{equation}\label{osp1}
(-1)^{p(F)p(\Upsilon)}\mathfrak{L}_{X_F}^{\lambda,\mu}
\Upsilon(X_G)-(-1)^{p(G)(p(F)+p(\Upsilon))}\mathfrak{L}_{X_G}^{\lambda,\mu}
\Upsilon(X_F)-\Upsilon([X_F,~X_G])=0,
\end{equation}
where $X_F,\,X_G\in ~\mathcal{K}(1).$  Thus, if $\Upsilon(X_F)=0$
for all $X_F\in\frak{osp}(1|2)$, the equation (\ref{osp1}) becomes
\begin{equation}\label{osp2}
(-1)^{p(F)p(\Upsilon)}\mathfrak{L}_{X_F}^{\lambda,\mu}
\Upsilon(X_G)-\Upsilon([X_F,~X_G])=0
\end{equation}
 expressing the $\frak{osp}(1|2)$-invariance  of
$\Upsilon$. \end{proof}

\begin{lemma}\label{sl2} (\cite{bab} Lemma 3.3.) Up to a coboundary, any 1-cocycle $\Upsilon\in
Z_{\mathrm{diff}}^1(\mathcal{K}(1);\frak{D}_{\lambda,\mu})$
vanishing on $\frak{sl}(2)$ is $\mathfrak{osp}(1|2)$-invariant.
That is, if $\Upsilon(X_1)=\Upsilon(X_x)=\Upsilon(X_{x^2})=0$,
then the restriction of $\Upsilon$ to $\frak{osp}(1|2)$ is
trivial.
\end{lemma} \begin{proof} Recall that, as
$\mathfrak{sl}(2)$-module, the subalgebra $\mathfrak{osp}(1|2)$ is
isomorphic to $\mathfrak{sl}(2)\oplus\mathfrak{a}$, where
$\mathfrak{a}=\text{Span}(X_\theta,\,X_{x\theta})$. Consider a
linear operator $A:
\mathfrak{a}\rightarrow\mathrm{D}_{\lambda,\mu} $. By a
straightforward computation, we show that if $A$ is
$\mathfrak{sl}(2)$-invariant, then $\mu=\lambda-\frac{1}{2}+k,$
where $k\in\mathbb{N} $ and the  corresponding operator $A_k$ has
the following expression
\begin{equation}\label{inv}
A_k(X_{h\theta})(fdx^{\lambda})=a_k\left(h f^{(k)}
+k(2\lambda+k-1)h'f^{(k-1)}\right)dx^{\lambda-\frac{1}{2}+k},
\end{equation}
where
\begin{equation*}k(k-1)(2\lambda+k-1)(2\lambda+k-2)a_k=0.\end{equation*}

Now, consider $\Upsilon\in
Z_{\mathrm{diff}}^1(\mathcal{K}(1);\frak{D}_{\lambda,\mu})$ such
that $\Upsilon(X_1)=\Upsilon(X_x)=\Upsilon(X_{x^2})=0$. The
1-cocycle relations give, for all $h$, $h_1$, $h_2$ polynomial
with degree 0 or 1 and $g$ polynomial with degree 0, 1 or 2, the
following equations
\begin{align}
 \label{sltriv1} &\mathfrak{L}^{\lambda,\mu}_{X_g}\Upsilon(X_{h\theta})-
 \Upsilon([X_g,X_{h\theta}])=0, \\
 \label{sltriv2} &\mathfrak{L}^{\lambda,\mu}_{X_{h_1\theta}}\Upsilon(X_{h_2\theta}) +
 \mathfrak{L}^{\lambda,\mu}_{X_{h_2\theta}}      \Upsilon(X_{h_1\theta})=0.
\end{align}
1) If $\Upsilon$ is an even 1-cocycle, then, according to
Propostion \ref{decom}, its restriction to $\mathfrak{a}$ is
decomposed into two maps:
$\mathfrak{a}\rightarrow\Pi(\mathrm{D}_{\lambda+\frac{1}{2},\mu})$
and
$\mathfrak{a}\rightarrow\Pi(\mathrm{D}_{\lambda,\mu+\frac{1}{2}})$.
The equation (\ref{sltriv1}) tell us that these maps are
$\mathfrak{sl}(2)$-invariant. Therefore, their expressions are
given by (\ref{inv}). So, we must have
$\mu=\lambda+k=(\lambda+\frac{1}{2})-\frac{1}{2}+k$ (and then
$\mu+\frac{1}{2}=\lambda-\frac{1}{2}+k+1$). More precisely, using
the equation (\ref{sltriv2}), we get (up to a scalar  factor):
\begin{equation*}
\Upsilon_{|\mathfrak{osp}(1|2)}=\left\{\begin{aligned} &0~~\text{
if }~k(k-1)(2\lambda+k)(2\lambda+k-1)\neq0~~\text{ or }~k=1 \text{
and }\lambda\notin\{0,\,-\frac{1}{2}\},\\
&\delta(\theta\partial_\theta\partial_x^k)~~\text{ if
}~(\lambda,\mu)= \left(\frac{-k}{2},\frac{k}{2}\right),\\ &\delta
(\partial_x^k-\theta\partial_\theta\partial_x^k)~~\text{ if
}~(\lambda,\mu)= \left(\frac{1-k}{2},\frac{1+k}{2}\right)~~\text{
or }~\lambda=\mu.
\end{aligned}\right.
\end{equation*}
2) Similarly, if $\Upsilon$ is an odd 1-cocycle, we get:
\begin{equation*}
\Upsilon_{|\mathfrak{osp}(1|2)}=\left\{\begin{array}{lll} 0&\text{
if }& k(k-1)(2\lambda+k-1)\neq0,\\ \delta(\theta)&\text{ if
}&\mu=\lambda-\frac{1}{2},\\ \delta(\partial_\theta)&\text{ if
}&\mu=\lambda+\frac{1}{2},\\ \delta(\theta\partial_x^k)&\text{ if
}&(\lambda,\mu)=(\frac{1-k}{2},\frac{k}{2}).
\end{array}\right.
\end{equation*}
\end{proof}

Now, we can compute $\mathrm{H}^1_{\rm
diff}(\mathcal{K}(1);\mathfrak{D}_{\lambda,\mu})$ and the
$\frak{osp}(1|2)$-relative cohomology\\
$\mathrm{H}^1_{\mathrm{diff}}(\mathcal{K}(1),\mathfrak{osp}(1|2);\frak{D}_{\lambda,\mu})$.
Let $\Upsilon$ be any 1-cocycle over $\mathcal{K}(1)$.
According to Proposition \ref{decom}, we have
\begin{equation*}
\Upsilon_{|\mathfrak{vect}(1)}\in\mathrm{H}_{\rm
diff}^1\left(\mathfrak{vect}(1);
\mathrm{D}_{\lambda,\mu}\right)\oplus\mathrm{H}_{\rm
diff}^1\left(\mathfrak{vect}(1);
\mathrm{D}_{\lambda+\frac{1}{2},\mu+\frac{1}{2}}\right)\quad\text{if}\quad\Upsilon\quad\text{is
even}
\end{equation*}
and
\begin{equation*}
\Upsilon_{|\mathfrak{vect}(1)}\in\mathrm{H}_{\rm
diff}^1\left(\mathfrak{vect}(1);\Pi(
\mathrm{D}_{\lambda,\mu+\frac{1}{2}})\right) \oplus\mathrm{H}_{\rm
diff}^1\left(\mathfrak{vect}(1);\Pi(
\mathrm{D}_{\lambda+\frac{1}{2},\mu})\right)\quad\text{if}\quad\Upsilon\quad\text{is
odd}.
\end{equation*}
We know that non-zero cohomology $\mathrm{H}^1_{\rm
diff}\left(\mathfrak{vect}(1);
\mathrm{D}_{\lambda,\lambda'}\right)$ only appear if
$\lambda'-\lambda\in\mathbb{N}$. Thus, according to Lemma
\ref{sa}, the following statements hold:
\begin{itemize}
  \item [i)] If
  $\mu-\lambda\notin{1\over2}(\mathbb{N}-1)$, then
  $\mathrm{H}^1_{\mathrm{diff}}(\mathcal{K}(1);\frak{D}_{\lambda,\mu})=0$.
  \item [ii)] If
  $\mu-\lambda$ is integer, then
  $\mathrm{H}^1_{\mathrm{diff}}(\mathcal{K}(1);\frak{D}_{\lambda,\mu})$ is
  spanned only by the cohomology classes of even cocycles.
  \item [iii)]  If
  $\mu-\lambda$ is semi-integer, then
  $\mathrm{H}^1_{\mathrm{diff}}(\mathcal{K}(1);\frak{D}_{\lambda,\mu})$ is
  spanned only by the cohomology classes of odd cocycles.
\end{itemize}

\subsection{The space
$\mathrm{H}^1_{\mathrm{diff}}(\mathcal{K}(1),\mathfrak{osp}(1|2);\frak{D}_{\lambda,\mu})$}
The main result of this subsection is the following:
\begin{theorem}
\label{th3}
$\rm{dim}\mathrm{H}^1_{\mathrm{diff}}(\mathcal{K}(1),\frak{
osp}(1|2);\frak{D}_{\lambda,\mu})=1$ if
\begin{equation*}
\begin{array}{llllll}
\mu-\lambda=\frac{3}{2} &\hbox{ and } \lambda\neq-{1\over2},\\
\mu-\lambda=2 &\hbox{ for all } \lambda,\\ \mu-\lambda=\frac{5}{2}
&\hbox{ and } \lambda\neq-1,\\ \mu-\lambda=3 &\hbox{ and }
\lambda\in\{0,\,-\frac{5}{2}\},\\ \mu-\lambda=4 &\hbox{ and }
\lambda=\frac{-7\pm\sqrt{33}}{4}.
\end{array}
\end{equation*}
Otherwise, $\mathrm{H}^1_{\mathrm{diff}}(\mathcal{K}(1),\frak{
osp}(1|2);\frak{D}_{\lambda,\mu})=0$.

The corresponding  spaces ${\mathrm
H}^1_{\mathrm{diff}}(\mathcal{K}(1),\frak{
osp}(1|2);\frak{D}_{\lambda,\lambda+\frac{k}{2}})$ are spanned by
the cohomology classes of
$\Upsilon_{\lambda,\lambda+\frac{k}{2}}=\frak{J}_{\frac{k}{2}+1}^{-1,\lambda}$,
where $k\in\{3,\,4,\,5,\,6,\,8\}$.
\end{theorem}

\begin{proof} Note that, by Lemma \ref{sa}, the $\frak{
osp}(1|2)$-relative cocycles are related to its homologous in the
classical setting, and by Lemma \ref{osp}, they are
supertransvectants. Bouarroudj and Ovsienko \cite{bo} showed that
\begin{equation}\label{h1}
\mathrm{H}^1_{\rm diff}(\mathfrak{vect}(1),{\rm\frak
sl}(2);\mathrm{D}_{\lambda,\lambda+k})\simeq \left\{
\begin{array}{lll}
\mathbb{K} & \hbox{if} ~~\left\{
\begin{array}{l}
k=2 \hbox{ and } \lambda\neq-{1\over2},\\ k=3 \hbox{ and }
\lambda\neq-1,\\ k=4 \hbox{ and } \lambda\neq-\frac{3}{2},\\ k=5
\hbox{ and } \lambda=0,-4,\\ k=6 \hbox{ and } \lambda=-\frac{5\pm
\sqrt{19}}{2},
\end{array}
\right.
\\[16pt]
0 &\hbox{otherwise}.
\end{array}
\right.
\end{equation}
 These spaces  are generated by the cohomology classes of
the following non-trivial $\frak{sl}(2)$-relative 1-cocycles,
$A_{\lambda,\lambda+k}:$
\begin{equation*}\begin{array}{lllllll}
A_{\lambda,\lambda+2}(X,f)=X^{(3)}f,\qquad
\lambda\neq-\frac{1}{2},\\[2pt]
A_{\lambda,\lambda+3}(X,f)=X^{(3)}f'-\frac{\lambda}{2}X^{(4)}f,\qquad \lambda\neq-1, \\[2pt]
A_{\lambda,\lambda+4}(X,f)=X^{(3)}f''-\frac{2\lambda+1}{2}X^{(4)}f'+\frac{\lambda(2\lambda+1)}{10}X^{(5)}f,
\qquad \lambda\neq-\frac{3}{2},\\[2pt]
A_{0,5}(X,f) =-3X^{(5)}f'+ 15X^{(4)}f'' -10X^{(3)}f^{(3)},\\[2pt]
 A_{-4,1}(X,f)=28 X^{(6)}f+63X^{(5)}f'+ 45X^{(4)}f''
+10X^{(3)}f^{(3)}\\[2pt]
A_{a_i,a_i+6}(X,f)=\alpha_iX^{(7)}f-14\beta_iX^{(6)}f'-126\gamma_iX^{(5)}f''
-210\tau_iX^{(4)}f^{(3)}+210X^{(3)}f^{4}\\
\end{array}\end{equation*}
where $\tau_1=-2+\sqrt{19}\quad\text{and}\quad
\tau_2=-2-\sqrt{19}$. The $a_i$, $\alpha_i$, $\beta_i$ and
$\gamma_i$ are those given in (\ref{cocycles}).

\medskip

So, we see first that if $2(\mu-\lambda)\notin\{3,\,\dots,\,13\}$,
then  by Lemma \ref{sa}, the corresponding cohomology ${\mathrm
H}^1_{\mathrm{diff}}(\mathcal{K}(1),\frak{
osp}(1|2);\frak{D}_{\lambda,\mu})$
vanish. Indeed, let $\Upsilon$ be any element of\\
$Z_{\mathrm{diff}}^1(\mathcal{K}(1),\frak{
osp}(1|2);\mathfrak{D}_{\lambda,\mu})$. Then by (\ref{h1}) and
(\ref{coho}), up to a coboundary, the restriction of $\Upsilon$ to
$\mathfrak{vect}(1)$ vanishes, so $\Upsilon=0$ by Lemma \ref{sa}.
By the same arguments, if $2(\mu-\lambda)>9$, generically, the
corresponding cohomology vanish.

For $2(\mu-\lambda)\in\{3,\,\dots,\,13\}$, we study the
supertranvectant $\frak{J}_{\mu-\lambda+1}^{-1,\lambda}$. If it is a
non-trivial 1-cocycle, then the corresponding cohomology space is
one-dimensional, otherwise it is zero. To study any supertranvectant
$\frak{J}_{\mu-\lambda+1}^{-1,\lambda}$ satisfying
$\delta(\frak{J}_{\mu-\lambda+1}^{-1,\lambda})=0$, we consider the
two components of its restriction to $\mathfrak{vect}(1)$ which we
compare with  $A_{\lambda,\mu}$ and
$A_{\lambda+{1\over2},\mu+{1\over2}}$ or $A_{\lambda+{1\over2},\mu}$
and $A_{\lambda,\mu+{1\over2}}$ depending on whether $\lambda-\mu$
is integer or semi-integer. For instance, we  show that
$\frak{J}_{5\over2}^{-1,\lambda}$ is a 1-cocycle. Moreover, it is
non-trivial for $\lambda\neq-{1\over2}$ since, for
$g,\,f\in\mathbb{K}[x]$, we have
$\frak{J}_{5\over2}^{-1,\lambda}(X_g)(f)=-\theta
A_{\lambda,\lambda+2}(g,f)$. More precisely, we get the following
non-trivial 1-cocycles:
\begin{equation*}\left\{\begin{array}{lllllllll}
\Upsilon_{\lambda,\lambda+\frac{3}{2}}(X_G)(F\alpha^{\lambda})&=
&\overline{\eta}(G'')F\alpha^{\lambda+\frac{3}{2}}\quad\text{
for } \lambda\neq-{1\over2},\\[4pt]
\Upsilon_{\lambda,\lambda+\frac{5}{2}}(X_G)(F\alpha^{\lambda})&=&\left(
2\lambda\overline{\eta}(G^{(3)})F-{3}\overline{\eta}(G'')F'
-(-1)^{p(G)}G^{(3)}\overline{\eta}(F)
\right)\alpha^{\lambda+\frac{5}{2}}~\text{ for }\lambda\neq-1,\\[4pt]
\Upsilon_{\lambda,\lambda+2}(X_G)(F\alpha^{\lambda})&=&\left({2\over3}\lambda
G^{(3)}F-(-1)^{p(G)}\overline{\eta}(G'') \overline{\eta}(F)
\right)\alpha^{\lambda+2}\quad\text{ for all }\lambda,\\[4pt]
\Upsilon_{\lambda,\lambda+3}(X_G)(F\alpha^{\lambda})&=&\Big((-1)^{p(G)}\overline{\eta}(G'')\overline{\eta}(F')
-{2\lambda+1\over3}\big((-1)^{p(G)}\overline{\eta}(G^{(3)})\overline{\eta}(F)
+G^{(3)}F'\big)\\[2pt]&~&+{\lambda(2\lambda+1)\over6}G^{(4)}F
\Big)\alpha^{\lambda+3}\quad \text{ for
}\lambda=0,\,-\frac{5}{2},\\[4pt]
\Upsilon_{\lambda,\lambda+4}(X_G)(F\alpha^{\lambda})&=&\Big((-1)^{p(G)}\overline{\eta}(G'')\overline{\eta}(F'')
-{2(\lambda+1)\over3}\big(2(-1)^{p(G)}\overline{\eta}(G^{(3)})\overline{\eta}(F')
+G^{(3)}F''\big)\\[2pt]&~&+{(\lambda+1)(2\lambda+1)\over6}\big((-1)^{p(G)}\overline{\eta}
(G^{(4)})\overline{\eta}(F)+2G^{(4)}F'\big)\\[2pt]&~&-
{\lambda(\lambda+1)(2\lambda+1)\over15}G^{(5)}F\Big)\alpha^{\lambda+4}\quad\text{
for } \lambda=\frac{-7\pm\sqrt{33}}{4}.\end{array}\right.
\end{equation*}
\end{proof}

\subsection{The  space $\mathrm{H}^1_{\rm
diff}(\mathcal{K}(1);\mathfrak{D}_{\lambda,\mu}$)}
\begin{theorem}
\label{th1} ${\rm dim}\mathrm{H}^1_{\rm
diff}(\mathcal{K}(1);\mathfrak{D}_{\lambda,\mu})=1$ if
\begin{equation*}
\begin{array}{llll}
\mu-\lambda=0 &\hbox{ for all } \lambda,\\
\mu-\lambda=\frac{3}{2}&\hbox{ for all }~\lambda,\\ \mu-\lambda=2
&\hbox{ for all } \lambda,\\ \mu-\lambda=\frac{5}{2}&\hbox{ for
all }~\lambda,\\ \mu-\lambda=3 &\hbox{ and }
\lambda\in\{0,\,-\frac{5}{2}\},\\ \mu-\lambda=4 &\hbox{ and }
\lambda=\frac{-7\pm\sqrt{33}}{4}.
\end{array}
\end{equation*}
${\rm dim}\mathrm{H}^1_{\rm
diff}(\mathcal{K}(1);\frak{D}_{0,\frac{1}{2}})=2$. Otherwise,
$\mathrm{H}^1_{\rm
diff}(\mathcal{K}(1);\frak{D}_{\lambda,\mu})=0$.

The spaces $\mathrm{H}^1_{\rm
diff}(\mathcal{K}(1);\frak{D}_{\lambda,\mu})$ are spanned by the
cohomology classes of the $1$-cocycles $\Upsilon_{\lambda,\mu}$
given in Theorem \ref{th3} and by the cohomology classes of the
following $1$-cocycles:
\begin{equation*}\left\{\begin{array}{llllll}\Upsilon_{\lambda,\lambda}(X_G)(F\alpha^{\lambda})&=&G'F\alpha^{\lambda},\\[4pt]
\Upsilon_{0,\frac{1}{2}}(X_G)(F)&=&
\overline{\eta}(G')F\alpha^{{1\over2}},\\[4pt]
\widetilde{\Upsilon}_{0,\frac{1}{2}}(X_G)(F)&=&
\overline{\eta}(G'F)\alpha^{{1\over2}},\\[4pt]
\Upsilon_{-\frac{1}{2},1}(X_G)(F\alpha^{-{1\over2}})&=&\Big(\overline{\eta}(G'')F
+\overline{\eta}(G')F'+(-1)^{p(G)}G''\overline{\eta}(F)\Big)\alpha\\[4pt]
\Upsilon_{-1,\frac{3}{2}}(X_G)(F\alpha^{-1})&=&
\Big((-1)^{p(G)}(G'''\overline{\eta}(F)+2G''\overline{\eta}(F'))+
2\overline{\eta}(G'')F'+\overline{\eta}(G')F''\Big)\alpha^{\frac{3}{2}}.
\end{array}\right.\end{equation*}\end{theorem}

\begin{proof} First, we recall the structure of the space ${\rm H}^1_\mathrm{diff}(\frak
{osp}(1|2);\mathfrak{D}_{\lambda,\mu})$ computed in \cite{bb}:
\begin{equation}\label{hosp}
{\rm H}^1_\mathrm{diff}({\frak
{osp}}(1|2),{\frak{D}}_{\lambda,\mu})\simeq\left\{
\begin{array}{ll}
\mathbb{K}&\makebox{ if }~~\lambda=\mu, \\[2pt]
\mathbb{K}^2 & \hbox{ if }~~\lambda=\frac{1-k}{2},~\mu=\frac{k}{2},~k\in\mathbb{N}\setminus\{0\},\\[2pt]
0&\makebox { otherwise. }
\end{array}
\right.
\end{equation}

Note that ${\mathrm H}^1_{\mathrm{diff}}(\mathcal{K}(1),\frak{
osp}(1|2);\frak{D}_{\lambda,\mu})\subset{\mathrm
H}^1_{\mathrm{diff}}(\mathcal{K}(1),\frak{D}_{\lambda,\mu}).$
Moreover, if $\mu\neq\lambda$, then by (\ref{hosp}) and Lemma
\ref{sl2} we can see that ${\mathrm
H}^1_{\mathrm{diff}}(\mathcal{K}(1),\frak{
osp}(1|2);\frak{D}_{\lambda,\mu})={\mathrm
H}^1_{\mathrm{diff}}(\mathcal{K}(1);\frak{D}_{\lambda,\mu})$,
except for
\begin{equation}
\label{sin} (\lambda,\mu)\in\left\{
\begin{array}{ll}
(0,\frac{1}{2}),(-\frac{1}{2},1), (-1,\frac{3}{2}),
(-\frac{3}{2},2),(-2,\frac{5}{2})
\end{array}
\right\}.
\end{equation}
Indeed, let $\Upsilon$ be any non-trivial element of
$Z_{\mathrm{diff}}^1(\mathcal{K}(1),\mathfrak{D}_{\lambda,\mu})$
where $\mu\neq\lambda$. If
$(\lambda,\mu)\neq(\frac{1-k}{2},\frac{k}{2})$ where
$k\in\mathbb{N}\setminus\{0\}$ then, by (\ref{hosp}), we can see
that $\Upsilon_{|\mathfrak{osp}(1|2)}$ is trivial, therefore, we
deduce by using Lemma \ref{osp} that the 1-cocycle $\Upsilon$
defines a non-trivial cohomology class in ${\mathrm
H}^1_{\mathrm{diff}}({\cal K}(1),\frak{
osp}(1|2);\frak{D}_{\lambda,\mu})$. If
$(\lambda,\mu)=(\frac{1-k}{2},\frac{k}{2})$ where
$k\in\mathbb{N}\setminus\{0\}$ then, by (\ref{CohSpace2}), we can
see that, up to a coboundary, generically the 1-cocycle $\Upsilon$
vanishes on $\mathfrak{vect}(1)$ and then we conclude  by using
Lemma \ref{sl2} since $\mathfrak{sl}(2)\subset\mathfrak{vect}(1).$

Thus, we need to study only the case $\mu=\lambda$ together with
the singular cases (\ref{sin}). According to Proposition
\ref{decom}, if $\mu-\lambda$ is integer, then
\begin{equation*}\mathrm{H}^1_{\rm diff}\left(\mathfrak{vect}(1);\mathfrak{D}_{\lambda,\mu}\right)
\simeq\mathrm{H}^1_{\rm diff}\left(\mathfrak{vect}(1);
\mathrm{D}_{\lambda,\mu}\right)\oplus\mathrm{H}^1_{\rm
diff}\left(\mathfrak{vect}(1);
\mathrm{D}_{\lambda+\frac{1}{2},\mu+\frac{1}{2}}\right),\end{equation*}
and if $\mu-\lambda$ is semi-integer, then
\begin{equation*}\mathrm{H}^1_{\rm
diff}\left(\mathfrak{vect}(1);\mathfrak{D}_{\lambda,\mu}\right)
\simeq\mathrm{H}^1_{\rm diff}\left(\mathfrak{vect}(1);\Pi(
\mathrm{D}_{\lambda+\frac{1}{2},\mu})\right)
\oplus\mathrm{H}^1_{\rm diff}\left(\mathfrak{vect}(1);\Pi(
\mathrm{D}_{\lambda,\mu+\frac{1}{2}})\right).\end{equation*} Thus,
we deduce $\mathrm{H}^1_{\rm
diff}(\mathfrak{vect}(1);\mathfrak{D}_{\lambda,\mu})$ from
(\ref{CohSpace2}).

Now, let $\Upsilon$ be a  $1$-cocycle from $\mathcal{K}(1)$ to
$\frak{D}_{\lambda,\lambda},$ that is, $\Upsilon$ is even. The map
$\Upsilon_{|\mathfrak{vect}(1)}$ is a 1-cocycle of
$\mathfrak{vect}(1)$. So, up to a coboundary, we have (here $\alpha,
\beta\in\mathbb{K}$)
\begin{equation}\label{decomp5}
\Phi_{\lambda,\lambda}\circ\Upsilon_{|\mathfrak{vect}(1)}=\alpha
C_{\lambda,\lambda}+ \beta
C_{\lambda+{1\over2},\lambda+{1\over2}}.
\end{equation}
By Lemma \ref{sa}, the 1-cocycle $\Upsilon$ is non-trivial if and
only if $(\alpha,\beta)\neq(0,0)$. By Lemma \ref{sd}, we can write
\begin{equation*}
\Upsilon(X_{h\theta})=\sum_{m,k}b_{m,k}h^{(k)}\theta\partial^m_x+
\sum_{m,k}\widetilde{b}_{m,k}h^{(k)}\partial_{\theta}\partial^m_x,
\end{equation*} where the coefficients $b_{m,k}$ and $\widetilde{b}_{m,k}$
are constants. Moreover, the map $\Upsilon$ must satisfy the
following equations
\begin{equation}\label{decomp6}\left\{\begin{array}{lllll}
\Upsilon([X_g,\,X_{
h\theta}])&=&\mathfrak{L}_{X_g}^{\lambda,\lambda}\Upsilon(X_{h\theta
})
-\mathfrak{L}_{X_{h\theta }}^{\lambda,\lambda}\Upsilon(X_g), \\
\Upsilon([X_{ h_1\theta},\,X_{h_2\theta }])&=&\mathfrak{L}_{X_{
h_1\theta}}^{\lambda,\lambda}\Upsilon(X_{h_2\theta })
+\mathfrak{L}_{X_{
h_2\theta}}^{\lambda,\lambda}\Upsilon(X_{h_1\theta }).
\end{array}\right.\end{equation}
We solve the equations (\ref{decomp5}) and (\ref{decomp6}) for
$\alpha,\,\beta,\,b_{k,m},\,\widetilde{b}_{m,k}$. We prove that
$\mathrm{H}^1_{\rm diff} (\mathcal{K}(1) ;
\frak{D}_{\lambda,\lambda})$ is spanned by the non-trivial cocycle
$\Upsilon_{\lambda,\lambda}$ corresponding to the cocycle
\begin{equation*}
\Phi_{\lambda,\lambda}^{-1}\circ
\left(C_{\lambda,\lambda}+C_{\lambda+{1\over2},\lambda+{1\over2}}\right)
\end{equation*}
via its restriction to $\mathfrak{vect}(1)$, see (\ref{cocycles}).

 For the singular cases (\ref{sin}), by the same arguments as above,
we get:
\begin{itemize}
  \item [i)] $\mathrm{H}^1_{\rm diff} (\mathcal{K}(1) ; \frak{D}_{0,{1\over2}})$
is spanned  by the non-trivial cocycles $\Upsilon_{0,{1\over2}}$
and $\widetilde{\Upsilon}_{0,{1\over2}}$ corresponding,
respectively, to the cocycles $\Phi_{0,{1\over2}}^{-1}\circ
\Pi\circ (-C_{0,1})$ and $ \Phi_{0,{1\over2}}^{-1}\circ\Pi\circ
(C_{{1\over2},{1\over2}}-\widetilde{C}_{0,1})$, via their
restrictions to $\mathfrak{vect}(1)$.
  \item [ii)] $\mathrm{H}^1_{\rm diff} (\mathcal{K}(1) ;
\frak{D}_{-{1\over2},1})$ is spanned  by the non-trivial cocycle
$\Upsilon_{-{1\over2},1}$ corresponding to the cocycle
$\Phi_{-{1\over2},1}^{-1}\circ\Pi\circ( C_{0,1}-
C_{-{1\over2},\frac{3}{2}})$ via its restriction to
$\mathfrak{vect}(1).$
\item[iii)] $\mathrm{H}^1_{\rm
diff} (\mathcal{K}(1) ; \frak{D}_{-1,\frac{3}{2}})$ is spanned  by
the non-trivial cocycle $\Upsilon_{-1,\frac{3}{2}}$ corresponding
to the cocycle $ \Phi_{-1,\frac{3}{2}}^{-1}\circ\Pi\circ(
C_{-{1\over2},\frac{3}{2}}-3C_{-1,2})$ via its restriction to
$\mathfrak{vect}(1)$.
\item[iv)] $\mathrm{H}^1_{\rm
diff} (\mathcal{K}(1) ;
\frak{D}_{-\frac{3}{2},2})=\mathrm{H}^1_{\rm diff} (\mathcal{K}(1)
; \frak{D}_{-2,\frac{5}{2}})=0$.
\end{itemize}
\end{proof}
\section{Deformation Theory and Cohomology}
Deformation theory of Lie algebra homomorphisms was first
considered with only one-parameter of deformation \cite{fi, nr,
r}. Recently, deformations of Lie (super)algebras with
multi-parameters were intensively studied ( see,  e.g.,
\cite{aalo, abbo, bbdo, bb, or1, or2}). Here we give an outline of
this theory.
\subsection{Infinitesimal deformations and the first cohomology}
Let $\rho_0 :\frak g \longrightarrow{\rm End}(V)$ be an action of
a Lie superalgebra $\frak g$ on a vector superspace $V$ and let
$\frak h$ be a subagebra of $\frak g$. When studying $\frak
h$-trivial deformations of the $\frak g$-action $\rho_0$, one
usually starts with {\it infinitesimal} deformations:
\begin{equation}\label{infdef}
\rho=\rho_0+t\,\Upsilon,
\end{equation}
where $\Upsilon:\frak g\to{\rm End}(V)$ is a linear map vanishing
on $\frak h$ and $t$ is a formal parameter with
$p(t)=p(\Upsilon)$. The homomorphism condition
\begin{equation}\label{homocond}
[\rho(x),\rho(y)]=\rho([x,y]),
\end{equation}
where $x,y\in\frak g$, is satisfied in order 1 in $t$ if and only
if $\Upsilon$ is a $\frak h$-relative 1-cocycle. That is, the map
$\Upsilon$ satisfies
\begin{equation*}
(-1)^{p(x)p(\Upsilon)}[\rho_0(x),\Upsilon(y)]-(-1)^{p(y)(p(x)+p(\Upsilon))}[\rho_0(y),
\Upsilon(x)]-\Upsilon([x,~y])=0.
\end{equation*}
Moreover, two $\frak h$-trivial infinitesimal deformations $
\rho=\rho_0+t\,\Upsilon_1, $ and $ \rho=\rho_0+t\,\Upsilon_2, $
are equivalents if and only if $\Upsilon_1-\Upsilon_2$ is $\frak
h$-relative coboundary:
\begin{equation*}(\Upsilon_1-\Upsilon_2)(x)=(-1)^{p(x)p(A)}[\rho_0(x),A]:=\delta
A(x),
\end{equation*}
where $A\in{\rm End}(V)^{\frak h}$ and $\delta$ stands for
differential of cochains on $\frak g$ with values in
$\mathrm{End}(V)$ (see, e.g., \cite{Fu, nr}). So, the space
$\mathrm{H}^1(\frak g,\frak h;{\rm End}(V))$ determines and
classifies infinitesimal deformations up to equivalence. If
$\dim{\mathrm{H}^1(\frak g,\frak h;{\rm End}(V))}=m$, then choose
1-cocycles $\Upsilon_1,\ldots,\Upsilon_m$ representing a basis of
$\mathrm{H}^1(\frak g,\frak h;{\rm End(V)})$ and consider the
infinitesimal deformation
\begin{equation}
\label{InDefGen2} \rho=\rho_0+\sum_{i=1}^m{}t_i\,\Upsilon_i,
\end{equation}
where $t_1,\ldots,t_m$ are independent parameters with
$p(t_i)=p(\Upsilon_i)$.

Since we are interested in the $\mathfrak{osp}(1|2)$- trivial
deformations of the $\mathcal{K}(1)$-action on ${\frak
S}^n_{\beta}$, we consider the space ${\mathrm H}^1_{\rm
diff}(\mathcal{K}(1),\mathfrak{osp}(1|2);\mathrm{End}({\frak
S}^n_{\beta}))$ spanned by the classes
$\Upsilon_{\lambda,\lambda+\frac{k}{2}}$, where $k=3,\,4,\,5$ and
$2(\beta-\lambda)\in\left\{k,\,k+1,\, \dots,\,2n\right\}$ for
generic $\beta$. Any infinitesimal $\mathfrak{osp}(1|2)$-trivial
deformation of the $\mathcal{K}(1)$-module ${\frak S}^n_{\beta}$
is then of the form
\begin{equation}\label{infdef1}
\widetilde{\frak L}_{X_F}=\frak{L}_{X_F}+{\frak L}^{(1)}_{X_F},
\end{equation}
where $\frak{L}_{X_F}$ is the Lie derivative of ${\frak
S}^n_{\beta}$ along the vector field $X_F$ defined by
(\ref{superaction}), and
\begin{equation}\label{infpart}
{\frak L}_{X_F}^{(1)}=\sum_{\lambda}\sum_{k=3,4,5}t_{\lambda,
\lambda+\frac{k}{2}}\,
\Upsilon_{\lambda,\lambda+\frac{k}{2}}(X_F),
\end{equation}
where the  $t_{\lambda,\lambda+\frac{k}{2}}$ are independent
parameters with
$p(t_{\lambda,\lambda+\frac{k}{2}})=p(\Upsilon_{\lambda,\lambda+\frac{k}{2}})$
and $2(\beta-\lambda)\in\left\{k,\,k+1,\, \dots,\,2n\right\}.$

\subsection{Integrability conditions and deformations\\ over
supercommutative algebras}

Consider the supercommutative associative superalgebra with unity
$\mathbb{C}[[t_1,\ldots,t_m]]$ and consider the problem of
integrability of infinitesimal deformations. Starting with the
infinitesimal deformation (\ref{InDefGen2}), we look for a formal
series
\begin{equation}
\label{BigDef2} \rho= \rho_0+\sum_{i=1}^m{}t_i\,\Upsilon_i+
\sum_{i,j}{}t_it_j\,\rho^{(2)}_{ij}+\cdots,
\end{equation}
where the higher order terms
$\rho^{(2)}_{ij},\rho^{(3)}_{ijk},\ldots$ are linear maps from
$\frak g$ to ${\rm End(V)}$ with $p(\rho^{(2)}_{ij})=p(t_it_j),
p(\rho^{(3)}_{ijk})=p(t_it_jt_k),\dots$ such that the map
\begin{equation} \label{map} \rho:\frak g\to
\mathbb{C}[[t_1,\ldots,t_m]]\otimes{\rm End(V)},
\end{equation}
satisfies the homomorphism condition (\ref{homocond}).

Quite often the above problem has no solution. Following \cite{fi}
and \cite{aalo}, we will impose extra algebraic relations on the
parameters $t_1,\ldots,t_m$. Let ${\cal R}$ be an ideal in
$\mathbb{C}[[t_1,\ldots,t_m]]$ generated by some set of relations,
and we can speak about  deformations with base ${\cal
A}=\mathbb{C}[[t_1,\ldots,t_m]]/{\cal R}$, (for details, see
\cite{fi}). The map (\ref{map}) sends $\frak g$ to ${\cal
A}\otimes{\rm End}(V)$.

Setting
\begin{equation*}\varphi_t = {\rho}- \rho_0,\,\,
\rho^{(1)}=\sum {}t_i\,\Upsilon_i,\,\,
\rho^{(2)}=\sum{}t_it_j\,\rho^{(2)}_{ij},\,\dots,
\end{equation*}
we can rewrite the relation (\ref{homocond}) in the following way:
\begin{equation}
\label{developping} [\varphi_t(x) , \rho_0(y) ] + [\rho_0(x) ,
\varphi_t(y) ] - \varphi_t([x , y]) +\sum_{i,j > 0}
\;[\rho^{(i)}(x) , \rho^{(j)}(y)] = 0.
\end{equation}
The first three terms are $(\delta\varphi_t) (x,y)$. For arbitrary
linear maps $\gamma_1,~ \gamma_2 :\mathfrak{g}
\longrightarrow\mathrm{End}(V)$, consider   the standard {\it
cup-product}: $[\![\gamma_1,\gamma_2]\!]:\frak g \otimes \frak g
\longrightarrow \mathrm{End}(V)$ defined by:
\begin{equation}
\label{maurrer cartan1} [\![\gamma_1 , \gamma_2]\!] (x , y) =
(-1)^{p(\gamma_2)(p(\gamma_1)+p(x))}[\gamma_1(x) , \gamma_2(y)] +
(-1)^{p(\gamma_1)p(x)}[\gamma_2(x) , \gamma_1(y)].
\end{equation}

The relation (\ref{developping}) becomes now equivalent to:
\begin{equation}
\label{maurrer cartan} \delta\varphi_t +{1\over2} [\![\varphi_t ,
\varphi_t]\!]= 0,
\end{equation}
Expanding (\ref{maurrer cartan}) in power series in $t_1,\dots,t_m
$, we obtain the following equation for $\rho^{(k)}$:
\begin{equation}
\label{maurrer cartank} \delta\rho^{(k)} + {1\over2}\sum_{i+j=k}
[\![\rho^{(i)} ,  \rho^{(j)}]\!] = 0.
\end{equation}

The first non-trivial relation $\delta{\rho^{(2)}} +{1\over2}
[\![\rho^{(1)} , \rho^{(1)}]\!] = 0 $ gives the first obstruction
to integration of an infinitesimal deformation. Thus, considering
the coefficient of $t_i\,t_j$, we  get
\begin{equation}
\label{nizar} \delta{\rho^{(2)}_{ij}}+{1\over2}[\![\Upsilon_{i} ,
\Upsilon_{j}]\!]=0.
\end{equation}
It is easy to check that for any two $1$-cocycles $\gamma_1$ and
$\gamma_2 \in Z^1 (\frak g ,\frak h; \mathrm{End}(V))$, the
bilinear map $[\![\gamma_1 , \gamma_2]\!]$ is a $\frak h$-relative
$2$-cocycle. The relation (\ref{nizar}) is precisely the condition
for this cocycle to be a coboundary. Moreover, if one of the
cocycles $\gamma_1$ or $\gamma_2$ is a $\frak h$-relative
coboundary, then $[\![\gamma_1 , \gamma_2]\!]$ is a $\frak
h$-relative $2$-coboundary. Therefore, we naturally deduce that
the operation (\ref{maurrer cartan1}) defines a bilinear map:
\begin{equation}
\label{cup-product} \mathrm{H}^1 (\frak g ,\frak h;\mathrm{End}(
V))\otimes \mathrm{H}^1 (\frak g ,\frak h; \mathrm{End}(
V))\longrightarrow \mathrm{H}^2 (\frak g ,\frak h; \mathrm{End}(
V)).
\end{equation}

All the obstructions lie in $\mathrm{H}^2 (\frak g, \frak
h;\mathrm{End}(V))$ and they are in the image of $\mathrm{H}^1
(\frak g,\frak h; \mathrm{End}(V))$ under the cup-product.

\subsection{Equivalence }

Two deformations, $\rho$ and $\rho'$ of a~$\frak g$-module $V$
over $\cal A$ are said to be {\it equivalent} (see \cite{fi}) if
there exists an inner automorphism $\Psi$ of the associative
superalgebra ${\cal A}\otimes{\rm End}(V)$ such that
\begin{equation*}
\Psi\circ\rho=\rho'\hbox{ and } \Psi(\mathbb{I})=\mathbb{I},
\end{equation*}
where $\mathbb{I}$ is the unity of the superalgebra ${\cal
A}\otimes{\rm End}(V).$

The following notion of miniversal deformation is fundamental. It
assigns to a~$\frak g$-module $V$ a canonical commutative
associative algebra $\mathcal{A}$ and a canonical deformation over
$\mathrm{A}$. A deformation (\ref{BigDef2}) over $\mathcal{A}$ is
said to be {\it miniversal} if
\begin{itemize}
  \item [(i)] for any other deformation $\rho'$ with base (local)
$\cal A'$, there exists a  homomorphism $\psi:{\cal A}'\to{\cal
A}$ satisfying $\psi(1)=1$, such that
\begin{equation*}
\rho=(\psi\otimes\mathrm{Id})\circ\rho'.
\end{equation*}
  \item [(ii)] under notation of (i), if $\rho$ is infinitesimal,
  then  $\psi$ is unique.
\end{itemize} If $\rho$ satisfies only the condition (i), then it is called
versal.   This definition does not depend on the choice 1-cocycles
$\Upsilon_1,\ldots,\Upsilon_m$ representing a basis of
$\mathrm{H}^1(\frak g,\frak h;{\rm End(V)})$.

The miniversal deformation corresponds to the smallest ideal
$\mathcal{R}$. We refer to \cite{fi} for a construction of
miniversal deformations of Lie algebras and to \cite{aalo} for
miniversal deformations of $\mathfrak{g}$-modules. Superization of
these results is immediate: by the Sign Rule.
\section{ Integrability Conditions}

In this section we obtain the integrability conditions for the
infinitesimal deformation(\ref{infdef1}).

\begin{proposition}\label{th2} The second-order integrability conditions of the
infinitesimal deformation~(\ref{infdef1}) are the following:
\begin{equation}
\label{ord2} t_{\lambda,
\lambda+\frac{5}{2}}\,t_{\lambda+\frac{5}{2},
\lambda+5}=0,\quad\text{where}\quad
2(\beta-\lambda)\in\left\{10,\, \dots,\,2n\right\}.
\end{equation}
\end{proposition}
To prove Proposition \ref{th2}, we need the following lemmas
\begin{lemma}\label{trans1} Consider a linear differential operator
$b:\mathcal{K}(1)\longrightarrow\frak{D}_{\lambda,\mu}$. If $b$
satisfies
\begin{equation*}
\delta(b)(X,Y)=b(X)=0\hbox{ for all }{X\in\frak {osp}}(1|2),
\end{equation*} then $b$ is a
supertransvectant.
\end{lemma}

\begin{proof} For all $X,\,Y\in\mathcal{K}(1)$ we have
\begin{equation*}
\delta(b)(X,Y):=(-1)^{p(X)p(b)}\mathfrak{L}^{\lambda,\mu}_{X}(b(Y))-
(-1)^{p(Y)(p(X)+p(b))}\mathfrak{L}^{\lambda,\mu}_{Y}(b(X))-b([X,Y]).
\end{equation*}
Since $ \delta(b)(X,Y)=b(X)=0\hbox{ for all }X\in\frak {osp}(1|2)
$ we deduce that \begin{equation*}
(-1)^{p(X)p(b)}\mathfrak{L}^{\lambda,\mu}_{X}(b(Y))-b([X,Y])=0.
\end{equation*}
Thus, the map ${b}$ is $\frak{ osp}(1|2)$-invariant.
\end{proof}

\medskip

\begin{lemma}\label{lth2}The map $B_{\lambda,\lambda+5}=[\![\Upsilon_{\lambda+\frac{5}{2},\lambda+5}
,\, \Upsilon_{\lambda,\lambda+\frac{5}{2}} ]\!]$ is a non-trivial
$\frak{osp}(1|2)$-relative 2-cocycle  for $\lambda\neq
0,\,-1,\,-\frac{7}{2},\,-\frac{9}{2}$.
\end{lemma}
\begin{proof} First, observe that for $\lambda=-1, -\frac{7}{2}$, the map
$B_{\lambda,\lambda+5}$ is not defined (see Theorem \ref{th3}).
The map $B_{\lambda,\lambda+5}$ is the cup-product of two
$\frak{osp}(1|2)$-relative 1-cocycles, so, $B_{\lambda,\lambda+5}$
is a $\frak{osp}(1|2)$-relative 2-cocycle:
$B_{\lambda,\lambda+5}\in
\mathrm{Z}^2(\mathcal{K}(1),\frak{osp}(1|2);\frak{D}_{\lambda,\lambda+5})$.
This 2-cocycle is trivial if and only if it is the coboundary of a
linear differential operator
\begin{equation*}b_{\lambda,\lambda+5}:\mathcal{K}(1)\longrightarrow
\frak{D}_{\lambda,\lambda+5}\end{equation*} vanishing on
$\mathfrak{osp}(1|2)$. Consider $b_{\lambda,\lambda+5}$ as a
bilinear map
$\mathfrak{F}_{-1}\otimes\mathfrak{F}_\lambda\longrightarrow\mathfrak{F}_{\lambda+5}$.
So, according to Lemma \ref{trans1} and Theorem \ref{main}, the
operator ${b}_{\lambda,\lambda+5}$ coincides (up to a scalar
factor) with the supertransvectant $\frak{J}_{6}^{-1,\lambda}$.
But, by a direct computation, we have, up to a multiple
\begin{equation*}
\begin{array}{llll}
B_{\lambda,\lambda+5}(X_{g_1},\,X_{g_2})(F\alpha^\lambda)&=&
\left(g_1^{(4)}g_2^{(3)}-g_1^{(3)}g_2^{(4)}\right)\left(2\lambda
f_0-(2\lambda+9)f_1\theta\right)\alpha^{\lambda+5},\\[10pt]
\delta(\frak{J}_{6}^{-1,\lambda})(X_{g_1},\,X_{g_2})(F\alpha^\lambda)&=&
\left(g_1^{(3)}g_2^{(4)}-g_1^{(4)}g_2^{(3)}\right)
\Big(\frac{\lambda(2\lambda+3)(\lambda^2+6\lambda+8)}{9}f_0+\\[10pt]
&&+\frac{(2\lambda+1)(\lambda+3)(4\lambda^2+28\lambda+45)}{36}f_1\theta\Big)\alpha^{\lambda+5},
\end{array}
\end{equation*} where $g_1, g_2\in\mathbb{K}[x]$ and
$F=f_0+f_1\theta\in\mathbb{K}[x,\theta]$. Therefore, the
restrictions of the maps $B_{\lambda,\lambda+5}$ and
$\delta(\frak{J}_{6}^{-1,\lambda})$ to
 $\mathfrak{vect}(1)\times \mathfrak{vect}(1)$ are linearly dependant if and only if
\begin{equation*}
\lambda(\lambda+1)(2\lambda+7)(2\lambda+9)(4\lambda+9)=0.
\end{equation*}
 Thus, the maps
$B_{\lambda,\lambda+5}$ and $\delta(\frak{J}_{6}^{-1,\lambda})$
are linearly independent for
$\lambda\neq0,-1,\,-\frac{7}{2},\,-\frac{9}{2},\,-\frac{9}{4}$.
Besides, we check that the maps $B_{-\frac{9}{4},\frac{11}{4}}$
and $\delta(\frak{J}_{6}^{-1,-\frac{9}{4}})$ are also linearly
independent although their restrictions to
$\mathfrak{vect}(1)\times \mathfrak{vect}(1)$ are linearly
dependant. Finally, for $\lambda=0,-\frac{9}{2}$, we check that
$B_{\lambda,\lambda+5}$ coincides (up to a scalar factor) with
$\delta(\frak{J}_{6}^{-1,\lambda})$. This completes the proof.
\end{proof}

\begin{remark}{\rm The map
$B_{\lambda,\lambda+5}:\mathcal{K}(1)\otimes\mathcal{K}(1)
\longrightarrow\frak{D}_{\lambda,\lambda+5}$ is a non-trivial
2-cocycle, so, $\mathrm{H}^2_{\mathrm{diff}}(\mathcal{K}(1),\frak{
osp}(1|2);\frak{D}_{\lambda,\lambda+5})\neq 0$ while  ${\mathrm
H}^2_{\rm diff}(\mathfrak{vect}(1),{\rm \frak{sl}}(2);{\rm
D}_{\lambda,\lambda+5})=0$ for generic $\l$ (see \cite{b}). Hence,
for the second cohomology, the analog of Lemma \ref{sa} is not
true.}
\end{remark}
\begin{proof} of Proposition \ref{th2}: Assume
that the infinitesimal deformation (\ref{infdef1}) can be
integrated to a formal deformation:
\begin{equation*}
\widetilde{\frak L}_{X_F}=\frak{L}_{X_F}+{\frak
L}^{(1)}_{X_F}+{\frak L}^{(2)}_{X_F}+\cdots
\end{equation*}
The homomorphism  condition gives, for the term ${\frak
L}^{(2)}_{\lambda,\mu,\lambda',\mu'}$ in
$t_{\lambda,\mu}t_{\lambda',\mu'}$, the following equation
\begin{equation}\label{cap}\delta({\frak
L}^{(2)}_{\lambda,\mu,\lambda',\mu'})=-[\![\Upsilon_{\lambda,\mu},\,
\Upsilon_{\lambda',\mu'}]\!].\end{equation}
For arbitrary $\lambda$, the right hand side of (\ref{cap}) yields
the following  2-cocycles:
\begin{equation}\label{2-cocyc}\begin{array}{llllll}
B_{\lambda,\lambda+3}&=&[\![\Upsilon_{\lambda+\frac{3}{2},\lambda+3},\,
\Upsilon_{\lambda,\lambda+\frac{3}{2}}
]\!]&:&\mathcal{K}(1)\otimes\mathcal{K}(1)\rightarrow\frak{D}_{\lambda,\lambda+3},\\
B_{\lambda,\lambda+{7\over2}}&=&[\![\Upsilon_{\lambda+\frac{3}{2},\lambda+{7\over2}},\,
\Upsilon_{\lambda,\lambda+\frac{3}{2}}
]\!]&:&\mathcal{K}(1)\otimes\mathcal{K}(1)\rightarrow\frak{D}_{\lambda,\lambda+{7\over2}},\\
\widetilde{B}_{\lambda,\lambda+{7\over2}}&=&[\![\Upsilon_{\lambda+2,\lambda+{7\over2}},\,
\Upsilon_{\lambda,\lambda+2}
]\!]&:&\mathcal{K}(1)\otimes\mathcal{K}(1)\rightarrow\frak{D}_{\lambda,\lambda+{7\over2}},\\
B_{\lambda,\lambda+4}&=&[\![\Upsilon_{\lambda+\frac{3}{2},\lambda+4}
,\, \Upsilon_{\lambda,\lambda+\frac{3}{2}}
]\!]&:&\mathcal{K}(1)\otimes\mathcal{K}(1)\rightarrow\frak{D}_{\lambda,\lambda+4},\\
\widetilde{B}_{\lambda,\lambda+4}&=&[\![\Upsilon_{\lambda+\frac{5}{2},\lambda+4}
,\, \Upsilon_{\lambda,\lambda+\frac{5}{2}}
]\!]&:&\mathcal{K}(1)\otimes\mathcal{K}(1)\rightarrow\frak{D}_{\lambda,\lambda+4},\\
\overline{B}_{\lambda,\lambda+4}&=&[\![\Upsilon_{\lambda+2,\lambda+4}
,\, \Upsilon_{\lambda,\lambda+2}
]\!]&:&\mathcal{K}(1)\otimes\mathcal{K}(1)\rightarrow\frak{D}_{\lambda,\lambda+4},\\
B_{\lambda,\lambda+{9\over2}}&=&[\![\Upsilon_{\lambda+\frac{5}{2},\lambda+{9\over2}}
,\, \Upsilon_{\lambda,\lambda+\frac{5}{2}}
]\!]&:&\mathcal{K}(1)\otimes\mathcal{K}(1)\rightarrow\frak{D}_{\lambda,\lambda+{9\over2}},\\
\widetilde{B}_{\lambda,\lambda+{9\over2}}&=&[\![\Upsilon_{\lambda+2,\lambda+{9\over2}}
,\, \Upsilon_{\lambda,\lambda+2}
]\!]&:&\mathcal{K}(1)\otimes\mathcal{K}(1)\rightarrow\frak{D}_{\lambda,\lambda+{9\over2}},\\
B_{\lambda,\lambda+5}&=&[\![\Upsilon_{\lambda+\frac{5}{2},\lambda+5}
,\, \Upsilon_{\lambda,\lambda+\frac{5}{2}}
]\!]&:&\mathcal{K}(1)\otimes\mathcal{K}(1)\rightarrow\frak{D}_{\lambda,\lambda+5}.\end{array}
\end{equation}

The necessary integrability conditions for the second-order terms
${\frak L}^{(2) }$  are that each 2-cocycle
$B_{\lambda,\lambda+k}$, where $2k=6,\,7,\,8,\,9,\,10$, must be a
coboundary of a linear differential operator
$b_{\lambda,\lambda+k}:\mathcal{K}(1)\longrightarrow
\frak{D}_{\lambda,\lambda+k},$ vanishing on $\mathfrak{osp}(1|2)$.
More precisely, as in the proof of Lemma \ref{lth2}, the operator
${b}_{\lambda,\lambda+k}$ coincides (up to a scalar factor) with
the supertransvectant $\frak{J}_{k+1}^{-1,\lambda}$. Clearly,
\begin{equation*}
\widetilde{B}_{\lambda,\lambda+4}=
B_{\lambda,\lambda+4}=3\overline{B}_{\lambda,\lambda+4},\quad
{B}_{\lambda,\lambda+{7\over2}}=-\widetilde{B}_{\lambda,\lambda+{7\over2}},\quad
{B}_{\lambda,\lambda+{9\over2}}=-\widetilde{B}_{\lambda,\lambda+{9\over2}}
\end{equation*}
and, by a direct computation, we check that
\begin{equation*}
\begin{array}{lll}
B_{\lambda,\lambda+3}(X_{G_1},X_{G_2})(F\alpha^\lambda)&=&\left(-2(-1)^{p(G_1)}\overline{\eta}(G_1'')\overline{\eta}(G_2'')F\right)\alpha^{\lambda+3},\\[10pt]
{B}_{\lambda,\lambda+{7\over2}}(X_{G_1},X_{G_2})(F\alpha^\lambda)&=&
\Big({2\lambda\over3}
\left((-1)^{p(G_1)}G_1^{(3)}\overline{\eta}(G_2'')-\overline{\eta}(G_1'')G_2^{(3)}\right)F+\\
&&2(-1)^{p(G_2)}\overline{\eta}(G_1'')\overline{\eta}(G_2'')\overline{\eta}(F)\Big)
\alpha^{\lambda+{7\over2}},\\[10pt]
B_{\lambda,\lambda+4}(X_{G_1},X_{G_2})(F\alpha^\lambda)&=&
\Big(-2\lambda(-1)^{p(G_1)}\left(\overline{\eta}(G_1^{(3)})\overline{\eta}(G_2'')+
\overline{\eta}(G_1'')\overline{\eta}(G_2^{(3)})\right)F+\\
&&(-1)^{p(G_2)}\left((-1)^{p(G_1)}\overline{\eta}(G_1'')G_2^{(3)}-G_1^{(3)}
\overline{\eta}(G_2'')\right)\overline{\eta}(F)+\\
&&6(-1)^{p(G_1)}\overline{\eta}(G_1'')\overline{\eta}(G_2'')F'\Big)
\alpha^{\lambda+4},\\[10pt]
{B}_{\lambda,\lambda+{9\over2}}(X_{G_1},X_{G_2})(F\alpha^\lambda)&=&\Big(
{2\lambda\over3}(2\lambda+5)\left((-1)^{p(G_1)}G_1^{(3)}\overline{\eta}(G_2^{(3)})-
\overline{\eta}(G_1^{(3)})G_2^{(3)}\right)F+\\ &&
2\lambda\left(\overline{\eta}(G_1'')G_2^{(4)}-(-1)^{p(G_1)p(G_2)}
\overline{\eta}(G_2'')G_1^{(4)}\right)F+\\
&&(2\lambda+1)(-1)^{p(G_2)}\left(\overline{\eta}(G_1'')\overline{\eta}(G_2^{(3)})+
\overline{\eta}(G_1^{(3)})\overline{\eta}(G_2'')\right)\overline{\eta}(F)+\\
&&(2\lambda+1)\left(\overline{\eta}(G_1'')G_2^{(3)}-(-1)^{p(G_1)}G_1^{(3)}
\overline{\eta}(G_2'')\right)F'-\\
&&6(-1)^{p(G_2)}\overline{\eta}(G_1'')\overline{\eta}(G_2'')
\overline{\eta}(F')\Big)\alpha^{\lambda+{9\over2}},
\end{array}
\end{equation*}
where $G_1,G_2,F\in\mathbb{K}[x,\theta].$ Besides, we can see that
\begin{gather*}
\zeta_\lambda
B_{\lambda,\lambda+3}=\delta(\frak{J}_{4}^{-1,\lambda}),\quad
\text{ where }\quad \zeta_\lambda= {\lambda(2\lambda+5)\over
4}\,\begin{pmatrix}2\lambda+3\\2\end{pmatrix},\\ \alpha_\lambda
B_{\lambda,\lambda+{7\over2}}=\delta(\frak{J}_{9\over2}^{-1,\lambda}),\quad\text{
where }\quad \alpha_\lambda=
{6\lambda+9\over4}\,\begin{pmatrix}2\lambda+4\\3\end{pmatrix},\\
\beta_\lambda
B_{\lambda,\lambda+4}=\delta(\frak{J}_{5}^{-1,\lambda}),\quad\text{
where }\quad \beta_\lambda=
{2\lambda^2+7\lambda+2\over6}\,\begin{pmatrix}2\lambda+4\\2\end{pmatrix},\\
\gamma_\lambda
B_{\lambda,\lambda+{9\over2}}=\delta(\frak{J}_{11\over2}^{-1,\lambda}),\quad\text{
where }\quad \gamma_\lambda= {3\lambda+6\over
2}\,\begin{pmatrix}2\lambda+5\\3\end{pmatrix}.
\end{gather*}

Now, by Lemma \ref{lth2},  $B_{\lambda,\lambda+5}$  is a
non-trivial $\frak{osp}(1|2)$-relative 2-cocycle, so, its
coefficient must vanish, that is, we get the first set of
necessary integrability conditions:
\begin{equation}\label{c1}t_{\lambda,
\lambda+\frac{5}{2}}\,t_{\lambda+\frac{5}{2},
\lambda+5}=0,\quad\text{where}\quad
2(\beta-\lambda)\in\left\{10,\, \dots,\,2n\right\}.\end{equation}
The equations (\ref{c1}) are the unique integrability conditions
for the 2nd order term ${\frak L}^{(2) }$. Under these conditions,
the second-order term ${\frak L}^{(2) }$ can be given by
\begin{equation*}\begin{array}{lll}{\frak L}^{(2) }=
&-\sum_\lambda\zeta_\lambda^{-1}t_{\lambda+\frac{3}{2},\lambda+3}
t_{\lambda,\lambda+\frac{3}{2}} \frak{J}_4^{-1,\lambda}\\[2pt]
&-\sum_\lambda\alpha_\lambda^{-1}(t_{\lambda+\frac{3}{2},\lambda+\frac{7}{2}}
t_{\lambda,\lambda+\frac{3}{2}}-
t_{\lambda+2,\lambda+\frac{7}{2}}t_{\lambda,\lambda+2})
\frak{J}_\frac{9}{2}^{-1,\lambda}\\[2pt]&-\sum_\lambda\beta_\lambda^{-1}(t_{\lambda+\frac{3}{2},\lambda+4}
t_{\lambda,\lambda+\frac{3}{2}}+
t_{\lambda+\frac{5}{2},\lambda+4}t_{\lambda,\lambda+\frac{5}{2}}+{1\over3}t_{\lambda+2,\lambda+4}t_{\lambda,\lambda+2})
\frak{J}_{5}^{-1,\lambda}\\[2pt]
&-\sum_\lambda\gamma_\lambda^{-1}(t_{\lambda+\frac{5}{2},\lambda+\frac{9}{2}}
t_{\lambda,\lambda+\frac{5}{2}}-
t_{\lambda+2,\lambda+\frac{9}{2}}t_{\lambda,\lambda+2})
\frak{J}_\frac{11}{2}^{-1,\lambda}.\end{array}\end{equation*}
\end{proof}

To compute the third term ${\frak L}^{(3) }$, we need the
following two lemmas which we can check by a direct computation
with the help of {\it Maple}.
\begin{lemma}
\label{benfraj1}
\begin{equation*}
\begin{array}{llllll}
1)~\xi_{\lambda}^{-1}\,\delta(\frak{J}_{11\over2}^{-1,\lambda})&=&
\left(\zeta_\lambda^{-1}[\![\Upsilon_{\lambda+3,\lambda+\frac{9}{2}},\
\frak{J}_4^{-1,\lambda}
]\!]+\zeta_{\lambda+\frac{3}{2}}^{-1}\,[\![\frak{J}_4^{-1,\lambda+\frac{3}{2}},
\,\Upsilon_{\lambda,\lambda+\frac{3}{2}}
]\!]\right),\\[8pt]
2)~
\alpha_{\lambda+\frac{3}{2}}^{-1}[\![\frak{J}_{\frac{9}{2}}^{-1,\lambda+\frac{3}{2}},\
\Upsilon_{\lambda,\lambda+\frac{3}{2}}
]\!]&=&\epsilon_{1,\lambda}\,\alpha_{\lambda}^{-1}\,[\![\Upsilon_{\lambda+\frac{7}{2},\lambda+5},\
\frak{J}_{\frac{9}{2}}^{-1,\lambda}
]\!]+\epsilon_{2,\lambda}\,\zeta_{\lambda}^{-1}\,[\![\Upsilon_{\lambda+3,\lambda+5},\
\frak{J}_{4}^{-1,\lambda}
]\!]+\\&&\epsilon_{3,\lambda}\,\delta(\frak{J}_{6}^{-1,\lambda}),\\[8pt]
3)~\beta_{\lambda+\frac{3}{2}}^{-1}[\![\frak{J}_{5}^{-1,\lambda+\frac{3}{2}},\
\Upsilon_{\lambda,\lambda+\frac{3}{2}}
]\!]&=&\epsilon_{4,\lambda}\beta_{\lambda}^{-1}[\![\Upsilon_{\lambda+4,\lambda+\frac{11}{2}},
\frak{J}_{5}^{-1,\lambda}
]\!]+\epsilon_{5,\lambda}\,\alpha_{\lambda}^{-1}\,[\![\Upsilon_{\lambda+\frac{7}{2},\lambda+\frac{11}{2}},\
\frak{J}_{\frac{9}{2}}^{-1,\lambda}
]\!]+\\&&\epsilon_{6,\lambda}\,\alpha_{\lambda+2}^{-1}\,[\![\frak{J}_{\frac{9}{2}}^{-1,\lambda+2},\
\Upsilon_{\lambda,\lambda+2}]\!],\\[8pt]
4)~
\alpha_\lambda^{-1}[\![\Upsilon_{\lambda+\frac{7}{2},\lambda+6},\
\frak{J}_{\frac{9}{2}}^{-1,\lambda}
]\!]&=&\epsilon_{7,\lambda}\,\beta_{\lambda+2}^{-1}\,[\![\frak{J}_5^{-1,\lambda+2},\
\Upsilon_{\lambda,\lambda+2}
]\!]+\epsilon_{8,\lambda}\,\alpha_{\lambda+\frac{5}{2}}^{-1}\,[\![\frak{J}_{\frac{9}{2}}^{-1,\lambda+\frac{5}{2}},\
\Upsilon_{\lambda,\lambda+\frac{5}{2}}
]\!]+\\&&\epsilon_{9,\lambda}\,\delta(\frak{J}_{7}^{-1,\lambda}),\\[8pt]
5)~\gamma_{\lambda+2}^{-1}[\![\frak{J}_{\frac{11}{2}}^{-1,\lambda+2},\
\Upsilon_{\lambda,\lambda+2}
]\!]&=&\epsilon_{10,\lambda}\,\beta_\lambda^{-1}\,[\![\Upsilon_{\lambda+4,\lambda+\frac{13}{2}},\
\frak{J}_{5}^{-1,\lambda}
]\!]+\epsilon_{11,\lambda}\,\beta_{\lambda+\frac{5}{2}}^{-1}[\![\frak{J}_{5}^{-1,\lambda+\frac{5}{2}},\
\Upsilon_{\lambda,\lambda+\frac{5}{2}} ]\!]+\\&&
\epsilon_{12,\lambda}\,\gamma_{\lambda}^{-1}[\![\Upsilon_{\lambda+\frac{9}{2},\lambda+\frac{13}{2}},\
\frak{J}_{\frac{11}{2}}^{-1,\lambda} ]\!],
\end{array}
\end{equation*}
where
\begin{equation*}\small{
\begin{array}{llllll}
\epsilon_{1,\lambda}&=&\frac{(2\lambda+11)(2\lambda+9)(\lambda+2)}{2(\lambda+3)(2\lambda^2+3\lambda-17)},\quad
&\epsilon_{7,\lambda}&=&-\frac{5(6\lambda^2+33\lambda+17)(2\lambda^2+15\lambda+24)}{2(\lambda+5)(\lambda+2)(2\lambda-3)(2\lambda^2+13\lambda+13)},\\[4pt]
\epsilon_{2,\lambda}&=&\frac{15(2\lambda+5)(\lambda+4)}{2(\lambda+3)(2\lambda^2+3\lambda-17)}\quad
&\epsilon_{8,\lambda}&=&\frac{(2\lambda+9)(2\lambda+5)(\lambda+7)(2\lambda^2+11\lambda+4)}{2(\lambda+5)(2\lambda-3)(\lambda+2)(2\lambda^2+13\lambda+13)},\\[4pt]
\epsilon_{3,\lambda}&=&\frac{48}{\lambda(\lambda+3)(2\lambda+3)(2\lambda+5)(2\lambda^2+3\lambda-17)},\quad
&\epsilon_{9,\lambda}&=&\frac{60}{(2\lambda+3)(2\lambda-3)(\lambda+3)(\lambda+4)(2\lambda^2+13\lambda+13)},\\[4pt]
\epsilon_{4,\lambda}&=&-\frac{(2\lambda+9)(2\lambda+3)(2\lambda^2+7\lambda+2)}{(2\lambda+7)(2\lambda+1)(2\lambda^2+13\lambda+17)}\quad
&\epsilon_{10,\lambda}&=&-\frac{(\lambda+5)(\lambda+2)(2\lambda^2+7\lambda+2)}{9(\lambda+4)^2(\lambda+1)},\\[4pt]
\epsilon_{5,\lambda}&=&-\frac{3(2\lambda+9)(2\lambda+3)^2}{2(2\lambda+7)(2\lambda+1)(2\lambda^2+13\lambda+17)}\quad
&\epsilon_{11,\lambda}&=&-\frac{(2\lambda^2+17\lambda+32)}{9(\lambda+4)},\\[4pt]
\epsilon_{6,\lambda}&=&-\frac{3(2\lambda+7)}{2(2\lambda^2+13\lambda+17)}\quad
&\epsilon_{12,\lambda}&=&-\frac{(\lambda+2)^2(\lambda+5)}{(\lambda+4)^2(\lambda+1)},\\[4pt]
\xi_{\lambda}&=& {3\over
16}\lambda(\lambda+4)(2\lambda+3)(2\lambda+5)\,\begin{pmatrix}2\lambda+5\\3\end{pmatrix}.
\end{array}}
\end{equation*}
\end{lemma}
\begin{lemma}
\label{benfraj2} Each of the following systems is linearly
independent
\begin{equation*}
\begin{array}{llllllllll}
1)\,&\left( [\![\Upsilon_{\lambda+\frac{7}{2},\lambda+5},\
\frak{J}_\frac{9}{2}^{-1,\lambda} ]\!],\,\,
[\![\Upsilon_{\lambda+3,\lambda+5},\ \frak{J}_4^{-1,\lambda}
]\!],\,\, [\![\frak{J}_4^{-1,\lambda+2},
\,\Upsilon_{\lambda,\lambda+2} ]\!],\,\,
\delta(\frak{J}_{6}^{-1,\lambda})\right),\\[2pt]
2)\,&\Big([\![\Upsilon_{\lambda+4,\lambda+\frac{11}{2}},\
\frak{J}_{5}^{-1,\lambda} ]\!],\,\,
[\![\Upsilon_{\lambda+\frac{7}{2},\lambda+\frac{11}{2}},\
\frak{J}_\frac{9}{2}^{-1,\lambda} ]\!],\,\,
[\![\frak{J}_\frac{9}{2}^{-1,\lambda+2},
\,\Upsilon_{\lambda,\lambda+2}
]\!],\,\,[\![\Upsilon_{\lambda+3,\lambda+\frac{11}{2}},\
\frak{J}_{4}^{-1,\lambda} ]\!]\\[2pt] &~~
[\![\frak{J}_{4}^{-1,\lambda+\frac{5}{2}},
\,\Upsilon_{\lambda,\lambda+\frac{5}{2}} ]\!],\,\,
 \delta(\frak{J}_{13\over2}^{-1,\lambda})\Big),\\[2pt]
3)\,&\Big([\![\Upsilon_{\lambda+\frac{9}{2},\lambda+6},\
\frak{J}_\frac{11}{2}^{-1,\lambda} ]\!],\,\,[\![
\frak{J}_\frac{11}{2}^{-1,\lambda+\frac{3}{2}},
\,\Upsilon_{\lambda,\lambda+\frac{3}{2}} ]\!],\,\,
[\![\Upsilon_{\lambda+4,\lambda+6},\ \frak{J}_{5}^{-1,\lambda}
]\!],\,\,[\![ \frak{J}_{5}^{-1,\lambda+2},
\,\Upsilon_{\lambda,\lambda+2}
]\!],\,\,\\[2pt]&~~
[\![\frak{J}_\frac{9}{2}^{-1,\lambda+\frac{5}{2}},
\,\Upsilon_{\lambda,\lambda+\frac{5}{2}} ]\!],\,\,
\delta(\frak{J}_{7}^{-1,\lambda})\Big),\\[2pt]
4)\,&\left([\![\Upsilon_{\lambda+4,\lambda+\frac{13}{2}},\
\frak{J}_{5}^{-1,\lambda} ]\!],\, \,[\![
\frak{J}_{5}^{-1,\lambda+\frac{5}{2}},
\,\Upsilon_{\lambda,\lambda+\frac{5}{2}}
]\!],\,\,[\![\Upsilon_{\lambda+{9\over2},\lambda+\frac{13}{2}},\
\frak{J}_{11\over2}^{-1,\lambda} ]\!],\,\,
\delta(\frak{J}_{15\over2}^{-1,\lambda})\right),\\[2pt]
5)\,&\left([\![\Upsilon_{\lambda+{9\over2},\lambda+7},\
\frak{J}_{11\over2}^{-1,\lambda} ]\!],\,\, [\![
\frak{J}_{11\over2}^{-1,\lambda+{5\over2}},
\,\Upsilon_{\lambda,\lambda+{5\over2}} ]\!],\,\,
\delta(\frak{J}_{8}^{-1,\lambda})\right).
\end{array}
\end{equation*}
\end{lemma}
Now, we are in position to exhibit the 3rd order integrability
conditions.
\begin{proposition}\label{pr3} The 3rd order integrability conditions of the
infinitesimal deformation~(\ref{infdef1}) are the following:
\begin{itemize}
\item [a)] For $2(\beta-\lambda)\in\left\{10,\,
\dots,\,2n\right\}:$ \begin{gather*}
t_{\lambda,\lambda+\frac{3}{2}}\left(\epsilon_{1,\lambda}\,
t_{\lambda+3,\lambda+5}
t_{\lambda+\frac{3}{2},\lambda+3}+(1-\epsilon_{1,\lambda})\,
t_{\lambda+\frac{7}{2},\lambda+5}
t_{\lambda+\frac{3}{2},\lambda+\frac{7}{2}}
\right)=0,\notag\\
t_{\lambda,\lambda+\frac{3}{2}}\left(\epsilon_{2,\lambda}\,t_{\lambda+\frac{7}{2},\lambda+5}
t_{\lambda+\frac{3}{2},\lambda+\frac{7}{2}}-(1+\epsilon_{2,\lambda})\,
 t_{\lambda+3,\lambda+5}
t_{\lambda+\frac{3}{2},\lambda+3}\right)=0,\notag\\
t_{\lambda+\frac{7}{2},\lambda+5}t_{\lambda+2,\lambda+\frac{7}{2}}
t_{\lambda,\lambda+2}=0.
\end{gather*}
  \item [b)] For $2(\beta-\lambda)\in\left\{11,\,
\dots,\,2n\right\}:$
\begin{gather*}
t_{\lambda+4,\lambda+\frac{11}{2}}\,t_{\lambda+\frac{5}{2},\lambda+4}
t_{\lambda,\lambda+\frac{5}{2}}= t_{\lambda,\lambda+\frac{3}{2}}
t_{\lambda+3,\lambda+\frac{11}{2}}
t_{\lambda+\frac{3}{2},\lambda+3}=0,\\
t_{\lambda+4,\lambda+\frac{11}{2}}\left(3(1+\epsilon_{4,\lambda})\,t_{\lambda+\frac{3}{2},\lambda+4}
t_{\lambda,\lambda+\frac{3}{2}}+
t_{\lambda+2,\lambda+4}t_{\lambda,\lambda+2}\right)+\\
\quad+\,\epsilon_{4,\lambda}\, t_{\lambda,\lambda+\frac{3}{2}}
t_{\lambda+\frac{7}{2},\lambda+\frac{11}{2}}t_{\lambda+\frac{3}{2},\lambda+\frac{7}{2}}
=0,\\
t_{\lambda+\frac{7}{2},\lambda+\frac{11}{2}}
\left((1+\frac{\epsilon_{5,\lambda}}{3})\,t_{\lambda+\frac{3}{2},\lambda+\frac{7}{2}}
t_{\lambda,\lambda+\frac{3}{2}}-
t_{\lambda+2,\lambda+\frac{7}{2}}t_{\lambda,\lambda+2}\right)+\\
\quad+\,\epsilon_{5,\lambda}\,t_{\lambda+4,\lambda+\frac{11}{2}}t_{\lambda+\frac{3}{2},\lambda+4}
t_{\lambda,\lambda+\frac{3}{2}}=0,\\
\epsilon_{6,\lambda}\,t_{\lambda,\lambda+\frac{3}{2}}\left(
t_{\lambda+4,\lambda+\frac{11}{2}}t_{\lambda+\frac{3}{2},\lambda+4}+\frac{1}{3}\,
t_{\lambda+\frac{7}{2},\lambda+\frac{11}{2}}t_{\lambda+\frac{3}{2},\lambda+\frac{7}{2}}\right)+\\
\quad
+t_{\lambda,\lambda+2}\left(t_{\lambda+\frac{7}{2},\lambda+\frac{11}{2}}
t_{\lambda+2,\lambda+\frac{7}{2}}-
t_{\lambda+4,\lambda+\frac{11}{2}}t_{\lambda+2,\lambda+4}\right)=0.
\end{gather*}
  \item [c)] For $2(\beta-\lambda)\in\left\{12,\,
\dots,\,2n\right\}:$\begin{gather*}
t_{\lambda+\frac{9}{2},\lambda+6}
\left(t_{\lambda+\frac{5}{2},\lambda+\frac{9}{2}}
t_{\lambda,\lambda+\frac{5}{2}}-
t_{\lambda+2,\lambda+\frac{9}{2}}t_{\lambda,\lambda+2}\right)=0,\\
t_{\lambda,\lambda+\frac{3}{2}} \left(t_{\lambda+4,\lambda+6}
t_{\lambda+\frac{3}{2},\lambda+4}-
t_{\lambda+\frac{7}{2},\lambda+6}t_{\lambda+\frac{3}{2},\lambda+\frac{7}{2}}\right)=0,\\
t_{\lambda+4,\lambda+6}\left(3\,t_{\lambda+\frac{3}{2},\lambda+4}
t_{\lambda,\lambda+\frac{3}{2}}+
3\,t_{\lambda+\frac{5}{2},\lambda+4}t_{\lambda,\lambda+\frac{5}{2}}+
t_{\lambda+2,\lambda+4}t_{\lambda,\lambda+2}\right)=0,\\
t_{\lambda,\lambda+2}\left((1-\epsilon_{7,\lambda})\,t_{\lambda+\frac{7}{2},\lambda+6}
t_{\lambda+2,\lambda+\frac{7}{2}}+
t_{\lambda+\frac{9}{2},\lambda+6}t_{\lambda+2,\lambda+\frac{9}{2}}+{1\over3}
t_{\lambda+4,\lambda+6}t_{\lambda+2,\lambda+4}\right)+\\
\quad +\,\epsilon_{7,\lambda}\,t_{\lambda+\frac{7}{2},\lambda+6}
t_{\lambda+\frac{3}{2},\lambda+\frac{7}{2}}
t_{\lambda,\lambda+\frac{3}{2}}=0,\\
\epsilon_{8,\lambda}\,t_{\lambda+\frac{7}{2},\lambda+6}
\left(t_{\lambda+\frac{3}{2},\lambda+\frac{7}{2}}
t_{\lambda,\lambda+\frac{3}{2}}-
t_{\lambda+2,\lambda+\frac{7}{2}}t_{\lambda,\lambda+2}\right)+\\
\quad +\left(t_{\lambda+4,\lambda+6}
t_{\lambda+\frac{5}{2},\lambda+4}-
t_{\lambda+\frac{9}{2},\lambda+6}t_{\lambda+\frac{5}{2},\lambda+\frac{9}{2}}\right)
t_{\lambda,\lambda+\frac{5}{2}}=0.\end{gather*}
  \item [d)] For $2(\beta-\lambda)\in\left\{13,\,
\dots,\,2n\right\}:$\begin{gather*}t_{\lambda+4,\lambda+\frac{13}{2}}\left(t_{\lambda+\frac{3}{2},\lambda+4}
t_{\lambda,\lambda+\frac{3}{2}}+
t_{\lambda+\frac{5}{2},\lambda+4}t_{\lambda,\lambda+\frac{5}{2}}+(\frac{1}{3}-\epsilon_{10,\lambda})
t_{\lambda+2,\lambda+4}t_{\lambda,\lambda+2}\right)+\\\quad
+\,\epsilon_{10,\lambda} t_{\lambda,\lambda+2}
t_{\lambda+\frac{9}{2},\lambda+\frac{13}{2}}
t_{\lambda+2,\lambda+\frac{9}{2}}=0,\\[0.5pt]
t_{\lambda,\lambda+\frac{5}{2}}\left(t_{\lambda+4,\lambda+\frac{13}{2}}
t_{\lambda+\frac{5}{2},\lambda+4}+\frac{1}{3}
\,t_{\lambda+\frac{9}{2},\lambda+\frac{13}{2}}t_{\lambda+\frac{5}{2},\lambda+\frac{9}{2}}\right)+
\\ \quad+\,\epsilon_{11,\lambda}t_{\lambda,\lambda+2}
\left(t_{\lambda+\frac{9}{2},\lambda+\frac{13}{2}}
t_{\lambda+2,\lambda+\frac{9}{2}}-
t_{\lambda+4,\lambda+\frac{13}{2}}t_{\lambda+2,\lambda+4}\right)=0,\\
t_{\lambda+\frac{9}{2},\lambda+\frac{13}{2}}
\left(t_{\lambda+\frac{5}{2},\lambda+\frac{9}{2}}
t_{\lambda,\lambda+\frac{5}{2}}+(\epsilon_{12,\lambda}-1)\,
t_{\lambda+2,\lambda+\frac{9}{2}}t_{\lambda,\lambda+2}\right)-\epsilon_{12,\lambda}
t_{\lambda,\lambda+2}t_{\lambda+4,\lambda+\frac{13}{2}}t_{\lambda+2,\lambda+4}=0.\end{gather*}
  \item [e)] For $2(\beta-\lambda)\in\left\{14,\,
\dots,\,2n\right\}:$\begin{gather*}t_{\lambda+\frac{9}{2},\lambda+7}
t_{\lambda+2,\lambda+\frac{9}{2}}t_{\lambda,\lambda+2}
=t_{\lambda,\lambda+\frac{5}{2}}
t_{\lambda+\frac{9}{2},\lambda+7}t_{\lambda+\frac{5}{2},\lambda+\frac{9}{2}}
=0.\end{gather*}
\end{itemize}
\end{proposition}
\begin{proof}
Considering again the homomorphism  condition, we compute the 3rd
order term ${\frak L}^{(3)}$ which is a solution of the
Maurer-Cartan equation:
\begin{equation}\label{cap3}\delta({\frak
L}^{(3)})=-{1\over2}\left([{\frak L}^{(1)},\, {\frak
L}^{(2)}]\!]+[\![{\frak L}^{(2)},\, {\frak
L}^{(1)}]\!]\right).\end{equation} The right hand side of
(\ref{cap3}) together with equation (\ref{ord2}) yield the
following maps:
\begin{equation*}\label{2-cocyc}\begin{array}{llllll}
\Omega_{\lambda,\lambda+\frac{9}{2}}&=& \varphi_1(t)\,
[\![\Upsilon_{\lambda+3,\lambda+\frac{9}{2}},\
\frak{J}_4^{-1,\lambda}
]\!]+\psi_1(t)\,[\![\frak{J}_4^{-1,\lambda+\frac{3}{2}},
\,\Upsilon_{\lambda,\lambda+\frac{3}{2}} ]\!]
&:&\mathcal{K}(1)\otimes\mathcal{K}(1)\rightarrow\frak{D}_{\lambda,\lambda+\frac{9}{2}},\\
\Omega_{\lambda,\lambda+5}&=&\varphi_2(t)\,[\![\Upsilon_{\lambda+\frac{7}{2},\lambda+5},\
\frak{J}_\frac{9}{2}^{-1,\lambda} ]\!]+\psi_2(t)\,
[\![\frak{J}_\frac{9}{2}^{-1,\lambda+\frac{3}{2}},
\,\Upsilon_{\lambda,\lambda+\frac{3}{2}}
]\!]&:&\mathcal{K}(1)\otimes\mathcal{K}(1)\rightarrow\frak{D}_{\lambda,\lambda+5},\\
\widetilde{\Omega}_{\lambda,\lambda+5}&=&\widetilde{\varphi}_2(t)\,[\![\Upsilon_{\lambda+3,\lambda+5},\
\frak{J}_4^{-1,\lambda} ]\!]+\widetilde{\psi}_2(t)\,
[\![\frak{J}_4^{-1,\lambda+2}, \,\Upsilon_{\lambda,\lambda+2}
]\!]&:&\mathcal{K}(1)\otimes\mathcal{K}(1)\rightarrow\frak{D}_{\lambda,\lambda+5},\\
\Omega_{\lambda,\lambda+{11\over2}}&=&\varphi_3(t)\,[\![\Upsilon_{\lambda+4,\lambda+\frac{11}{2}},\
\frak{J}_{5}^{-1,\lambda} ]\!]+\psi_3(t)\,
[\![\frak{J}_{5}^{-1,\lambda+\frac{3}{2}},
\,\Upsilon_{\lambda,\lambda+\frac{3}{2}}
]\!]&:&\mathcal{K}(1)\otimes\mathcal{K}(1)\rightarrow\frak{D}_{\lambda,\lambda+\frac{11}{2}},\\
\widetilde{\Omega}_{\lambda,\lambda+{11\over2}}&=&\widetilde{\varphi}_3(t)\,[\![\Upsilon_{\lambda+\frac{7}{2},\lambda+\frac{11}{2}},\
\frak{J}_\frac{9}{2}^{-1,\lambda} ]\!]+\widetilde{\psi}_3(t)\,
[\![\frak{J}_\frac{9}{2}^{-1,\lambda+2},
\,\Upsilon_{\lambda,\lambda+2}
]\!]&:&\mathcal{K}(1)\otimes\mathcal{K}(1)\rightarrow\frak{D}_{\lambda,\lambda+\frac{11}{2}},\\
\overline{\Omega}_{\lambda,\lambda+{11\over2}}&=&\overline{\varphi}_3(t)
\,[\![\Upsilon_{\lambda+3,\lambda+\frac{11}{2}},\
\frak{J}_4^{-1,\lambda} ]\!]+\overline{\psi}_3(t)\,
[\![\frak{J}_4^{-1,\lambda+\frac{5}{2}},
\,\Upsilon_{\lambda,\lambda+\frac{5}{2}}
]\!]&:&\mathcal{K}(1)\otimes\mathcal{K}(1)\rightarrow\frak{D}_{\lambda,\lambda+\frac{11}{2}},\\
\Omega_{\lambda,\lambda+6}&=&{\varphi}_4(t)\,[\![\Upsilon_{\lambda+\frac{9}{2},\lambda+6},\
\frak{J}_\frac{11}{2}^{-1,\lambda} ]\!]+{\psi}_4(t) \,[\![
\frak{J}_\frac{11}{2}^{-1,\lambda+\frac{3}{2}},
\,\Upsilon_{\lambda,\lambda+\frac{3}{2}}
]\!]&:&\mathcal{K}(1)\otimes\mathcal{K}(1)
\rightarrow\frak{D}_{\lambda,\lambda+6},\\
\widetilde{\Omega}_{\lambda,\lambda+6}&=&\widetilde{\varphi}_4(t)\,[\![\Upsilon_{\lambda+4,\lambda+6},\
\frak{J}_{5}^{-1,\lambda} ]\!]+\widetilde{\psi}_4(t) \,[\![
\frak{J}_{5}^{-1,\lambda+2}, \,\Upsilon_{\lambda,\lambda+2}
]\!]&:&\mathcal{K}(1)\otimes\mathcal{K}(1)
\rightarrow\frak{D}_{\lambda,\lambda+6},\\
\overline{\Omega}_{\lambda,\lambda+6}&=&\overline{\varphi}_4(t)\,[\![\Upsilon_{\lambda+\frac{7}{2},\lambda+6},\
\frak{J}_\frac{9}{2}^{-1,\lambda} ]\!]+\overline{\psi}_4(t)\,[\![
\frak{J}_\frac{9}{2}^{-1,\lambda+\frac{5}{2}},
\,\Upsilon_{\lambda,\lambda+\frac{5}{2}}
]\!]&:&\mathcal{K}(1)\otimes\mathcal{K}(1)
\rightarrow\frak{D}_{\lambda,\lambda+6},\\
\Omega_{\lambda,\lambda+{13\over2}}&=&{\varphi}_5(t)\,[\![\Upsilon_{\lambda+4,\lambda+\frac{13}{2}},\
\frak{J}_{5}^{-1,\lambda} ]\!]+{\psi}_5(t) \,[\![
\frak{J}_{5}^{-1,\lambda+\frac{5}{2}},
\,\Upsilon_{\lambda,\lambda+\frac{5}{2}}
]\!]&:&\mathcal{K}(1)\otimes\mathcal{K}(1)
\rightarrow\frak{D}_{\lambda,\lambda+\frac{13}{2}},\\
\widetilde{\Omega}_{\lambda,\lambda+{13\over2}}&=&\widetilde{\varphi}_5(t)\,[\![\Upsilon_{\lambda+{9\over2},\lambda+\frac{13}{2}},\
\frak{J}_{11\over2}^{-1,\lambda} ]\!]+\widetilde{\psi}_5(t)\, [\![
\frak{J}_{11\over2}^{-1,\lambda+2}, \,\Upsilon_{\lambda,\lambda+2}
]\!]&:&\mathcal{K}(1)\otimes\mathcal{K}(1)
\rightarrow\frak{D}_{\lambda,\lambda+\frac{13}{2}},\\
\Omega_{\lambda,\lambda+7}&=&\varphi_6(t)\,[\![\Upsilon_{\lambda+{9\over2},\lambda+7},\
\frak{J}_{11\over2}^{-1,\lambda} ]\!]+\psi_6(t)\, [\![
\frak{J}_{11\over2}^{-1,\lambda+{5\over2}},
\,\Upsilon_{\lambda,\lambda+{5\over2}}
]\!]&:&\mathcal{K}(1)\otimes\mathcal{K}(1)
\rightarrow\frak{D}_{\lambda,\lambda+7},\end{array}
\end{equation*}
where \begin{equation*}
\begin{array}{lll}
\varphi_1(t)=\zeta_\lambda^{-1}t_{\lambda+3,\lambda+\frac{9}{2}}
t_{\lambda+\frac{3}{2},\lambda+3}
t_{\lambda,\lambda+\frac{3}{2}},\\
\psi_1(t)=\zeta_{\lambda+\frac{3}{2}}^{-1}t_{\lambda+3,\lambda+\frac{9}{2}}
t_{\lambda+\frac{3}{2},\lambda+3}
t_{\lambda,\lambda+\frac{3}{2}},\\
\varphi_2(t)=\alpha_\lambda^{-1}t_{\lambda+\frac{7}{2},\lambda+5}
\left(t_{\lambda+\frac{3}{2},\lambda+\frac{7}{2}}
t_{\lambda,\lambda+\frac{3}{2}}-
t_{\lambda+2,\lambda+\frac{7}{2}}t_{\lambda,\lambda+2}\right),\\
\psi_2(t)=\alpha_{\lambda+\frac{3}{2}}^{-1}
\left(t_{\lambda+3,\lambda+5} t_{\lambda+\frac{3}{2},\lambda+3}-
t_{\lambda+\frac{7}{2},\lambda+5}t_{\lambda+\frac{3}{2},\lambda+\frac{7}{2}}\right)
t_{\lambda,\lambda+\frac{3}{2}},\\
\widetilde{\varphi}_2(t)=\zeta_\lambda^{-1}
t_{\lambda+3,\lambda+5}
t_{\lambda+\frac{3}{2},\lambda+3}t_{\lambda,\lambda+\frac{3}{2}},\\
\widetilde{\psi}_2(t)=\zeta_{\lambda+2}^{-1}
t_{\lambda+\frac{7}{2},\lambda+5}t_{\lambda+2,\lambda+\frac{7}{2}}
t_{\lambda,\lambda+2},\\
\varphi_3(t)=\beta_\lambda^{-1}t_{\lambda+4,\lambda+\frac{11}{2}}\left(t_{\lambda+\frac{3}{2},\lambda+4}
t_{\lambda,\lambda+\frac{3}{2}}+
t_{\lambda+\frac{5}{2},\lambda+4}t_{\lambda,\lambda+\frac{5}{2}}+
{1\over3}t_{\lambda+2,\lambda+4}t_{\lambda,\lambda+2}\right),\\
\psi_3(t)=\beta_{\lambda+\frac{3}{2}}^{-1}t_{\lambda,\lambda+\frac{3}{2}}\left(t_{\lambda+3,\lambda+\frac{11}{2}}
t_{\lambda+\frac{3}{2},\lambda+3}+
t_{\lambda+4,\lambda+\frac{11}{2}}t_{\lambda+\frac{3}{2},\lambda+4}+
{1\over3}t_{\lambda+\frac{7}{2},\lambda+\frac{11}{2}}
t_{\lambda+\frac{3}{2},\lambda+\frac{7}{2}}\right),\\
\widetilde{\varphi}_3(t)=\alpha_\lambda^{-1}t_{\lambda+\frac{7}{2},\lambda+\frac{11}{2}}
\left(t_{\lambda+\frac{3}{2},\lambda+\frac{7}{2}}
t_{\lambda,\lambda+\frac{3}{2}}-
t_{\lambda+2,\lambda+\frac{7}{2}}t_{\lambda,\lambda+2}\right),\\
\widetilde{\psi}_3(t)=\alpha_{\lambda+2}^{-1}t_{\lambda,\lambda+2}
\left(t_{\lambda+\frac{7}{2},\lambda+\frac{11}{2}}
t_{\lambda+2,\lambda+\frac{7}{2}}-
t_{\lambda+4,\lambda+\frac{11}{2}}t_{\lambda+2,\lambda+4}\right),\\
\overline{\varphi}_3(t)=\zeta_\lambda^{-1}t_{\lambda+3,\lambda+\frac{11}{2}}
t_{\lambda+\frac{3}{2},\lambda+3}
t_{\lambda,\lambda+\frac{3}{2}},\\
\overline{\psi}_3(t)=\zeta_{\lambda+\frac{5}{2}}^{-1}t_{\lambda+4,\lambda+\frac{11}{2}}
t_{\lambda+\frac{5}{2},\lambda+4}
t_{\lambda,\lambda+\frac{5}{2}},\\
\varphi_4(t)=\gamma_\lambda^{-1}t_{\lambda+\frac{9}{2},\lambda+6}
\left(t_{\lambda+\frac{5}{2},\lambda+\frac{9}{2}}
t_{\lambda,\lambda+\frac{5}{2}}-
t_{\lambda+2,\lambda+\frac{9}{2}}t_{\lambda,\lambda+2}\right),\\\end{array}
\end{equation*}\begin{equation*}
\begin{array}{llllllllllllllll}
\psi_4(t)=\gamma_{\lambda+\frac{3}{2}}^{-1}
\left(t_{\lambda+4,\lambda+6} t_{\lambda+\frac{3}{2},\lambda+4}-
t_{\lambda+\frac{7}{2},\lambda+6}t_{\lambda+\frac{3}{2},\lambda+\frac{7}{2}}\right)
t_{\lambda,\lambda+\frac{3}{2}},\\
\widetilde{\varphi}_4(t)=\beta_\lambda^{-1}t_{\lambda+4,\lambda+6}\left(t_{\lambda+\frac{3}{2},\lambda+4}
t_{\lambda,\lambda+\frac{3}{2}}+
t_{\lambda+\frac{5}{2},\lambda+4}t_{\lambda,\lambda+\frac{5}{2}}+
{1\over3}t_{\lambda+2,\lambda+4}t_{\lambda,\lambda+2}\right),\\
\widetilde{\psi}_4(t)=\beta_{\lambda+2}^{-1}t_{\lambda,\lambda+2}\left(t_{\lambda+\frac{7}{2},\lambda+6}
t_{\lambda+2,\lambda+\frac{7}{2}}+
t_{\lambda+\frac{9}{2},\lambda+6}t_{\lambda+2,\lambda+\frac{9}{2}}+
{1\over3}t_{\lambda+4,\lambda+6}t_{\lambda+2,\lambda+4}\right),\\
\overline{\varphi}_4(t)=\alpha_\lambda^{-1}t_{\lambda+\frac{7}{2},\lambda+6}
\left(t_{\lambda+\frac{3}{2},\lambda+\frac{7}{2}}
t_{\lambda,\lambda+\frac{3}{2}}-
t_{\lambda+2,\lambda+\frac{7}{2}}t_{\lambda,\lambda+2}\right),\\
\overline{\psi}_4(t)=\alpha_{\lambda+\frac{5}{2}}^{-1}
\left(t_{\lambda+4,\lambda+6} t_{\lambda+\frac{5}{2},\lambda+4}-
t_{\lambda+\frac{9}{2},\lambda+6}t_{\lambda+\frac{5}{2},\lambda+\frac{9}{2}}\right)
t_{\lambda,\lambda+\frac{5}{2}},\\
{\varphi}_5(t)=\beta_\lambda^{-1}t_{\lambda+4,\lambda+\frac{13}{2}}\left(t_{\lambda+\frac{3}{2},\lambda+4}
t_{\lambda,\lambda+\frac{3}{2}}+
t_{\lambda+\frac{5}{2},\lambda+4}t_{\lambda,\lambda+\frac{5}{2}}+
{1\over3}t_{\lambda+2,\lambda+4}t_{\lambda,\lambda+2}\right),
\\
{\psi}_5(t)=\beta_{\lambda+\frac{5}{2}}^{-1}t_{\lambda,\lambda+\frac{5}{2}}
\left(t_{\lambda+4,\lambda+\frac{13}{2}}
t_{\lambda+\frac{5}{2},\lambda+4}+
{1\over3}t_{\lambda+\frac{9}{2},\lambda+\frac{13}{2}}t_{\lambda+\frac{5}{2},\lambda+\frac{9}{2}}\right),\\
\widetilde{\varphi}_5(t)=\gamma_\lambda^{-1}t_{\lambda+\frac{9}{2},\lambda+\frac{13}{2}}
\left(t_{\lambda+\frac{5}{2},\lambda+\frac{9}{2}}
t_{\lambda,\lambda+\frac{5}{2}}-
t_{\lambda+2,\lambda+\frac{9}{2}}t_{\lambda,\lambda+2}\right),
\\\widetilde{\psi}_5(t)={\gamma}_{\lambda+2}^{-1}t_{\lambda,\lambda+2}
\left(t_{\lambda+\frac{9}{2},\lambda+\frac{13}{2}}
t_{\lambda+2,\lambda+\frac{9}{2}}-
t_{\lambda+4,\lambda+\frac{13}{2}}t_{\lambda+2,\lambda+4}\right),\\
\varphi_6(t)=\gamma_\lambda^{-1}t_{\lambda+\frac{9}{2},\lambda+7}
\left(t_{\lambda+\frac{5}{2},\lambda+\frac{9}{2}}
t_{\lambda,\lambda+\frac{5}{2}}-
t_{\lambda+2,\lambda+\frac{9}{2}}t_{\lambda,\lambda+2}\right),
\\\psi_6(t)=\gamma_{\lambda+\frac{5}{2}}^{-1}t_{\lambda,\lambda+\frac{5}{2}}
t_{\lambda+\frac{9}{2},\lambda+7}t_{\lambda+\frac{5}{2},\lambda+\frac{9}{2}}.
\end{array}
\end{equation*}
 Now, the same arguments,
as in the proof of Proposition \ref{th2}, show that we must have:
\begin{equation*}
\label{E}
\begin{array}{llllll}
\Omega_{\lambda,\lambda+{9\over2}}&=&
\omega_1(t)\delta\left(\frak{J}_{{11\over2}}^{-1,\lambda}\right),\\
\Omega_{\lambda,\lambda+5}+\widetilde{\Omega}_{\lambda,\lambda+5}&=&
\omega_2(t)\delta\left(\frak{J}_{6}^{-1,\lambda}\right),\\
\Omega_{\lambda,\lambda+{11\over2}}+\widetilde{\Omega}_{\lambda,\lambda+{11\over2}}+
\overline{\Omega}_{\lambda,\lambda+{11\over2}}&=&
\omega_3(t)\delta\left(\frak{J}_{13\over2}^{-1,\lambda}\right),\\
\Omega_{\lambda,\lambda+6}+\widetilde{\Omega}_{\lambda,\lambda+6}+\overline{\Omega}_{\lambda,\lambda+6}&=&
\omega_4(t)\delta\left(\frak{J}_{7}^{-1,\lambda}\right),\\
\Omega_{\lambda,\lambda+{13\over2}}+\widetilde{\Omega}_{\lambda,\lambda+{13\over2}}&=&
\omega_5(t)\delta\left(\frak{J}_{15\over2}^{-1,\lambda}\right),\\
\Omega_{\lambda,\lambda+7}&=&
\omega_6(t)\delta\left(\frak{J}_{8}^{-1,\lambda}\right),\end{array}\end{equation*}
where $\omega_1,\,\dots,\,\omega_5$ are some functions. So, by
Lemma \ref{benfraj1} and Lemma \ref{benfraj2}, we obtain for the
nonzero $\varphi_i(t), \widetilde{\varphi}_i(t),
\overline{\varphi}_i, \psi_i(t), \widetilde{\psi}_i(t),
\overline{\psi}_i(t)$ and $\omega_i(t)$ the following relation:
\begin{equation*}
\begin{array}{llllll}
\alpha_{\lambda}\,\varphi_2(t)+\epsilon_{1,\lambda}\,
\alpha_{\lambda+\frac{3}{2}}\,\psi_2(t)=0,\quad
&\zeta_\lambda\,\widetilde{\varphi}_2(t)+\epsilon_{2,\lambda}\,
\alpha_{\lambda+\frac{3}{2}}\,\psi_2(t)=0,\\[1pt]
\beta_{\lambda}\,\varphi_3(t)+\epsilon_{4,\lambda}\,
\beta_{\lambda+\frac{3}{2}}\,\psi_3(t)=0,\quad
&\alpha_{\lambda}\,\widetilde{\varphi}_3(t)+\epsilon_{5,\lambda}\,
\beta_{\lambda+\frac{3}{2}}\,\psi_3(t)=0,\\[1pt]
\alpha_{\lambda+2}\,\widetilde{\psi}_3(t)+\epsilon_{6,\lambda}\,
\beta_{\lambda+\frac{3}{2}}\,\psi_3(t)=0,\quad
&\beta_{\lambda+2}\,\widetilde{\psi}_4(t)+\epsilon_{7,\lambda}\,
\alpha_{\lambda}\,\overline{\varphi}_4(t)=0,\\[1pt]
\alpha_{\lambda+\frac{5}{2}}\,\overline{\psi}_4(t)+\epsilon_{8,\lambda}\,
\alpha_{\lambda}\,\overline{\varphi}_4(t)=0,\quad &
\beta_{\lambda}\,\varphi_5(t)+\epsilon_{10,\lambda}\,\gamma_{\lambda+2}\widetilde{\psi}_5(t)=0,\\[1pt]
\beta_{\lambda+\frac{5}{2}}\,\psi_5(t)+\epsilon_{11,\lambda}\,\gamma_{\lambda+2}\widetilde{\psi}_5(t)=0,\quad
&\gamma_{\lambda}\,\widetilde{\varphi}_5(t)+\epsilon_{12,\lambda}
\,\gamma_{\lambda+2}\widetilde{\psi}_5(t)=0,\\[1pt]
\omega_4(t)=\epsilon_{9,\lambda}\,\alpha_{\lambda}\,\overline{\varphi}_4(t),\quad
&\omega_2(t)=\epsilon_{3,\lambda}\,
\alpha_{\lambda+\frac{3}{2}}\,\psi_2(t),\\[1pt]
\omega_1(t)=\xi_\lambda^{-1}\,
\zeta_\lambda\,\varphi_1(t)=\xi_\lambda^{-1}\,\zeta_{\lambda+\frac{3}{2}}\,\psi_1(t).
\end{array}
\end{equation*}
Therefore, we get the necessary integrability conditions for
${\frak L}^{(3)}$. Under these conditions, the third-order term
${\frak L}^{(3)}$ can be given by:
\begin{equation*}
\begin{array}{llll}
{\frak L}^{(3) }=&
\sum_\lambda\xi_\lambda^{-1}t_{\lambda+3,\lambda+\frac{9}{2}}
t_{\lambda+\frac{3}{2},\lambda+3}
t_{\lambda,\lambda+\frac{3}{2}}\,
\frak{J}_{11\over2}^{-1,\lambda}+\\
&\sum_\lambda\epsilon_{3,\lambda} \left(t_{\lambda+3,\lambda+5}
t_{\lambda+\frac{3}{2},\lambda+3}-
t_{\lambda+\frac{7}{2},\lambda+5}t_{\lambda+\frac{3}{2},\lambda+\frac{7}{2}}\right)
t_{\lambda,\lambda+\frac{3}{2}}\,
\frak{J}_{6}^{-1,\lambda}+\\
&\sum_\lambda\epsilon_{9,\lambda}\,t_{\lambda+\frac{7}{2},\lambda+6}
\left(t_{\lambda+\frac{3}{2},\lambda+\frac{7}{2}}
t_{\lambda,\lambda+\frac{3}{2}}-
t_{\lambda+2,\lambda+\frac{7}{2}}t_{\lambda,\lambda+2}\right)\frak{J}_{7}^{-1,\lambda}.
\end{array}
\end{equation*}
\end{proof}

\begin{proposition}\label{th4} The 4th order integrability conditions of the
infinitesimal deformation~(\ref{infdef1}) are the following:
\begin{itemize}
 \item [a)] For $2(\beta-\lambda)\in\left\{12,\,
\dots,\,2n\right\}:$
\begin{gather*}
t_{\lambda+\frac{9}{2},\lambda+6}
t_{\lambda+3,\lambda+\frac{9}{2}}t_{\lambda+\frac{3}{2},\lambda+3}
t_{\lambda,\lambda+\frac{3}{2}} =0.
\end{gather*}
 \item [b)] For $2(\beta-\lambda)\in\left\{13,\,
\dots,\,2n\right\}:$\begin{gather*}
t_{\lambda+5,\lambda+\frac{13}{2}}
t_{\lambda+\frac{7}{2},\lambda+5}t_{\lambda+\frac{3}{2},\lambda+\frac{7}{2}}
t_{\lambda,\lambda+\frac{3}{2}}=0,\\
t_{\lambda+\frac{3}{2},\lambda+3} t_{\lambda,\lambda+\frac{3}{2}}
t_{\lambda+5,\lambda+\frac{13}{2}}t_{\lambda+3,\lambda+5}=0,\\
t_{\lambda+\frac{9}{2},\lambda+\frac{13}{2}}
t_{\lambda+3,\lambda+\frac{9}{2}}
t_{\lambda+\frac{3}{2},\lambda+3}
t_{\lambda,\lambda+\frac{3}{2}}=0.
\end{gather*}
\item [c)] For $2(\beta-\lambda)\in\left\{14,\,
\dots,\,2n\right\}:$
\begin{gather*}
t_{\lambda+\frac{9}{2},\lambda+7}
t_{\lambda+3,\lambda+\frac{9}{2}}
t_{\lambda+\frac{3}{2},\lambda+3}
t_{\lambda,\lambda+\frac{3}{2}}=0,\\
t_{\lambda,\lambda+\frac{3}{2}}t_{\lambda+\frac{3}{2},\lambda+3}
t_{\lambda+5,\lambda+7}t_{\lambda+3,\lambda+5}
=0,\\
t_{\lambda,\lambda+2}t_{\lambda+2,\lambda+\frac{7}{2}}
t_{\lambda+\frac{11}{2},\lambda+7}t_{\lambda+\frac{7}{2},\lambda+\frac{11}{2}}
=0,\\
t_{\lambda,\lambda+\frac{3}{2}}t_{\lambda+\frac{3}{2},\lambda+\frac{7}{2}}
t_{\lambda+\frac{11}{2},\lambda+7}t_{\lambda+\frac{7}{2},\lambda+\frac{11}{2}}
=0,\\
t_{\lambda,\lambda+\frac{3}{2}}t_{\lambda+5,\lambda+7}
t_{\lambda+\frac{7}{2},\lambda+5}t_{\lambda+\frac{3}{2},\lambda+\frac{7}{2}}
=0,\\
t_{\lambda+\frac{11}{2},\lambda+7}
t_{\lambda+4,\lambda+\frac{11}{2}}\left(3t_{\lambda+\frac{3}{2},\lambda+4}
t_{\lambda,\lambda+\frac{3}{2}}+
t_{\lambda+2,\lambda+4}t_{\lambda,\lambda+2}\right)=0.
\end{gather*}
\item [d)] For $2(\beta-\lambda)\in\left\{15,\,
\dots,\,2n\right\}:$\begin{gather*}
t_{\lambda,\lambda+\frac{3}{2}}
t_{\lambda+\frac{3}{2},\lambda+3}t_{\lambda+5,\lambda+\frac{15}{2}}
t_{\lambda+3,\lambda+5}=0,\\
t_{\lambda,\lambda+\frac{3}{2}}
t_{\lambda+5,\lambda+\frac{15}{2}}t_{\lambda+\frac{7}{2},\lambda+5}
t_{\lambda+\frac{3}{2},\lambda+\frac{7}{2}}=0,\\
t_{\lambda,\lambda+\frac{5}{2}} t_{\lambda+\frac{5}{2},\lambda+4}
t_{\lambda+6,\lambda+\frac{15}{2}}t_{\lambda+4,\lambda+6}=0,\\
t_{\lambda+\frac{11}{2},\lambda+\frac{15}{2}}
t_{\lambda+\frac{7}{2},\lambda+\frac{11}{2}}\left(t_{\lambda+\frac{3}{2},\lambda+\frac{7}{2}}
t_{\lambda,\lambda+\frac{3}{2}}-
t_{\lambda+2,\lambda+\frac{7}{2}}t_{\lambda,\lambda+2}\right)=0,\\
t_{\lambda+\frac{11}{2},\lambda+\frac{15}{2}}
t_{\lambda+4,\lambda+\frac{11}{2}}\left(3t_{\lambda+\frac{3}{2},\lambda+4}
t_{\lambda,\lambda+\frac{3}{2}}+
t_{\lambda+2,\lambda+4}t_{\lambda,\lambda+2}\right)=0,\\
t_{\lambda+6,\lambda+\frac{15}{2}}
t_{\lambda+\frac{7}{2},\lambda+6}\left(t_{\lambda+\frac{3}{2},\lambda+\frac{7}{2}}
t_{\lambda,\lambda+\frac{3}{2}}-
t_{\lambda+2,\lambda+\frac{7}{2}}t_{\lambda,\lambda+2}\right)=0.
\end{gather*}
\item [e)] For $2(\beta-\lambda)\in\left\{16,\,
\dots,\,2n\right\}:$\begin{equation*} \begin{array}{lllllll}
\epsilon_{14,\lambda}\,t_{\lambda+6,\lambda+8}t_{\lambda+\frac{7}{2},\lambda+6}
\left(t_{\lambda+\frac{3}{2},\lambda+\frac{7}{2}}
t_{\lambda,\lambda+\frac{3}{2}}-
t_{\lambda+2,\lambda+\frac{7}{2}}t_{\lambda,\lambda+2}\right)+\\
\quad~+
t_{\lambda+\frac{13}{2},\lambda+8}t_{\lambda+\frac{9}{2},\lambda+\frac{13}{2}}
\left( t_{\lambda+\frac{5}{2},\lambda+\frac{9}{2}}
t_{\lambda,\lambda+\frac{5}{2}}-
t_{\lambda+2,\lambda+\frac{9}{2}}t_{\lambda,\lambda+2}\right)=0
,\\
\epsilon_{16,\lambda}\,t_{\lambda+6,\lambda+8}t_{\lambda+\frac{7}{2},\lambda+6}
\left(t_{\lambda+\frac{3}{2},\lambda+\frac{7}{2}}
t_{\lambda,\lambda+\frac{3}{2}}-
t_{\lambda+2,\lambda+\frac{7}{2}}t_{\lambda,\lambda+2}\right)+\\
\quad~+ t_{\lambda,\lambda+2}t_{\lambda+\frac{11}{2},\lambda+8}
\left(t_{\lambda+\frac{7}{2},\lambda+\frac{11}{2}}
t_{\lambda+2,\lambda+\frac{7}{2}}-
t_{\lambda+4,\lambda+\frac{11}{2}}t_{\lambda+2,\lambda+4}\right)=0,\\
\left((1+\epsilon_{15,\lambda})\,t_{\lambda+6,\lambda+8}t_{\lambda+\frac{7}{2},\lambda+6}-
t_{\lambda+\frac{11}{2},\lambda+8}t_{\lambda+\frac{7}{2},\lambda+\frac{11}{2}}
\right)\times\\
\quad~\times \left(t_{\lambda+\frac{3}{2},\lambda+\frac{7}{2}}
t_{\lambda,\lambda+\frac{3}{2}}-
t_{\lambda+2,\lambda+\frac{7}{2}}t_{\lambda,\lambda+2}\right)=0,\\\end{array}\end{equation*}
\begin{equation*}
\begin{array}{llllll}
\epsilon_{13,\lambda}\,t_{\lambda+6,\lambda+8}t_{\lambda+\frac{7}{2},\lambda+6}
\left(t_{\lambda+\frac{3}{2},\lambda+\frac{7}{2}}
t_{\lambda,\lambda+\frac{3}{2}}-
t_{\lambda+2,\lambda+\frac{7}{2}}t_{\lambda,\lambda+2}\right)+\\
\quad~+t_{\lambda+\frac{11}{2},\lambda+8}
t_{\lambda+4,\lambda+\frac{11}{2}}
\left(t_{\lambda+\frac{3}{2},\lambda+4}
t_{\lambda,\lambda+\frac{3}{2}}+{1\over3}
t_{\lambda+2,\lambda+4}t_{\lambda,\lambda+2}\right)+\\
\quad~+ t_{\lambda+\frac{13}{2},\lambda+8}
t_{\lambda+4,\lambda+\frac{13}{2}}
\left(t_{\lambda+\frac{3}{2},\lambda+4}
t_{\lambda,\lambda+\frac{3}{2}}+
t_{\lambda+\frac{5}{2},\lambda+4}t_{\lambda,\lambda+\frac{5}{2}}+
\frac{1}{3}t_{\lambda+2,\lambda+4}t_{\lambda,\lambda+2}\right)=0.\end{array}
\end{equation*}
\item [f)] For $2(\beta-\lambda)\in\left\{17,\,
\dots,\,2n\right\}:$\begin{gather*}
t_{\lambda,\lambda+\frac{5}{2}} t_{\lambda+6,\lambda+\frac{17}{2}}
t_{\lambda+4,\lambda+6} t_{\lambda+\frac{5}{2},\lambda+4}=0,\\
t_{\lambda+\frac{13}{2},\lambda+\frac{17}{2}}
t_{\lambda+\frac{9}{2},\lambda+\frac{13}{2}}
\left(t_{\lambda+\frac{5}{2},\lambda+\frac{9}{2}}
t_{\lambda,\lambda+\frac{5}{2}}-
t_{\lambda+2,\lambda+\frac{9}{2}}t_{\lambda,\lambda+2}\right)=0,\\
t_{\lambda+6,\lambda+\frac{17}{2}}
t_{\lambda+\frac{7}{2},\lambda+6}
\left(t_{\lambda+\frac{3}{2},\lambda+\frac{7}{2}}
t_{\lambda,\lambda+\frac{3}{2}}-
t_{\lambda+2,\lambda+\frac{7}{2}}t_{\lambda,\lambda+2}\right)=0,\\
t_{\lambda+\frac{13}{2},\lambda+\frac{17}{2}}
t_{\lambda+4,\lambda+\frac{13}{2}}
\left(3\,t_{\lambda+\frac{3}{2},\lambda+4}
t_{\lambda,\lambda+\frac{3}{2}}+
3\,t_{\lambda+\frac{5}{2},\lambda+4}t_{\lambda,\lambda+\frac{5}{2}}+
t_{\lambda+2,\lambda+4}t_{\lambda,\lambda+2}\right)0.
\end{gather*}
\item [g)] For $2(\beta-\lambda)\in\left\{18,\,
\dots,\,2n\right\}:$\begin{gather*}
t_{\lambda+\frac{13}{2},\lambda+9}t_{\lambda+\frac{9}{2},\lambda+\frac{13}{2}}
\left(t_{\lambda+\frac{5}{2},\lambda+\frac{9}{2}}
t_{\lambda,\lambda+\frac{5}{2}}-
t_{\lambda+2,\lambda+\frac{9}{2}}t_{\lambda,\lambda+2}\right)=0.
\end{gather*}
\end{itemize}
\end{proposition}

To prove Proposition \ref{th4}, we need the following two lemmas
which we can check by a direct computation with the help of {\it
Maple}.
\begin{lemma}
\label{benfraj3} We have
\begin{equation*}
\begin{array}{llllll}
\epsilon_{9,\lambda}\,[\![\Upsilon_{\lambda+6,\lambda+8},\
\frak{J}_{7}^{-1,\lambda}
]\!]&=&\epsilon_{13,\lambda}\,\beta_{\lambda+4}^{-1}\,\beta_{\lambda}^{-1}\,
[\![\frak{J}_5^{-1,\lambda+4},\ \frak{J}_5^{-1,\lambda}
]\!]+\epsilon_{14,\lambda}\,\alpha_{\lambda+\frac{9}{2}}^{-1}\,\gamma_{\lambda}^{-1}
[\![\frak{J}_{\frac{9}{2}}^{-1,\lambda+\frac{9}{2}},\
\frak{J}_{\frac{11}{2}}^{-1,\lambda}
]\!]+\\&&\epsilon_{15,\lambda}\,\alpha_{\lambda}^{-1}\,\gamma_{\lambda+\frac{7}{2}}^{-1}
[\![\frak{J}_{\frac{11}{2}}^{-1,\lambda+\frac{7}{2}},\
\frak{J}_{\frac{9}{2}}^{-1,\lambda} ]\!]+
\epsilon_{16,\lambda}\,\epsilon_{9,\lambda+2}\,[\![\frak{J}_{7}^{-1,\lambda+2},\
 \Upsilon_{\lambda,\lambda+2}]\!]+\\&&\epsilon_{17,\lambda}\,
\delta(\frak{J}_{9}^{-1,\lambda}),
\end{array}
\end{equation*}
where
\begin{equation*}\small{
\begin{array}{llllll}
\epsilon_{13,\lambda}&=&\frac{\epsilon_{9,\lambda}\,(2\lambda+3)
(2\lambda+5)(2\lambda+9)(\lambda+2)(\lambda+3)(\lambda+5)(2\lambda^2+7\lambda+2)
(2\lambda^2+23\lambda+62)(16\lambda^4+240\lambda^3+1034\lambda^2+1005\lambda+300)}{36(\lambda+4)(2\lambda+7)
(32\lambda^6+784\lambda^5+7156\lambda^4+29576\lambda^3+53961\lambda^2+40281\lambda+11760)},\\[4pt]
\epsilon_{14,\lambda}&=&
\frac{\epsilon_{9,\lambda}\,(2\lambda+3)(2\lambda+5)(2\lambda-5)(2\lambda+9)(\lambda-4)(\lambda+2)(\lambda+7)(\lambda+9)}
{(\lambda+4)(32\lambda^6+784\lambda^5+7156\lambda^4+29576\lambda^3+53961\lambda^2+40281\lambda+11760)},\\[4pt]
\epsilon_{15,\lambda}&=&
-\frac{\epsilon_{9,\lambda}\,(2\lambda-3)(2\lambda+1)(2\lambda+3)(2\lambda+6)
(2\lambda+23)(\lambda+2)(\lambda+5)(\lambda+10)}
{(2\lambda+7)(32\lambda^6+784\lambda^5+7156\lambda^4+29576\lambda^3+53961\lambda^2+40281\lambda+11760)},\\[4pt]
\epsilon_{16,\lambda}&=&-\frac{\epsilon_{9,\lambda}\,(2\lambda+3)(\lambda+2)
(32\lambda^6+656\lambda^5+4756\lambda^4+14104\lambda^3+14901\lambda^2+7059\lambda+240)}
{\epsilon_{9,\lambda+2}(2\lambda+11)(\lambda+6)(32\lambda^6+784\lambda^5+
7156\lambda^4+29576\lambda^3+53961\lambda^2+40281\lambda+11760)},\\[4pt]
\epsilon_{17,\lambda}&=&-\frac{\epsilon_{9,\lambda}\,(2\lambda+5)(2\lambda+9)(2\lambda+6)
(\lambda+5)(16\lambda^4+240\lambda^3+1034\lambda^2+1005\lambda+420)}
{(2\lambda+11)(2\lambda+7)(\lambda+4)(\lambda+6)(32\lambda^6+784\lambda^5+7156\lambda^4+29576\lambda^3+
53961\lambda^2+40281\lambda+11760)}
\end{array}}
\end{equation*}
\end{lemma}
\begin{lemma}
\label{benfraj4} Each of the following systems is linearly
independent
\begin{equation*}
\begin{array}{lll}
1)\,&\left(\delta(\frak{J}_{7}^{-1,\lambda}),\,\,
\zeta_{\lambda+3}^{-1}\zeta_{\lambda}^{-1}[\![\frak{J}_{4}^{-1,\lambda+3},\,
\frak{J}_{4}^{-1,\lambda}]\!]+\xi_{\lambda}^{-1}[\![\Upsilon_{\lambda+\frac{9}{2},\lambda+6}
,\,\frak{J}_{\frac{11}{2}}^{-1,\lambda}]\!]+\xi_{\lambda+\frac{3}{2}}^{-1}
[\![\frak{J}_{\frac{11}{2}}^{-1,\lambda+\frac{3}{2}},\,
\Upsilon_{\lambda,\lambda+\frac{3}{2}}]\!]\right),\\[2pt]
2)\,&\Big(\delta(\frak{J}_{\frac{15}{2}}^{-1,\lambda}),\,\,
[\![\frak{J}_{4}^{-1,\lambda+\frac{7}{2}},\,
\frak{J}_{\frac{9}{2}}^{-1,\lambda}]\!],\,\,
[\![\Upsilon_{\lambda+5,\lambda+\frac{13}{2}},\,
\frak{J}_{6}^{-1,\lambda}]\!],\,\,[\![\Upsilon_{\lambda+\frac{9}{2},\lambda+\frac{13}{2}},\,
\frak{J}_{\frac{11}{2}}^{-1,\lambda}]\!],\,\,\\&~~~
\zeta^{-1}_\lambda\alpha_{\lambda+3}^{-1}
[\![\frak{J}_{\frac{9}{2}}^{-1,\lambda+3},\
\frak{J}_{4}^{-1,\lambda}
]\!]+\epsilon_{3,\lambda+\frac{3}{2}}[\![\frak{J}_{6}^{-1,\lambda+\frac{3}{2}},\
\Upsilon_{\lambda,\lambda+\frac{3}{2}} ]\!]\Big),\\[2pt]
3)\,&\Big(\delta(\frak{J}_{8}^{-1,\lambda}),\,\,
[\![\frak{J}_{\frac{9}{2}}^{-1,\lambda+\frac{7}{2}},\,
\frak{J}_{\frac{9}{2}}^{-1,\lambda}]\!],\,\,
[\![\frak{J}_{4}^{-1,\lambda+4},\,
\frak{J}_{5}^{-1,\lambda}]\!],\,\,[\![\frak{J}_{5}^{-1,\lambda+3},\,
\frak{J}_{4}^{-1,\lambda}]\!],\,\,[\![\frak{J}_{6}^{-1,\lambda+2},\,
\Upsilon_{\lambda,\lambda+2}]\!],\\&~~~[\![\Upsilon_{\lambda+\frac{9}{2},\lambda+7},\,
\frak{J}_{\frac{11}{2}}^{-1,\lambda}]\!],\,\,[\![\Upsilon_{\lambda+5,\lambda+7},\,
\frak{J}_{6}^{-1,\lambda}]\!]
 \Big),\\[2pt]
4)\,&\Big(\delta(\frak{J}_{\frac{17}{2}}^{-1,\lambda}),\,\,
[\![\frak{J}_{5}^{-1,\lambda+\frac{7}{2}},\,
\frak{J}_{\frac{9}{2}}^{-1,\lambda}]\!],\,\,
[\![\frak{J}_{\frac{9}{2}}^{-1,\lambda+4},\,
\frak{J}_{5}^{-1,\lambda}]\!],\,\,
[\![\frak{J}_{\frac{11}{2}}^{-1,\lambda+3},\,
\frak{J}_{4}^{-1,\lambda}]\!],\,\,
[\![\Upsilon_{\lambda+6,\lambda+\frac{15}{2}},\,
\frak{J}_{7}^{-1,\lambda}]\!],\\&~~~
[\![\frak{J}_{6}^{-1,\lambda+\frac{5}{2}},\,
\Upsilon_{\lambda,\lambda+\frac{5}{2}}]\!],\,\,
\epsilon_{3,\lambda} [\![\Upsilon_{\lambda+5,\lambda+\frac{15}{2}}
,\,\frak{J}_{6}^{-1,\lambda}]\!]+\epsilon_{9,\lambda+\frac{3}{2}}
[\![\frak{J}_{7}^{-1,\lambda+\frac{3}{2}},\
\Upsilon_{\lambda,\lambda+\frac{3}{2}} ]\!] \Big),
\\[2pt]
5)\,&\left(\delta(\frak{J}_{9}^{-1,\lambda}),\,\,
[\![\frak{J}_{5}^{-1,\lambda+4},\,
\frak{J}_{5}^{-1,\lambda}]\!],\,\,[\![\frak{J}_{\frac{9}{2}}^{-1,\lambda+\frac{9}{2}},\,
\frak{J}_{\frac{11}{2}}^{-1,\lambda}]\!],\,\,
[\![\frak{J}_{\frac{11}{2}}^{-1,\lambda+\frac{7}{2}},\,
\frak{J}_{\frac{9}{2}}^{-1,\lambda}]\!],\,\,[\![\frak{J}_{7}^{-1,\lambda+2},\,
\Upsilon_{\lambda,\lambda+2}]\!]\right),
\\[2pt]
6)\,&\left(\delta(\frak{J}_{\frac{19}{2}}^{-1,\lambda}),\,\,[\![\frak{J}_{5}^{-1,\lambda+\frac{9}{2}},\,
\frak{J}_{\frac{11}{2}}^{-1,\lambda}]\!],\,\,
[\![\frak{J}_{\frac{11}{2}}^{-1,\lambda+4},\,
\frak{J}_{5}^{-1,\lambda}]\!],\,\,[\![\frak{J}_{7}^{-1,\lambda+\frac{5}{2}},\,
\Upsilon_{\lambda,\lambda+\frac{5}{2}}]\!],\,\,[\![\Upsilon_{\lambda+6,\lambda+\frac{17}{2}},\,
\frak{J}_{7}^{-1,\lambda}]\!]\right),\\[2pt]
7)\,&\left(\delta(\frak{J}_{10}^{-1,\lambda}),\,\,[\![\frak{J}_{\frac{11}{2}}^{-1,\lambda+\frac{9}{2}},\,
\frak{J}_{\frac{11}{2}}^{-1,\lambda}]\!]\right).
\end{array}
\end{equation*}
\end{lemma}
\begin{proof} (Proposition \ref{th4})
The fourth order integrability conditions of the infinitesimal
deformation~(\ref{infdef1}) follow from Lemma \ref{benfraj3} and
Lemma \ref{benfraj4} together with Proposition \ref{th2} and
Proposition \ref{pr3} and arguments similar to those from the
proof of proposition \ref{pr3}. Under these conditions, the
fourth-order term ${\frak L}^{(4)}$ can be given by:
\begin{equation*}
{\frak L}^{(4)
}=-\epsilon_{17,\lambda}\,t_{\lambda+6,\lambda+8}t_{\lambda+\frac{7}{2},\lambda+6}
\left(t_{\lambda+\frac{3}{2},\lambda+\frac{7}{2}}
t_{\lambda,\lambda+\frac{3}{2}}-
t_{\lambda+2,\lambda+\frac{7}{2}}t_{\lambda,\lambda+2}\right)\frak{J}_{9}^{-1,\lambda}.
\end{equation*}
\end{proof}
\begin{proposition}\label{th5} The 5th order integrability conditions of the
infinitesimal deformation~(\ref{infdef1}) are the following:
\begin{itemize}
 \item [a)] For $2(\beta-\lambda)\in\left\{19,\,
\dots,\,2n\right\}:$
\begin{gather*}
t_{\lambda+8,\lambda+\frac{19}{2}}
t_{\lambda+6,\lambda+8}t_{\lambda+\frac{7}{2},\lambda+6}
\left(t_{\lambda+\frac{3}{2},\lambda+\frac{7}{2}}
t_{\lambda,\lambda+\frac{3}{2}}-
t_{\lambda+2,\lambda+\frac{7}{2}}t_{\lambda,\lambda+2}\right) =0,
\end{gather*}
\item [b)] For $2(\beta-\lambda)\in\left\{20,\,
\dots,\,2n\right\}:$
\begin{gather*}
t_{\lambda+8,\lambda+10}
t_{\lambda+6,\lambda+8}t_{\lambda+\frac{7}{2},\lambda+6}
\left(t_{\lambda+\frac{3}{2},\lambda+\frac{7}{2}}
t_{\lambda,\lambda+\frac{3}{2}}-
t_{\lambda+2,\lambda+\frac{7}{2}}t_{\lambda,\lambda+2}\right)=0,\\
t_{\lambda,\lambda+2}t_{\lambda+8,\lambda+10}
t_{\lambda+\frac{11}{2},\lambda+8}
\left(t_{\lambda+\frac{7}{2},\lambda+\frac{11}{2}}
t_{\lambda+2,\lambda+\frac{7}{2}}-
t_{\lambda+4,\lambda+\frac{11}{2}}t_{\lambda+2,\lambda+4}\right)=0,
\end{gather*}
\item [c)] For $2(\beta-\lambda)\in\left\{21,\,
\dots,\,2n\right\}:$
\begin{gather*}
t_{\lambda+8,\lambda+\frac{21}{2}}
t_{\lambda+6,\lambda+8}t_{\lambda+\frac{7}{2},\lambda+6}
\left(t_{\lambda+\frac{3}{2},\lambda+\frac{7}{2}}
t_{\lambda,\lambda+\frac{3}{2}}-
t_{\lambda+2,\lambda+\frac{7}{2}}t_{\lambda,\lambda+2}\right)=0,\\
t_{\lambda+8,\lambda+\frac{21}{2}}t_{\lambda+\frac{13}{2},\lambda+8}
t_{\lambda+\frac{9}{2},\lambda+\frac{13}{2}}
\left(t_{\lambda+\frac{5}{2},\lambda+\frac{9}{2}}
t_{\lambda,\lambda+\frac{5}{2}}-
t_{\lambda+2,\lambda+\frac{9}{2}}t_{\lambda,\lambda+2}\right)=0,\\
t_{\lambda,\lambda+\frac{5}{2}}t_{\lambda+\frac{17}{2},\lambda+\frac{21}{2}}
t_{\lambda+6,\lambda+\frac{17}{2}}
t_{\lambda+\frac{9}{2},\lambda+6}t_{\lambda+\frac{5}{2},\lambda+\frac{9}{2}}=0.
\end{gather*}
\end{itemize}
\end{proposition}
To prove Proposition \ref{th5}, we need the following lemma which
we can check by a direct computation.
\begin{lemma}\label{benfraj5}
Each of the following systems is linearly independent
\begin{equation*}
\begin{array}{lll}
1)\,&\left(\delta(\frak{J}_{\frac{21}{2}}^{-1,\lambda}),\,\,
\alpha_{\lambda+6}^{-1}\epsilon_{9,\lambda}
[\![\frak{J}_{\frac{9}{2}}^{-1,\lambda+6},\,
\frak{J}_{7}^{-1,\lambda}]\!]-\epsilon_{17,\lambda}
[\![\Upsilon_{\lambda+8,\lambda+\frac{19}{2}}
,\,\frak{J}_{9}^{-1,\lambda}]\!]\right),\\[2pt]
2)\,&\left(\delta(\frak{J}_{11}^{-1,\lambda}),\,\,
[\![\frak{J}_{9}^{-1,\lambda+2},\,
\Upsilon_{\lambda,\lambda+2}]\!],\,\, \epsilon_{17,\lambda}
[\![\Upsilon_{\lambda+8,\lambda+10},\,
\frak{J}_{9}^{-1,\lambda}]\!]+\frac{1}{3}\beta^{-1}_{\lambda+6}\epsilon_{9,\lambda}
[\![\frak{J}_{5}^{-1,\lambda+6},\,
\frak{J}_{7}^{-1,\lambda}]\!]\right),\\[2pt]
3)\,&\Big(\delta(\frak{J}_{\frac{23}{2}}^{-1,\lambda}),\,\,
[\![\frak{J}_{7}^{-1,\lambda+\frac{9}{2}},\,
\frak{J}_{\frac{11}{2}}^{-1,\lambda}]\!],\,\,
[\![\frak{J}_{9}^{-1,\lambda+\frac{5}{2}},\,
\Upsilon_{\lambda,\lambda+\frac{5}{2}}]\!],\,\,\\
&\quad~~~~ \epsilon_{9,\lambda}\gamma_{\lambda+6}^{-1}
[\![\frak{J}_{\frac{11}{2}}^{-1,\lambda+6},\,
\frak{J}_{7}^{-1,\lambda}]\!]-
\epsilon_{17,\lambda}[\![\Upsilon_{\lambda+8,\lambda+\frac{21}{2}},\,
\frak{J}_{9}^{-1,\lambda}]\!] \Big).
\end{array}
\end{equation*}
\end{lemma}
\begin{proof}
(Proposition \ref{th5}) Using the same arguments as in proof of
proposition \ref{pr3} together with Lemma \ref{benfraj5},
Proposition \ref{pr3} and Proposition \ref{th4}, we get the
necessary integrability
conditions for ${\frak L}^{(5)}$.
Under these conditions, it can be easily checked that
$\delta({\frak L}^{(m)})=0$ for $m=5,6,7,8.$
\end{proof}
\vskip3mm The main result in this section is the following
theorem.
\begin{theorem}\label{threc}
The conditions given in Propositions \ref{th2}, \ref{pr3},
\ref{th4}, \ref{th5} are necessary and sufficient for the
integrability of the infinitesimal deformation~(\ref{infdef1}).
Moreover, any formal $\mathfrak{osp}(1|2)$-trivial deformation of
the $\mathcal{K}(1)$-module ${\frak S}^n_{\beta}$  is equivalent
to a polynomial one of degree $\leq4$.
\end{theorem}
\begin{proof} Of course these conditions are necessary.
Now, we show that these conditions are sufficient. The solution
$\frak{L}^{(m)}$ of the Maurer-Cartan equation is defined up to a
1-cocycle and it has been shown in~\cite{fi,aalo} that different
choices of solutions of the Maurer-Cartan equation correspond to
equivalent deformations. Thus, we can always reduce
$\frak{L}^{(m)},$ for $m=5,6,7,8,$ to zero by equivalence. Then,
by recurrence, the terms $\frak{L}^{(m)}$, for $m\geq9$, satisfy
the equation $\delta(\frak{L}^{(m)})=0$ and can also be reduced to
the identically zero map.
\end{proof}

{\remark {\rm There are no integrability conditions of any
infinitesimal $\mathfrak{osp}(1|2)$-trivial deformation of the
$\mathcal{K}(1)$-module ${\frak S}^n_{\beta}$ if $n<5$. In this
case, any formal $\mathfrak{osp}(1|2)$-trivial deformation is
equivalent to its infinitesimal part.  }}
\section{Examples}
We study formal $\mathfrak{osp}(1|2)$-trivial deformations of
$\mathcal{K}(1)$-modules ${\frak S}^{n}_{\lambda+n}$ for some
$n\in{1\over2}\mathbb{N}$ and for arbitrary generic
$\lambda\in\mathbb{K}.$ For $n<5$, each of these deformations is
equivalent to its infinitesimal one, without any integrability
condition.

{\example \label{Example1}{\rm The $\mathcal{K}(1)$-module ${\frak
S}^{5}_{\lambda+5}$.}}
\begin{proposition} The $\mathcal{K}(1)$-module
${\frak S}^{5}_{\lambda+5}$ admits six formal
$\mathfrak{osp}(1|2)$-trivial deformations with 18 independent
parameters. These deformations are polynomial of degree 3.
\end{proposition}
\begin{proof} In this case, any
$\mathfrak{osp}(1|2)$-trivial deformation is given by
\begin{equation}
\widetilde{\frak L}_{X_F}=\frak{L}_{X_F}+{\frak
L}^{(1)}_{X_F}+{\frak L}^{(2)}_{X_F}+{\frak L}^{(3)}_{X_F},
\end{equation}
where
\begin{equation*}\begin{array}{ll}
{\frak L}^{(1)}=&t_{\lambda, \lambda+\frac{3}{2}}\,
\Upsilon_{\lambda,\lambda+\frac{3}{2}}+t_{\lambda, \lambda+2}\,
\Upsilon_{\lambda,\lambda+2}+t_{\lambda, \lambda+\frac{5}{2}}\,
\Upsilon_{\lambda,\lambda+\frac{5}{2}}+t_{\lambda+{1\over2},
\lambda+2}\,
\Upsilon_{\lambda+{1\over2},\lambda+2}\\&+t_{\lambda+{1\over2},
\lambda+\frac{5}{2}}\,\Upsilon_{\lambda+{1\over2},\lambda+\frac{5}{2}}+t_{\lambda+{1\over2},
\lambda+3}\,\Upsilon_{\lambda+{1\over2},\lambda+3}+t_{\lambda+1,
\lambda+\frac{5}{2}}\,
\Upsilon_{\lambda+1,\lambda+\frac{5}{2}}+t_{\lambda+1,
\lambda+3}\, \Upsilon_{\lambda+1,\lambda+3}\\&+t_{\lambda+1,
\lambda+\frac{7}{2}}\,
\Upsilon_{\lambda+1,\lambda+\frac{7}{2}}+t_{\lambda+\frac{3}{2},
\lambda+3}\,
\Upsilon_{\lambda+\frac{3}{2},\lambda+3}+t_{\lambda+\frac{3}{2},
\lambda+\frac{7}{2}}\,
\Upsilon_{\lambda+\frac{3}{2},\lambda+\frac{7}{2}}+t_{\lambda+\frac{3}{2},
\lambda+4}\,
\Upsilon_{\lambda+\frac{3}{2},\lambda+4}\\&+t_{\lambda+2,
\lambda+\frac{7}{2}}\,
\Upsilon_{\lambda+2,\lambda+\frac{7}{2}}+t_{\lambda+2,
\lambda+4}\, \Upsilon_{\lambda+2, \lambda+4}+t_{\lambda+2,
\lambda+\frac{9}{2}}\, \Upsilon_{\lambda+2,
\lambda+\frac{9}{2}}+t_{\lambda+\frac{5}{2}, \lambda+4}\,
\Upsilon_{\lambda+\frac{5}{2},
\lambda+4}\\&+t_{\lambda+\frac{5}{2}, \lambda+\frac{9}{2}}\,
\Upsilon_{\lambda+\frac{5}{2},
\lambda+\frac{9}{2}}+t_{\lambda+\frac{5}{2}, \lambda+5}\,
\Upsilon_{\lambda+\frac{5}{2}, \lambda+5}+t_{\lambda+3,
\lambda+\frac{9}{2}}\, \Upsilon_{\lambda+3,
\lambda+\frac{9}{2}}+t_{\lambda+3, \lambda+5}\,
\Upsilon_{\lambda+3, \lambda+5}\\&+t_{\lambda+\frac{7}{2},
\lambda+5}\, \Upsilon_{\lambda+\frac{7}{2}, \lambda+5},\\[10pt]
{\frak L}^{(2) }=
&-\sum_\mu\zeta_\mu^{-1}t_{\mu+\frac{3}{2},\mu+3}
t_{\mu,\mu+\frac{3}{2}}\frak{J}_4^{-1,\mu}\\[2pt]
&-\sum_\nu\alpha_\nu^{-1}(t_{\nu+\frac{3}{2},\nu+\frac{7}{2}}
t_{\nu,\nu+\frac{3}{2}}- t_{\nu+2,\nu+\frac{7}{2}}t_{\nu,\nu+2})
\frak{J}_\frac{9}{2}^{-1,\nu}\\[2pt]&-\sum_\varepsilon\beta_\varepsilon^{-1}
(t_{\varepsilon+\frac{3}{2},\varepsilon+4}
t_{\varepsilon,\varepsilon+\frac{3}{2}}+
t_{\varepsilon+\frac{5}{2},\varepsilon+4}
t_{\varepsilon,\varepsilon+\frac{5}{2}}+{1\over3}t_{\varepsilon+2,\varepsilon+4}t_{\varepsilon,\varepsilon+2})
\frak{J}_{5}^{-1,\varepsilon}\\[2pt]
&-\sum_\ell\gamma_\ell^{-1}(t_{\ell+\frac{5}{2},\ell+\frac{9}{2}}
t_{\ell,\ell+\frac{5}{2}}-
t_{\ell+2,\ell+\frac{9}{2}}t_{\ell,\ell+2})
\frak{J}_\frac{11}{2}^{-1,\ell},\\[10pt]
{\frak L}^{(3) }=
&\sum_\ell\xi_\ell^{-1}t_{\ell+3,\ell+\frac{9}{2}}
t_{\ell+\frac{3}{2},\ell+3} t_{\ell,\ell+\frac{3}{2}}\,
\frak{J}_{11\over2}^{-1,\ell}\\[2pt]&+\,\epsilon_{3,\lambda}\,
\left(t_{\lambda+3,\lambda+5} t_{\lambda+\frac{3}{2},\lambda+3}-
t_{\lambda+\frac{7}{2},\lambda+5}t_{\lambda+\frac{3}{2},\lambda+\frac{7}{2}}\right)
t_{\lambda,\lambda+\frac{3}{2}}
\frak{J}_{6}^{-1,\lambda}\end{array}
\end{equation*}
with
$\mu\in\{\lambda,\,\lambda+{1\over2},\,\lambda+1,\,\lambda+\frac{3}{2},\lambda+2\}$,
$\nu\in\{\lambda,\,\lambda+{1\over2},\,\lambda+1,\,\lambda+\frac{3}{2}\}$,
$\varepsilon\in\{\lambda,\,\lambda+{1\over2},\,\lambda+1\}$ and
$\ell\in\{\lambda,\,\lambda+{1\over2}\}$ . The following equations
\begin{gather}
\label{co1}t_{\lambda,
\lambda+\frac{5}{2}}\,t_{\lambda+\frac{5}{2}, \lambda+5}=0,\\
t_{\lambda,\lambda+\frac{3}{2}}\left(\epsilon_{1,\lambda}\,
t_{\lambda+3,\lambda+5}
t_{\lambda+\frac{3}{2},\lambda+3}+(1-\epsilon_{1,\lambda})\,
t_{\lambda+\frac{7}{2},\lambda+5}
t_{\lambda+\frac{3}{2},\lambda+\frac{7}{2}}
\right)=0,\\
t_{\lambda,\lambda+\frac{3}{2}}\left(\epsilon_{2,\lambda}\,t_{\lambda+\frac{7}{2},\lambda+5}
t_{\lambda+\frac{3}{2},\lambda+\frac{7}{2}}-(1+\epsilon_{2,\lambda})\,
 t_{\lambda+3,\lambda+5}
t_{\lambda+\frac{3}{2},\lambda+3}\right)=0,\\
t_{\lambda+\frac{7}{2},\lambda+5}t_{\lambda+2,\lambda+\frac{7}{2}}
t_{\lambda,\lambda+2}=0.\label{co2}
\end{gather} are the integrability conditions of the
infinitesimal deformation. The formal deformations with the
greatest number of independent parameters are those corresponding
to
$t_{\lambda,\lambda+\frac{5}{2}}t_{\lambda+\frac{5}{2},\lambda+5}=t_{\lambda+\frac{7}{2},\lambda+5}t_{\lambda+2,\lambda+\frac{7}{2}}
t_{\lambda,\lambda+2}=t_{\lambda,\lambda+\frac{3}{2}}=0$. So, we
must kill at least three parameters and there are six choices.
Thus, there are only six deformations with eighteen independent
parameters. Of course, there are many formal deformations with
less then eighteen independent parameters. The deformation $
\widetilde{\frak L}_{X_F}=\frak{L}_{X_F}+{\frak
L}^{(1)}_{X_F}+{\frak L}^{(2) }_{X_F}+{\frak L}^{(3) }_{X_F}, $ is
the miniversal $\mathfrak{osp}(1|2)$-trivial deformation of
${\frak S}^{5}_{\lambda+5}$ with base $\mathcal{A}=
\mathbb{C}[t]/\mathcal{R}$, where $t=(t_{\lambda,
\lambda+\frac{3}{2}},\dots)$ is the family of all parameters given
in the expression of ${\frak L}^{(1)}$ and $\mathcal{R}$ is the
ideal generated by the left hand sides of
(\ref{co1})--(\ref{co2}).
\end{proof}

{\example \label{Example2} {\rm The $\mathcal{K}(1)$-module
${\frak S}^\frac{11}{2}_{\lambda+\frac{11}{2}}$.}}
\begin{proposition} The $\mathcal{K}(1)$-module
${\frak S}^\frac{11}{2}_{\lambda+\frac{11}{2}}$ admits 36
$\mathfrak{osp}(1|2)$-trivial deformations with 17 independent
parameters. These deformations are polynomial of degree 3.
\end{proposition}
\begin{proof}
Any $\mathfrak{osp}(1|2)$-trivial deformation of ${\frak
S}^\frac{11}{2}_{\lambda+\frac{11}{2}}$ is given by
\begin{equation}
\widetilde{\frak L}_{X_F}=\frak{L}_{X_F}+{\frak
L}^{(1)}_{X_F}+{\frak L}^{(2)}_{X_F}+{\frak L}^{(3)}_{X_F},
\end{equation}
where
\begin{equation*}\begin{array}{lll}
{\frak L}^{(1)}=&t_{\lambda, \lambda+\frac{3}{2}}\,
\Upsilon_{\lambda,\lambda+\frac{3}{2}}+t_{\lambda, \lambda+2}\,
\Upsilon_{\lambda,\lambda+2}+t_{\lambda, \lambda+\frac{5}{2}}\,
\Upsilon_{\lambda,\lambda+\frac{5}{2}}+t_{\lambda+{1\over2},
\lambda+2}\,
\Upsilon_{\lambda+{1\over2},\lambda+2}\\&+t_{\lambda+{1\over2},
\lambda+\frac{5}{2}}\,\Upsilon_{\lambda+{1\over2},\lambda+\frac{5}{2}}+t_{\lambda+{1\over2},
\lambda+3}\,\Upsilon_{\lambda+{1\over2},\lambda+3}+t_{\lambda+1,
\lambda+\frac{5}{2}}\,
\Upsilon_{\lambda+1,\lambda+\frac{5}{2}}+t_{\lambda+1,
\lambda+3}\, \Upsilon_{\lambda+1,\lambda+3}\\&+t_{\lambda+1,
\lambda+\frac{7}{2}}\,
\Upsilon_{\lambda+1,\lambda+\frac{7}{2}}+t_{\lambda+\frac{3}{2},
\lambda+3}\,
\Upsilon_{\lambda+\frac{3}{2},\lambda+3}+t_{\lambda+\frac{3}{2},
\lambda+\frac{7}{2}}\,
\Upsilon_{\lambda+\frac{3}{2},\lambda+\frac{7}{2}}+t_{\lambda+\frac{3}{2},
\lambda+4}\,
\Upsilon_{\lambda+\frac{3}{2},\lambda+4}\\&+t_{\lambda+2,
\lambda+\frac{7}{2}}\,
\Upsilon_{\lambda+2,\lambda+\frac{7}{2}}+t_{\lambda+2,
\lambda+4}\, \Upsilon_{\lambda+2, \lambda+4}+t_{\lambda+2,
\lambda+\frac{9}{2}}\, \Upsilon_{\lambda+2,
\lambda+\frac{9}{2}}+t_{\lambda+\frac{5}{2}, \lambda+4}\,
\Upsilon_{\lambda+\frac{5}{2},
\lambda+4}\\&+t_{\lambda+\frac{5}{2}, \lambda+\frac{9}{2}}\,
\Upsilon_{\lambda+\frac{5}{2},
\lambda+\frac{9}{2}}+t_{\lambda+\frac{5}{2}, \lambda+5}\,
\Upsilon_{\lambda+\frac{5}{2}, \lambda+5}+t_{\lambda+3,
\lambda+\frac{9}{2}}\, \Upsilon_{\lambda+3,
\lambda+\frac{9}{2}}+t_{\lambda+3, \lambda+5}\,
\Upsilon_{\lambda+3, \lambda+5}\\&+t_{\lambda+3,
\lambda+\frac{11}{2}}\, \Upsilon_{\lambda+3,
\lambda+\frac{11}{2}}+t_{\lambda+{7\over2}, \lambda+5}\,
\Upsilon_{\lambda+{7\over2}, \lambda+5}+t_{\lambda+\frac{7}{2},
\lambda+\frac{11}{2}}\, \Upsilon_{\lambda+\frac{7}{2},
\lambda+\frac{11}{2}}\\&+t_{\lambda+4, \lambda+\frac{11}{2}}\,
\Upsilon_{\lambda+4, \lambda+\frac{11}{2}},\\[10pt]
{\frak L}^{(2) }=
&-\sum_\mu\zeta_\mu^{-1}t_{\mu+\frac{3}{2},\mu+3}
t_{\mu,\mu+\frac{3}{2}}\frak{J}_4^{-1,\mu}\\[2pt]
&-\sum_\nu\alpha_\nu^{-1}(t_{\nu+\frac{3}{2},\nu+\frac{7}{2}}
t_{\nu,\nu+\frac{3}{2}}- t_{\nu+2,\nu+\frac{7}{2}}t_{\nu,\nu+2})
\frak{J}_\frac{9}{2}^{-1,\nu}\\[2pt]&-\sum_\varepsilon\beta_\varepsilon^{-1}
(t_{\varepsilon+\frac{3}{2},\varepsilon+4}
t_{\varepsilon,\varepsilon+\frac{3}{2}}+
t_{\varepsilon+\frac{5}{2},\varepsilon+4}
t_{\varepsilon,\varepsilon+\frac{5}{2}}+{1\over3}t_{\varepsilon+2,\varepsilon+4}t_{\varepsilon,\varepsilon+2})
\frak{J}_{5}^{-1,\varepsilon}\\[2pt]
&-\sum_\ell\gamma_\ell^{-1}(t_{\ell+\frac{5}{2},\ell+\frac{9}{2}}
t_{\ell,\ell+\frac{5}{2}}-
t_{\ell+2,\ell+\frac{9}{2}}t_{\ell,\ell+2})
\frak{J}_\frac{11}{2}^{-1,\ell},\\[10pt]
{\frak L}^{(3) }=
&\sum_\ell\xi_\ell^{-1}t_{\ell+3,\ell+\frac{9}{2}}
t_{\ell+\frac{3}{2},\ell+3} t_{\ell,\ell+\frac{3}{2}}\,
\frak{J}_{11\over2}^{-1,\ell}+\\&\sum_\iota\epsilon_{3,\iota}\,
\left(t_{\iota+3,\iota+5} t_{\iota+\frac{3}{2},\iota+3}-
t_{\iota+\frac{7}{2},\iota+5}t_{\iota+\frac{3}{2},\iota+\frac{7}{2}}\right)
t_{\iota,\iota+\frac{3}{2}}
\frak{J}_{6}^{-1,\iota}\\\end{array}\end{equation*} with
$\mu\in\{\lambda,\,\lambda+{1\over2},\,\lambda+1,\,\lambda+\frac{3}{2},\lambda+2,\,\lambda+\frac{5}{2}\}$,
$\nu\in\{\lambda,\,\lambda+{1\over2},\,\lambda+1,\,\lambda+\frac{3}{2},\,\lambda+2\}$,
$\varepsilon\in\{\lambda,\,\lambda+{1\over2},\,\lambda+1,\,\lambda+\frac{3}{2}\}$,
$\ell\in\{\lambda,\,\lambda+{1\over2},\,\lambda+1\}$ and
$\iota\in\{\lambda,\,\lambda+{1\over2}\}$. The integrability
conditions of this infinitesimal deformation vanishing of the
following polynomials, where in the first four lines
$\mu\in\left\{\lambda,\, \lambda+{1\over2}\right\}:$
\begin{equation}\label{2-cocy}\begin{array}{llllllll}
t_{\mu,\mu+\frac{5}{2}}\,t_{\mu+\frac{5}{2}, \mu+5},\\
t_{\mu+\frac{7}{2},\mu+5}t_{\mu+2,\mu+\frac{7}{2}}
t_{\mu,\mu+2},\\
t_{\mu,\mu+\frac{3}{2}}\left(\epsilon_{1,\mu}\, t_{\mu+3,\mu+5}
t_{\mu+\frac{3}{2},\mu+3}+(1-\epsilon_{1,\mu})\,
t_{\mu+\frac{7}{2},\mu+5} t_{\mu+\frac{3}{2},\mu+\frac{7}{2}}
\right),\\
t_{\mu,\mu+\frac{3}{2}}\left(\epsilon_{2,\mu}\,t_{\mu+\frac{7}{2},\mu+5}
t_{\mu+\frac{3}{2},\mu+\frac{7}{2}}-(1+\epsilon_{2,\mu})\,
 t_{\mu+3,\mu+5}
t_{\mu+\frac{3}{2},\mu+3}\right),\\
t_{\lambda,\lambda+\frac{3}{2}} t_{\lambda+3,\lambda+\frac{11}{2}}
t_{\lambda+\frac{3}{2},\lambda+3},\\
t_{\lambda+4,\lambda+\frac{11}{2}}t_{\lambda+\frac{5}{2},\lambda+4}
t_{\lambda,\lambda+\frac{5}{2}},\\
t_{\lambda+4,\lambda+\frac{11}{2}}\left(3(1+\epsilon_{4,\lambda})\,t_{\lambda+\frac{3}{2},\lambda+4}
t_{\lambda,\lambda+\frac{3}{2}}+
t_{\lambda+2,\lambda+4}t_{\lambda,\lambda+2}\right)+
\epsilon_{4,\lambda}\, t_{\lambda,\lambda+\frac{3}{2}}
t_{\lambda+\frac{7}{2},\lambda+\frac{11}{2}}t_{\lambda+\frac{3}{2},\lambda+\frac{7}{2}},\\
t_{\lambda+\frac{7}{2},\lambda+\frac{11}{2}}
\left((1+\frac{\epsilon_{5,\lambda}}{3})\,t_{\lambda+\frac{3}{2},\lambda+\frac{7}{2}}
t_{\lambda,\lambda+\frac{3}{2}}-
t_{\lambda+2,\lambda+\frac{7}{2}}t_{\lambda,\lambda+2}\right)+\epsilon_{5,\lambda}\,t_{\lambda+4,\lambda+\frac{11}{2}}t_{\lambda+\frac{3}{2},\lambda+4}
t_{\lambda,\lambda+\frac{3}{2}},\\
\epsilon_{6,\lambda}\,t_{\lambda,\lambda+\frac{3}{2}}\left(
t_{\lambda+4,\lambda+\frac{11}{2}}t_{\lambda+\frac{3}{2},\lambda+4}+\frac{1}{3}\,
t_{\lambda+\frac{7}{2},\lambda+\frac{11}{2}}t_{\lambda+\frac{3}{2},\lambda+\frac{7}{2}}\right)+\\
\quad
t_{\lambda,\lambda+2}\left(t_{\lambda+\frac{7}{2},\lambda+\frac{11}{2}}
t_{\lambda+2,\lambda+\frac{7}{2}}-
t_{\lambda+4,\lambda+\frac{11}{2}}t_{\lambda+2,\lambda+4}\right).
\end{array}
\end{equation}
These deformations are with 24 parameters $t_{\mu,\nu}$ which are
subject to conditions (\ref{2-cocy}). Obviously, we can construct
many $\mathfrak{osp}(1|2)$-trivial deformation of
$\mathfrak{S}^\frac{11}{2}_{\lambda+\frac{11}{2}}$ with
independent parameters. But, to have the greatest number of
independent parameters, we see that we must kill at least seven
parameters, that is , we put
\begin{equation*}
t_{\lambda,\lambda+\frac{3}{2}}=t_{\lambda+{1\over2},\lambda+2}=t_{\lambda,\lambda+2}=0
\quad\text{ and }\quad
t_{\mu,\mu+\frac{5}{2}}\,t_{\mu+\frac{5}{2},
\mu+5}=t_{\mu+\frac{7}{2},\mu+5}t_{\mu+2,\mu+\frac{7}{2}}
t_{\mu,\mu+2}=0
\end{equation*}
where $\mu=\lambda$ or $\lambda+{1\over2}$. So, there are 36
possible choices of such parameters.
\end{proof}

\section*{Acknowledgments}
We would like to thank Dimitry Leites and
Valentin Ovsienko for helpful discussions. We are also grateful to
the referee for his comments and suggestion.

\end{document}